\documentclass[a4paper,10pt]{article}
\usepackage{amsmath,amssymb,amsthm}
\usepackage{picins} 
\theoremstyle{definition}
 \newtheorem{dfn}{Definition}[section]
 \newtheorem{remark}[dfn]{Remark}  
\theoremstyle{plain}
 \newtheorem{thm}[dfn]{Theorem}
 
 \newtheorem{lem}[dfn]{Lemma}
 \newtheorem{cor}[dfn]{Corollary}

\numberwithin{equation}{section}
\newcommand{\DV}{{\rm Div}\,}
\newcommand{\dv}{{\rm div}\,}

\newcommand{\BR}{{\mathbb R}}
\newcommand{\BC}{{\mathbb C}} 

\newcommand{\BN}{{\mathbb N}}

\newcommand{\BF}{{\mathbb F}}

\newcommand{\CA}{{\mathcal A}}

\newcommand{\CE}{{\mathcal E}}
\newcommand{\CF}{{\mathcal F}}
\newcommand{\CI}{{\mathcal I}}

\newcommand{\CL}{{\mathcal L}}
\newcommand{\CM}{{\mathcal M}}

\newcommand{\CP}{{\mathcal P}}
\newcommand{\CR}{{\mathcal R}}
\newcommand{\CS}{{\mathcal S}}
\newcommand{\CT}{{\mathcal T}}
\newcommand{\CH}{{\mathcal H}}

\newcommand{\CU}{{\mathcal U}}

\newcommand{\CX}{{\mathcal X}}
\newcommand{\CY}{{\mathcal Y}}
\newcommand{\CZ}{{\mathcal Z}}
\newcommand{\bD}{{\mathbf D}}
\newcommand{\bF}{{\mathbf F}}
\newcommand{\bI}{{\mathbf I}}
\newcommand{\bK}{{\mathbf K}}

\newcommand{\bM}{{\mathbf M}}
\newcommand{\bG}{{\mathbf G}}

\newcommand{\bS}{{\mathbf S}}
\newcommand{\bT}{{\mathbf T}}
\newcommand{\bH}{{\mathbf H}}

\newcommand{\bh}{{\mathbf h}}
\newcommand{\ba}{{\mathbf a}}
\newcommand{\bb}{{\mathbf b}}

\newcommand{\bn}{{\mathbf n}}

\newcommand{\bw}{{\mathbf w}}

\newcommand{\bu}{{\mathbf u}}
\newcommand{\bv}{{\mathbf v}}
\newcommand{\bg}{{\mathbf g}}
\newcommand{\bff}{{\mathbf f}}

\newcommand{\pd}{\partial}
\newcommand{\HSp}{\BR^N_+}
\newcommand{\HSn}{\BR^N_-}

\newcommand{\Hol}{{\rm Hol}}
\newcommand{\re}{{\rm Re}\,}
\newcommand{\im}{{\rm Im}\,}
\newcommand{\bng}{{\bn_{\Gamma(t)}}}

\newenvironment{cases*}%
{%
\left\{
\begin{array}{@{}r@{\;}l@{\quad}r@{}}
}%
{\end{array}\right.}
\begin{document}
\title{On the $\CR$-boundedness for the two phase problem \\
with phase transition: \\
compressible-incompressible model problem}
\author{
%Katharina Schade\thanks{Department of Mathematics, 
%TU Darmstadt, \endgraf
%e-mail address: schade@mathematik.tu-darmstadt.de \endgraf
%Partially supported by DFG IRTG 1529 at TU Darmstadt}
%\enskip and \enskip 
Yoshihiro Shibata
\thanks{Department of Mathematics and Research Institute of 
Science and Engineering \endgraf
Waseda University, 
Ohkubo 3-4-1, Shinjuku-ku, Tokyo 169-8555, Japan. \endgraf
e-mail address: yshibata@waseda.jp \endgraf
Partially supported by JST CREST, JSPS
Grant-in-aid for Scientific Research (S) \# 24224004.
%and JSPS Japanese-Germann Graduate externship at Waseda University. }
}}
\date{}
\maketitle
%%%%%%%%%%%%%%% 
\begin{abstract}
In this paper, we prove the maximal $L_p$-$L_q$ regularity of the compressible
and incompressible two phase flow with phase transition in the model problem 
case with the help of $\CR$-bounded solution operators corresponding 
to generalized resolvent problem. 
The problem arises from the mathematical study of the 
motion of two-phase flows having gaseous phase and liquid phase separated 
by a sharp interface with phase transition. Using the result obtained in
this paper, in \cite{S0} we proved 
the local well-posedness of free boundary problem for the compressible and 
incompressible two phase flow separated by sharp interface
 with phase transition.
\end{abstract}
\vskip1pc\noindent
{\small \textbf{Mathematics Subject Classification (2012).
Primary: 35R35; Secondary: 35Q30, 76T10.} 
}

\vskip0.5pc\noindent
{\small \textbf{Keywords.}two phase flow, compressible and incompressible viscous flow, surface tension, phase transitions, maximal $L_p$-$L_q$ regularity,
$\CR$-bounded solution operator } 
%%%%%%%%%%%%%%%%%%%%%%%
\section{Introduction} \label{sec:1}
In this paper, we prove the maximal $L_p$-$L_q$ regularity of the compressible
and incompressible two phase flow with phase transition in the model problem 
case with the help of $\CR$-bounded solution operators 
corresponding to generalized
resolvent problem.  The problem arises from the mathematical study of the 
motion of two-phase flows having gaseous phase and liquid phase separated 
by a sharp interface with phase transition, which is formulated as follows: 
Let $\BR^N$ be the $N$ dimensional Euclidean space. 
 Let $\Omega_-$ be a domain in $\BR^N$ 
with boundary $\Gamma$, which is occupied by the liquid.  
Set $\Omega_+ = 
\BR^N-\overline{\Omega_-}$,  where $\overline{\Omega_-}$ stands for 
the closure of domain $\Omega_-$, which is occupied by the gas.  
Let  $\varphi=\varphi(\xi, t)
= (\varphi_1(\xi, t), \ldots, \varphi_N(\xi, t))$ be a function 
defined on $\BR^N$ for each $t \in (0, T)$,
$\xi = (\xi_1, \ldots, \xi_N)$ being the reference coordinate system.  
We assume that the correspondence: $\xi \to \varphi(\xi, t)$ is one to
one map from $\BR^N$ onto itself  for each $t \in (0, T)$. 
Set $(\pd_t\varphi)(\xi, t) = \bv(x, t)$ with $x = \varphi(\xi, t)$, 
\begin{align*}
\Omega_\pm(t) &= \{x = \varphi(\xi, t) \mid \xi \in \Omega_\pm\}, \quad
\Gamma(t)  = \{x = \varphi(\xi, t) \mid \xi \in \Gamma\}. 
\end{align*}
Let $\bng$ be the unit outer normal to $\Gamma(t)$ pointed
from $\Omega_-(t)$ to $\Omega_+(t)$. 
For any $x_0 \in \Gamma(t)$, we set 
$$[[v]](x_0) =\lim_{x\to x_0\atop x\in \Omega_-(t)}v(x)-
\lim_{x\to x_0\atop x\in \Omega_+(t)}v(x)
 \quad\text{(the jump of $v$ accross $\Gamma(t)$)}$$
for any $v$ defined on $\dot\Omega(t)=\Omega_+(t) \cup\Omega_-(t)$.  
In the sequel,  
we write $v_\pm = v|_{\Omega_\pm(t)}$. Moreover, given $v_\pm$ 
defined on $\Omega_\pm(t)$, we define $v$ by
$v(x) = v_\pm(x)$ for $x \in \Omega_\pm(t)$. 

Let $\bu: \dot\Omega(t) \to \BR^N$ be the velocity fields, 
$\rho: \dot\Omega(t) \to (0, \infty)$ the mass field,
$\pi: \dot\Omega(t) \to \BR$ the pressure field ,
$\bT:\dot\Omega(t) \to \BR^N$ 
the stress tensor field, $\theta: \dot\Omega(t)
\to(0, \infty)$ the thermal fields,
$\psi: \dot\Omega(t)\to \BR$ the free energy, 
$\eta: \dot\Omega(t) \to \BR$ the entropy,
$H_\Gamma$  the mean curvature of $\Gamma(t)$,
$\kappa$ the specific heat, $d$ the thermal conductivity
and $\j$ the phase flux.

We assume that $\rho_+\not=\rho_-$. Then, 
the motion of two-phase flows having gaseous phase and liquid phase 
separated 
by a sharp interface with phase transition
 is described as follows: 
\begin{equation}\label{eq1}\begin{split}
\begin{cases}
\rho_+(\pd_t\bu_++ \bu_+\cdot\nabla\bu_+) - \DV\bT_+  = 0, \quad 
\pd_t\rho_+ + \dv(\rho_+\bu_+)  = 0, \\[0.5pc]
\rho_+ \kappa_+
(\pd_t\theta_+ + \bu_+\cdot\nabla\theta_+) 
-\dv(d_+\nabla\theta_+)- \bT_+:\nabla\bu_+
-\dfrac{\pi_+}{\rho}\dv\bu_+ = 0
\end{cases}
&\quad\text{for $x \in \Omega_+(t)$, $t>0$}, \\
\begin{cases}
\rho_{*-}(\pd_t\bu_-+ \bu_-\cdot\nabla\bu_-) - \DV\bT_-  = 0, \quad 
\dv\bu_-  = 0, \\
\rho_{*-} \kappa_-
(\pd_t\theta_- + \bu_-\cdot\nabla\theta_-) 
-\dv(d_-\nabla\theta_-)- \bT_-:\nabla\bu_- = 0
\end{cases}
&\quad\text{for $x \in \Omega_-(t)$, $t>0$}
\end{split}\end{equation}
subject to the interface conditions: for $x \in \Gamma(t)$ and $t > 0$,
\begin{equation}\label{jump:1}
\begin{cases}
[[\dfrac1\rho]]\j^2\bn_\Gamma - [[\bT\bng]] 
= -\sigma H_\Gamma\bng, 
&\quad [[\bu-(\bu\cdot\bng)\bng]]  =0, \\[0.5pc]
\j[[\theta\eta]] - [[d(\nabla\theta)\cdot\bng]]  = 0, &\quad
[[\theta]] = 0, \\[0.5pc]
[[\psi]] + [[\dfrac{1}{2\rho^2}]]\j^2 - [[\frac{1}{\rho}
\bng\cdot\bT\bng]]  = 0, &\quad
\bv\cdot\bng = \dfrac{[[\rho\bu]]\cdot\bng}{[[\rho]]},\\[0.8pc]
\quad
\j=\dfrac{[[\rho\bu]]\cdot\bng}{[[\rho]]},
\end{cases}
\end{equation}
and the initial conditions:
\begin{equation}\label{initial}\begin{split}
&(\rho_+, \bu_+, \theta_+)|_{t=0} = (\rho_{*+}+\rho_{0+}, \bu_{0+},
\theta_{*} + \theta_{0+}) \enskip\text{in $\Omega_+$}, \\
&(\bu_-, \theta_-)|_{t=0} = (\bu_{0+},
\theta_{*} + \theta_{0+}) \enskip\text{in $\Omega_-$}, 
\quad 
h|_{t=0} = h_0 \enskip\text{on $\Gamma$}.
\end{split}\end{equation}
Here,  
$\rho_{*\pm}$, $\theta_{*}$ and $\sigma$ 
are positive constants describing the 
reference mass densities of $\Omega_\pm$, the reference temperature
of both $\Omega_\pm$ and the coefficient of 
the surface tension, respectively. Moreover, 
 $\bT_\pm = \bS_\pm - \pi_\pm\bI$ with
\begin{align*}
&\bS_+ =\bS_+(\bu_+, \rho_+, \theta_+))
= \mu_+\bD(\bu_+) + (\nu_+-\mu_+)\dv\bu\bI, \\
&\bS_-=\bS_-(\bu_-, \theta_-) = \mu_-\bD(\bu_-);
\end{align*}
for any scalor field $\theta$ we set 
$\nabla\theta = (\pd_1\theta, \ldots, \pd_N\theta)$, 
where $\pd_j = \pd/\pd x_j$; 
for any vector field $\bu=(u_1, \ldots, u_N)$
$\nabla\bu$ is the $N\times N$ matrix whose $(i,j)$
component is $\pd_iu_j$,   
$\bD(\bu)$ the deformation tensor whose
$(j,k)$ components are $D_{jk}(\bu) = \frac12(\pd_ju_k
+ \pd_ku_j)$ and $\dv\bu = \sum_{j=1}^N\pd_ju_j$;   
and $\bI$ is the $N\times N$ identity
matrix. Finally, for any matrix field $\bK$ with components 
$K_{ij}$, the quantity $\DV \bK$ is 
an $N$-vector with $i$-component
$\sum_{j=1}^N\pd_jK_{ij}$ and $\bK:\nabla\bu = \sum_{i,j=1}^N
K_{ij}\pd_iu_j$.

We assume that $d$, $\mu$, $\nu_+$, $\kappa$, $\psi$ and 
$\eta$ are 
given as follows: $d_+ =d_+(\rho, \theta)$, 
$\mu_+ = \mu_+(\rho, \theta)$, 
$\nu_+ = \nu_+(\rho, \theta)$, 
$\kappa_+ = \kappa_+(\rho, \theta)$
 are positive 
$C^\infty$ functions with respect to $(\rho, \theta) 
\in (0, \infty)\times(0, \infty)$, and 
$\psi_+=\psi_+(\theta, \rho)$ and $\eta_+=\eta_+(\theta, \rho)$ 
are  real valued $C^\infty$
functions with respect to $(\rho, \theta) \in (0, \infty)\times(0, \infty)$,
while  $d_- =d_-(\theta)$, $\mu_- = \mu_-(\theta)$, 
$\kappa_- = \kappa_-(\theta)$
 are positive 
$C^\infty$ functions with respect to $\theta 
\in (0, \infty)$, and 
$\psi_-=\psi_-(\theta)$ and $\eta_-=\eta_-(\theta)$ are real valued $C^\infty$
functions with respect to $\theta \in (0, \infty)$.
And also,  
we assume that $\pi_+$ is given by $\pi_+ = P_+(\rho, \theta)$,
where $P_+$ is a $C^\infty$ function
with respect to $(\rho, \theta) \in (0, \infty)\times(0, \infty)$
such that $\dfrac{\pd P_+}{\pd\rho} > 0$ for any 
$(\rho, \theta) \in (0, \infty)\times(0, \infty)$. 

Since we prove the local well-posedness of problem 
\eqref{eq1}, \eqref{jump:1} and \eqref{initial} with 
the help of the maximal $L_p$-$L_q$ regularity results 
for the linearized equations, representing $\rho_+$ 
by the integration of the equation:
$\pd_t\rho_+ + \dv(\rho_+\bu_+) = 0$ along the 
characteristic curve generated by $\bu_+$ we eliminate
$\rho_+$ from the equations in $\Omega_+(t)$ and $\Gamma(t)$
\footnote{The idea follows from Tani \cite{T1}.},
so that we have the nonlinear parabolic equations. After this
procedure, as the linearized problem we have the following 
two problems as model problems: In the sequel, 
for any $x_0 \in \BR^N_0$ we define
$f|_\pm(x_0)$ by 
$$f|_\pm(x_0) = \lim_{x\to x_0 \atop x\in\BR^N_\pm}
f(x),$$
where we have set 
$$\BR^N_\pm = \{x = (x_1, \ldots, x_N) \in \BR^N \mid \pm x_N > 0\},
\quad \BR^N_0 = \{x \in \BR^N \mid x_N=0\}. $$
One is  the interface problem for the Stokes system: 
\begin{alignat}2
&\rho_{*+}\pd_t\bu_+ - \DV\bS_{*+}(\bu_+) = \bff_+
&\quad&\text{in $\BR^N_+\times(0, T)$}\nonumber \\
&\rho_{*-}\pd_t\bu_- - \DV\bS_{*-}(\bu_-)+\nabla\pi_- = \bff_-,
\quad\dv \bu_- = f_\dv = \dv\bff_\dv
&\quad&\text{in $\BR^N_-\times(0, T)$}\label{stokes:1} 
\end{alignat}
subject to the interface condition: for $x \in \BR^N_0$ and $t \in (0, T)$
\begin{equation}\label{stokes:2}\begin{split}
&\mu_{*-}D_{iN}(\bu_-)|_- - \mu_{*+}D_{iN}(\bu_+)|_+ = g_i
\quad(i = 1, \ldots, N-1), \\
&(\mu_{*-}D_{NN}(\bu_-) - \pi_-)|_- - (\mu_{*+}D_{NN}(\bu_+)
+ (\nu_{*+} - \mu_{*+})\dv\bu_+)|_+-\sigma\Delta'H = g_N, \\
&\frac{1}{\rho_{*-}}(\mu_{*-}D_{NN}(\bu_-) - \pi_-)|_-
-\frac{1}{\rho_{*+}}(\mu_{*+}D_{NN}(\bu_+)
+ (\nu_{*+} - \mu_{*+})\dv\bu_+)|_+ = g_{N+1}, \\
&u_{-i}|_- - u_{+i}|_+ = h_i \quad(i = 1, \ldots, N-1), \\
&\pd_tH - \Bigl(\frac{\rho_{*-}}{\rho_{*-}-\rho_{*+}}u_{-N}
-\frac{\rho_{*+}}{\rho_{*-} - \rho_{*+}}u_{+N}
\Bigr) = d,
\end{split}\end{equation}
and the initial condition:
\begin{equation}\label{stokes:3}
\bu_\pm|_{t=0} = \bu_{0\pm} \quad\text{in $\BR^N_\pm$},\quad
H|_{t=0} = H_0 \quad\text{in $\BR^N$}, 
\end{equation}
where we have set $\bu_\pm = (u_{\pm1}, \ldots, u_{\pm N})$, 
$\mu_{*+} = \mu(\rho_{*+}, \theta_*)$, 
$\nu_{*+} = \nu_+(\rho_{*+}, \theta_*)$, $\mu_{*-}
= \mu_-(\theta_*)$, $\Delta'H = \sum_{j=1}^{N-1}\pd_j^2 H$,
$S_{*+}(\bu_+) = \mu_{*+}\bD(\bu) + (\nu_{*+}-\mu_{*+})\dv\bu_+\bI$,
and $S_{*-}(\bu_-) = \mu_{*-}\bD(\bu_-)$.

And, another is  the interface problem for the heat equations:
\begin{equation}\label{heat:1}\begin{split}
\rho_{*+}\kappa_{*+}\pd_t\theta_+ - d_{*+}\Delta\theta_+
= \tilde f_+ 
&\quad\text{in $\BR^N_+\times(0, T)$}, \\
\rho_{*-}\kappa_{*-}\pd_t\theta_- - d_{*-}\Delta\theta_-
= \tilde f_- 
&\quad\text{in $\BR^N_-\times(0, T)$}, 
\end{split}\end{equation}
subject to the interface condition: for $x \in \BR^N_0$ and $t \in (0, T)$
\begin{equation}\label{heat:2}
\theta_-|_- - \theta_+|_+=0, \quad
d_{*+}\pd_N\theta_-|_- - d_{*+}\pd_N\theta_+|_+ = \tilde g, 
\end{equation}
and the initial condition:
\begin{equation}\label{heat:3}
\theta_\pm|_{t=0} = \theta_{0\pm}
\quad\text{on $\BR^N_\pm$},
\end{equation}
where we have set $d_{*+} = d(\rho_{*+}, \theta_*)$, $\kappa_{*+} = 
\kappa(\rho_{*+}, \theta_*)$, $d_{*-} = d_-(\theta_*)$ and 
$\kappa_{*-} = \kappa_-(\theta_*)$. 
Note that the interface condition \eqref{stokes:2} is equivalent to the 
following interface condition:
\begin{equation}\label{stokes:2*}\begin{split}
&\mu_{*-}D_{iN}(\bu_-)|_- - \mu_{*+}D_{iN}(\bu_+)|_+ = g_i
\quad(i = 1, \ldots, N-1), \\
&(\mu_{*-}D_{NN}(\bu_-) - \pi_-)|_- =\frac{\rho_{*-}}{\rho_{*-}-\rho_{*+}}
(\sigma\Delta'H + g_N - \rho_{*+}g_{N+1}), \\
&(\mu_{*+}D_{NN}(\bu_+)
+ (\nu_{*+} - \mu_{*+})\dv\bu_+)|_+ = 
\frac{\rho_{*+}}{\rho_{*-}-\rho_{*+}}
(\sigma\Delta'H + g_N - \rho_{*-}g_{N+1}), \\
&u_{-i}|_- - u_{+i}|_+ = h_i \quad(i = 1, \ldots, N-1), \\
&\pd_tH - \Bigl(\frac{\rho_{*-}}{\rho_{*-}-\rho_{*+}}u_{-N}
-\frac{\rho_{*+}}{\rho_{*-} - \rho_{*+}}u_{+N}
\Bigr) = d
\end{split}\end{equation}
The purpose of this paper is to prove 
the following theorem about the maximal 
$L_p$-$L_q$ regularity for problem \eqref{stokes:1},
\eqref{stokes:2}, \eqref{stokes:3}. 
\begin{thm}\label{thm:stokes}
Let $1 < p, q < \infty$ and $0 < T < \infty$. 
Assume that $\rho_{*-} \not= \rho_{*+}$. 
Then, given right-hand sides
of \eqref{stokes:1} and \eqref{stokes:2} 
\begin{align*}
&\bff_\pm \in L_p((0, T), L_q(\BR^N_\pm)^N), \quad 
f_\dv \in L_p((0, T), W^1_q(\BR^N_-)) \cap W^1_p((0, T), 
W^{-1}_q(\BR^N_-))\\ 
&\bff_\dv
\in W^1_p((0, T), L_q(\BR^N_-)^N), \quad 
g_i \in L_p((0, T), W^1_q(\BR^N)) \cap W^1_p((0, T), W^{-1}_q(\BR^N))
\quad(i = 1, \ldots, N+1), \\
&h_j \in L_p((0, T), W^2_q(\BR^N)) \cap W^1_p((0, T), L_p(\BR^N))
\quad(j=1, \ldots, N-1), 
\quad d \in L_p((0, T), W^2_q(\BR^N)),
\end{align*}
and initial data $\bu_{0\pm} \in B^{2(1-1/p)}_{q,p}(\BR^N_\pm)^N$ 
and $H_0 \in B^{3-1/p}_{q,p}(\BR^N)$ 
satisfying the compatibility conditions:
\begin{alignat}2
&\dv \bu_{0-} = f_-|_{t=0}, \quad \bu_{0-} = \bff_\dv|_{t=0} 
&\quad&\text{in $\BR^N_-$}, \nonumber\\
&\mu_{*-}D_{iN}(\bu_{0-})|_- - \mu_{*+}D_{iN}(\bu_{0+})|_+ 
= g_i|_{t=0}
\quad(i = 1, \ldots, N-1) &\quad
&\text{on $\BR^N_0$}, \nonumber\\
&
(\mu_{*+}D_{NN}(\hat\bu_{0+}) + (\nu_{*+}-\mu_{*+})
\dv\hat\bu_{0+})|_+\nonumber\\
&\quad = \frac{\rho_{*+}}{\rho_{*-}-\rho_{*+}}
(\sigma\Delta'H_0 + g_N|_{t=0}- \rho_{*-}g_{N+1}|_{t=0})
&\quad&\text{on $\BR^N_0$}, \nonumber\\
&u_{0-i}|_- - u_{0+i}|_+ = h_i|_{t=0} \quad(i = 1, \ldots, N-1)
&\quad
&\text{on $\BR^N_0$}. \label{comp:1}
\end{alignat}
then, problem \eqref{stokes:1}, \eqref{stokes:2}, \eqref{stokes:3} 
admits unique solutions $\bu_\pm$, $\pi_-$ and $H$ with
\begin{align*}
&\bu_\pm \in L_p((0, T), W^2_q(\BR^N_\pm)^N)\cap
W^1_p((0, T), L_q(\BR^N_\pm)^N), \\
&\pi \in L_p((0, T), \hat W^1_q(\BR^N_-)), \\
&H \in L_p((0, T), W^3_q(\BR^N)) \cap W^1_p((0, T), 
W^2_q(\BR^N))
\end{align*}
possessing the estimates: $\CI_{p,q}(\bu_\pm, \pi_-, H)(t)
\leq Ce^{\gamma t}\BF_{p,q}(t)$ 
for any $t \in (0, T)$ with some positive constants
$C$ and $\gamma$ independent of $t$ and $T$,
where we have set 
\begin{align*}
&\CI_{p,q}(\bu_\pm, \pi_-, H)(t) \\
&= \|\bu_+\|_{L_p((0, t), W^2_q(\BR^N_+))}
+ \|\pd_t\bu_+\|_{L_p((0, t), L_q(\BR^N_+))}
+ \|\bu_-\|_{L_p((0, t), W^2_q(\BR^N_-))}
+ \|\pd_t\bu_-\|_{L_p((0, t), L_q(\BR^N_-))} \\
&+ \|\nabla\pi_-\|_{L_p((0, t), L_q(\BR^N_-))}
+ \|H\|_{L_p((0, t), W^3_q(\BR^N))} 
+ \|\pd_tH\|_{L_p((0, t), W^2_q(\BR^N))}, \\
&\BF_{p,q}(t) \\
&=\{\sum_{\ell=\pm}(\|\bu_{0\ell}\|_{B^{2(1-1/p)}_{q,p}(\BR^N_\ell)}
+ \|\bff_\ell\|_{L_p((0, t), L_q(\BR^N_\ell))})
+ \|f_\dv\|_{L_p((0, t), W^1_q(\BR^N_-))} 
+ \|\pd_tf_\dv\|_{L_p((0, t), W^{-1}_q(\BR^N_-))}\\
&+\|\bff_\dv\|_{L_p((0, T), L_q(\BR^N_-))} 
+ \sum_{i=1}^{N+1}(\|g_i\|_{L_p((0, t), W^1_q(\BR^N))}
+ \|\pd_tg_i\|_{L_p((0, t), W^{-1}_q(\BR^N))})\\
&
+ \sum_{j=1}^{N-1} (\|h_j\|_{L_p((0, t), W^2_q(\BR^N))}
+ \|\pd_th_j\|_{L_p((0, t), L_q(\BR^N))})
+ \|d\|_{L_p((0, t), W^2_q(\BR^N))}\}.
\end{align*}
\end{thm}
%\begin{remark} The space $W^{-1}_q(\BR^N)$ is defined by
%$$W^{-1}_q(\BR^N) = \{f \in 
%L_{1, {\rm loc}}(\BR^N) \cap \CS'(\BR^N) \mid
%\CF^{-1}[(1 + |\xi|^2)^{-1/2}\CF[f](\xi)] \in L_q(\BR^N)\},$$  
%where $\CF$ and $\CF^{-1}$ denote the Fourier transform and 
%the inverse Fourier transform on $\BR^N$. 
% And, for 
%$\Omega= \BR^N_\pm$, $W^{-1}_q(\BR^N_\pm) = \{f \in L_{1, {\rm loc}}
%(\BR^N_\pm) \mid e_\pm f \in W^{-1}_q(\BR^N)\}$, where $e_\pm$ are
%the extension map defined by 
%$e_\pm f(x) = f(x)$ for $x \in \BR^N_\pm$ and 
%$e_\pm f(x) = f(x', -x_N)$ for $x \in \BR^N_\mp$, 
%where $x' = (x_1, \ldots, x_{N-1})$,
%\end{remark}
\begin{remark} The maximal $L_p$-$L_q$ regularity theorem
for problem \eqref{heat:1}, \eqref{heat:2} and 
\eqref{heat:3} seems to be known and employing the 
similar argumentation to that in the proof 
of Theorem \ref{thm:stokes} given in the sequel we can prove
it too, so that we do not consider problem \eqref{heat:1}, \eqref{heat:2} and 
\eqref{heat:3} in this paper. 
\end{remark} 
\begin{remark}
\thetag1 The mathematical study of the compressible and 
incompressible two phase problem is quite rare as far as 
the author knows.  First Denisova \cite{Denisova} studied the evolution 
of the compressible and incompressible two phase flow 
with sharp interface without phase transition under
some restriction on the viscosity coefficients.
Recently, Kubo, Shibata and Soga \cite{KSS} studied the same problem as 
in \cite{Denisova} without surface tension and  phase transition
and proved the maximal $L_p$-$L_q$ regularity 
under the assumption that viscosity coefficients 
are positive constants.  The derivation and the local 
well-posedness of  
problem \eqref{eq1}, \eqref{jump:1} and 
\eqref{initial} are treated in Shibata \cite{S0} and 
all the results of this paper and in \cite{S0}  were 
announced in the abstract of $39^{\rm th}$ Sapporo symposium
on PDE at Hokkaido University (cf. \cite{S3}). 
The incompressible and incompressible two phase problem
with phase transition was studied by J.~Pr\"uss, 
et al. \cite{PSSS, PS, PSW}. 
\end{remark}
\textbf{ Notation}~ Here, we summarize our symbols used throughout the paper. 
$\BN$, $\BR$ and $\BC$ denote the sets of all natural numbers,
real numbers and complex
numbers, respectively. We set $\BN_0 = \BN \cup \{0\}$. 
For any multi-index $\kappa = (\kappa_1, \ldots, \kappa_N) 
\in \BN_0^N$, we write $|\kappa| = \kappa_1 + \cdots + \kappa_N$ and
$\pd_x^\kappa = \pd_1^{\kappa_1}\cdots\pd_N^{\kappa_N}$ with 
$x = (x_1, \ldots, x_N)$ and $\pd_j = \pd/\pd x_j$. 
For any scalar function $f$ and $N$-vector of functions 
$\bg=(g_1, \ldots, g_N)$, we set 
\begin{alignat*}2
\nabla f &= (\pd_1f, \ldots, \pd_Nf), &\enskip
\nabla \bg &= (\pd_ig_j \mid i, j=1, \ldots, N), \\
\nabla^2f & = (\pd^\alpha f \mid |\alpha = 2 ),
&\enskip 
\nabla^2\bg & = (\pd^\alpha g_i \mid |\alpha|=2, i=1, \ldots, N).
\end{alignat*} 
We use  bold small letters to denote $N$-vector or $N$-vector valued
functions and  bold capital letters to denote $N\times N$ matrix or  
$N\times N$ matrix valued functions, respectively.
For any domain $D$ and $1\leq q\leq \infty$, $L_q(D)$, $W^m_q(D)$ 
and $B^\ell_{q,p}(D)$ 
denote the standard Lebesgue space,  
Sobolev space and Besov space, 
while $\|\cdot\|_{L_q(D)}$, $\|\cdot\|_{W^m_q(D)}$
and $\|\cdot\|_{B^\ell_{q,p}(D)}$ 
denote their norms, respectively.  We set $W^0_q(D) = L_q(D)$.
$\hat W^1_q(D)$ is a homogeneous space defined by 
$\hat W^1_q(D) = \{f \in L_{q, {\rm loc}}(D) \mid 
\nabla f \in L_q(D)\}$.  $W^{-1}_q(\BR^N)$ denotes the usual Bessel 
potential space of order $-1$ on $\BR^N$ and $W^{-1}_q(\BR^N_\pm) 
= \{ f \in L_{1, {\rm loc}}(\BR^N_\pm)
 \mid f = g \text{ in $\BR^N_\pm$ for some
$g \in W^{-1}_q(\BR^N)$}\}$. 
For any Banach space $X$ with norm $\|\cdot\|_X$ and 
$1 \leq p \leq \infty$, $L_p((a, b), X)$ and 
$W^m_p((a, b), X)$ denote the usual  Lebesgue space and Sobolev space of
$X$-valued functions defined on an interval $(a, b)$, 
while $\|\cdot\|_{L_p((a, b), X)}$ and 
$\|\cdot\|_{W^m_p((a, b), X)}$ denote their norms, respectively. 
For any $\gamma \in \BR$ we set
$\|e^{\gamma t}f\|_{L_p((a,b),X)} = 
\Bigl(\int^b_a(e^{\gamma t}\|f(t)\|_X)^p\,dt\Bigr)^{1/p}$. 
The $d$-product space of $X$ is defined by  
$X^d = \{f = (f, \ldots, f_d) \mid 
f_i \in X \, (i=1, \ldots, d)\}$, while its norm is denoted by 
$\|\cdot\|_X$ instead of $\|\cdot\|_{X^d}$ for the sake of 
simplicity. For any two Banach spaces $X$ and $Y$,  
$\CL(X, Y)$ denotes the set of all bounded linear operators
from $X$ into $Y$.  $\Hol(U, X)$
denotes the set of all $X$- valued holomorphic functions 
defined on $U$. 
%%%%%%%%%%%%%%%%% 
For $\ba = (a_1, \ldots, a_N)$ and $\bb=(b_1, \ldots, b_N)
\in \BR^N$, we set $\ba\cdot\bb = <\ba, \bb> = 
\sum_{j=1}^Na_jb_j$. For scalar functions $f$, $g$
and $N$-vectors of functions $\bff$, $\bg$, we set 
$(f, g)_D = \int_D f(x)g(x)\,dx$ and $(\bff, \bg)_D
= \int_D\bff(x)\cdot\bg(x)\,dx$. 
The letter $C$ denotes generic constants and 
the constant $C_{a, b, \cdots}$ depends on $a$, $b$, $\cdots$.
The values of constants $C$ and $C_{a, b, \cdots}$ may change from line to 
line. 
%%%%%%%%%%%%%%%%%%

The paper is organized as follows.  In Sect.2, we state the main 
results for the $\CR$ bounded solution operators to the corresponding
resolvent problem to time dependent problem \eqref{stokes:1},
\eqref{stokes:2*} and \eqref{stokes:3}. From Sect.3 through
Sect.5, we consider the problem without surface tension.  In Sect.3, we give 
an exact solution formulas to the resolvent problem. 
In Sect.4, we give some estimates for the multipliers appearing 
in the solution formula. In Sect.5 we analyze the Lopatinski determinant
In Sect.6, we prove the main result for the $\CR$ bounded solution
operators. In Sect.7, using the $\CR$ bounded solution operator,
we prove Theorem \ref{thm:stokes}. 
%%%%%%%%%%%%%%%%%%%
\section{Main results for the $\CR$ bounded solution operators}
In the sequel, for notational simplicity 
$\rho_{*\pm}$, $\mu_{*\pm}$ and $\nu_{*+}$ are denoted by 
$\rho_{\pm}$, $\mu_{\pm}$ and $\nu_+$, respectively.  And,  
$\bS_\pm(\bu_\pm)$ are redefined by 
$$\bS_+(\bu_+) = \mu_+\bD(\bu_+) + (\nu_+-\mu_+)\dv\bu_+\bI,\quad
\bS_-(\bu_-) = \mu_-\bD(\bu_-).$$
To prove Theorem \ref{thm:stokes}, we consider the following 
generalized resolvent problem:
\begin{alignat}2
&\rho_+\lambda\bu_+ - \DV\bS_+(\bu_+) = \bff_+ &\quad&\text{in $\BR^N_+$},
\nonumber\\
&\rho_-\lambda\bu_- - \DV\bS_-(\bu_-) + \nabla\pi_-
= \bff_-, \enskip \dv\bu_- = \tilde f_- 
=\dv\tilde\bff_- &\quad&\text{in $\BR^N_-$},
\nonumber\\
&\mu_-\bD_{mN}(\bu_-)|_- - \mu_+\bD_{mN}(\bu_+)|_+ = g_m, \nonumber \\
&(\mu_-D_{NN}(\bu_-)-\pi_-)|_- = \frac{\rho_-\sigma}{\rho_--\rho_+}\Delta'H
+ g_N, \nonumber \\
&(\mu_+\bD_{NN}(\bu_+) + (\nu_+-\mu_+)\dv\bu_+)|_+
= \frac{\rho_+\sigma}{\rho_--\rho_+}\Delta'H + g_{N+1}, \nonumber\\
&u_{-m}|_- - u_{+m}|_+ = h_m, \nonumber\\
&\lambda H - \Bigl(\frac{\rho_-}{\rho_--\rho_+}u_{-N}|_-
- \frac{\rho_+}{\rho_--\rho_+}u_{+N}|_+\Bigr)
= d, \label{eq:2}
\end{alignat}
which is corresponding to the time dependent problem \eqref{stokes:1},
\eqref{stokes:2*} and \eqref{stokes:3}. Here and in the sequel,
$m$ runs from $1$ through $N-1$. 

Before stating the main result of this section, 
we first introduce the definition of
$\CR$-boundedness.
\begin{dfn}\label{def:2.1}
\label{def:2.1}  A family of operators $\CT \subset 
\CL(X, Y)$ is called  $\CR$-bounded on $\CL(X, Y)$, 
if there exist constants $C > 0$ and $p \in [1, \infty)$ such 
that for any $n\in \BN$, $\{T_j\}_{j=1}^n 
\subset \CT$, $\{f_j\}_{j=1}^n \subset X$ and sequences 
$\{r_j(u)\}_{j=1}^n$ of independent, symmetric, 
$\{-1, 1\}$-valued random variables on $[0, 1]$ there holds the
inequality:
$$\Bigl\{\int^1_0\|\sum_{j=1}^nr_j(u)T_jf_j\|_Y^p\,du\Bigr\}^{\frac1p}
\leq C\Bigl\{\int^1_0\|\sum_{j=1}^n r_j(u)f_j\|_X^p\,du\Bigr\}^{\frac1p}.$$
The smallest such $C$ is called $\CR$-bound of $\CT$, 
which is denoted by $\CR_{\CL(X, Y)}(\CT)$. Here and in the following, 
$\CL(X, Y)$ denotes the set of all bounded linear operators from 
$X$ into $Y$. 
\end{dfn}
The following theorem is a main result for problem \eqref{eq:2}.
\begin{thm}\label{main:r-bound} Let $1 < q < \infty$
 and $0 < \epsilon < \pi/2$.  Set 
\begin{align*}
&\Sigma_\epsilon = \{\lambda = \gamma + i\tau 
\in \BC\setminus\{0\} \mid 
|\arg\lambda| \leq \pi-\epsilon\}, \quad
\Sigma_{\epsilon, \lambda_0} = \{\lambda \in \Sigma_\epsilon \mid
|\lambda| > \lambda_0\} \enskip(\lambda_0 \geq 0), \\
&X_q = \{(\bff_+, \bff_-, f_\dv, \bff_\dv, \bg, \bh, d)\mid
\bff_+ \in L_q(\BR^N_+)^N, \enskip \bff_-, \bff_\dv \in L_q(\BR^N_-)^N, 
\enskip 
f_\dv \in W^1_q(\BR^N), \enskip \\
&\qquad  \bg=(g_1, \ldots, g_{N+1}) \in W^1_q(\BR^N)^{N+1}, 
\enskip \bh=(h_1, \ldots, h_{N-1}) \in W^2_q(\BR^N)^{N-1}, 
\enskip d \in W^2_q(\BR^N)\}, \\
&\CX_q = \{F=(F_{+1}, F_{-1}, F_{-2}, F_{-3}., F_{-4}, F_1, F_2,
F_3, F_4, F_5, F_6) \mid F_{\pm 1} \in L_q(\BR^N_\pm), \enskip \\
&\quad F_{-2}\in L_q(\BR^N_-), F_{-3}, F_{-4} \in L_q(\BR^N_-)^N, \enskip 
F_1\in L_q(\BR^N)^{N+1}, F_2 \in L_q(\BR^N)^{(N+1)N}, \\
&\quad F_3 \in L_q(\BR^N)^{N-1}, 
F_4 \in L_q(\BR^N)^{(N-1)N}, F_5 \in L_q(\BR^N)^{(N-1)N^2},
\enskip F_6 \in W^2_q(\BR^N)\}.
\end{align*}
Then, there exist a constant $\lambda_0 > 0$ and operator families 
$\CA_\pm(\lambda) \in \Hol(\Sigma_{\epsilon, \lambda_0}, 
\CL(X_q, W^2_q(\BR^N_\pm)^N))$,
$\CP_- \in \Hol(\Sigma_{\epsilon, \lambda_0}, \CL(X_q, \hat W^1_q(\BR^N_-)))$, 
and 
$\CH(\lambda) \in  \Hol(\Sigma_{\epsilon, \lambda_0}, \CL(X_q, W^3_q(\BR^N)))$
such that for any $\lambda \in \Sigma_{\epsilon, \lambda_0}$ and 
$\bF = (\bff_+, \bff_-, f_\dv, \bff_\dv, \bg, \bh, d) \in X_q$,
$\bu_\pm = \CA_\pm(\lambda)\bF_\lambda$, $\pi_- = \CP_-(\lambda)\bF_\lambda$
and $H = \CH(\lambda)\bF_\lambda$ 
are unique solutions of problem \eqref{eq:2} and we have 
\begin{align*}
\CR_{\CL(\CX_q, L_q(\BR^N_\pm)^{N+N^2+N^3})}(\{(\tau\pd_\tau)^\ell
G^1_\lambda\CA_\pm(\lambda) \mid \lambda \in \Sigma_{\epsilon, \lambda_0}\})
\leq c \quad(\ell = 0, 1),\\
\CR_{\CL(\CX_q, L_q(\BR^N_-)^N)}(\{(\tau\pd_\tau)^\ell
\nabla\CP_-(\lambda) \mid \lambda \in \Sigma_{\epsilon, \lambda_0}\})
\leq c \quad(\ell = 0, 1), \\
\CR_{\CL(\CX_q, W^2_q(\BR^N)^{N+1})}(\{(\tau\pd_\tau)^\ell
G^2_\lambda\CH(\lambda) \mid \lambda \in \Sigma_{\epsilon, \lambda_0}\})
\leq c \quad(\ell = 0, 1)
\end{align*}
with some constant $c$.  Here, $G^1_\lambda\CA_\pm(\lambda) = 
(\lambda\CA_\pm(\lambda), \lambda^{1/2}\nabla\CA_\pm(\lambda), 
\nabla^2\CA_\pm
(\lambda))$, $G^2_\lambda \CH(\lambda) 
= (\lambda\CH(\lambda), \nabla\CH(\lambda))$, and 
$\bF_\lambda = (F_+, F_-, \lambda^{1/2}f_\dv, \nabla f_\dv, 
\lambda\bff_\dv, \lambda^{1/2}\bg, \nabla\bg, \lambda\bh, 
\lambda^{1/2}\nabla\bh, \nabla^2\bh, d)$.
\end{thm}
\begin{remark}
$F_\pm$, $F_{-2}$, $F_{-3}$, $F_{-4}$, $F_1$, $F_2$, $F_3$, 
$F_4$, $F_5$ and $F_6$ are corresponding variables to 
$\bff_\pm$, $\lambda^{1/2}f_\dv$, $\nabla f_\dv$, $\lambda\bff_\dv$,
$\lambda^{1/2}\bg$, $\nabla\bg$, $\lambda\bh$, $\lambda^{1/2}\nabla\bh$,
$\nabla^2\bh$ and $d$, respectively. 
\end{remark}
To prove Theorem \ref{main:r-bound}, as auxiliary problem, we consider the 
following two equations. 
\begin{align}
&\begin{cases}
\rho_+\lambda\bu_+ - \DV\bS_+(\bu_+) = \bff_+ &\quad\text{in $\BR^N_+$}, \\
\mu_+\bD_{mN}(\bu_+)|_+ = 0, \quad
(\mu_+\bD_{NN}(\bu_+) + (\nu_+-\mu_+)\dv\bu_+)|_+
= g_{N+1}
\end{cases} \label{eq:au4}
\intertext{and}
&\begin{cases}
\rho_-\lambda\bu_- - \DV\bS_-(\bu_-) + \nabla\pi_-
= \bff_-, \enskip \dv\bu_- = \tilde f_- 
=\dv\tilde\bff_- &\quad\text{in $\BR^N_-$}, \\
\mu_-\bD_{mN}(\bu_-)|_- = g_m, \quad 
(\mu_-D_{NN}(\bu_-)-\pi_-)|_- = g_N.
\end{cases}
\label{eq:au5}
\end{align}
The existence of $\CR$ bounded solution operators of \eqref{eq:au4}
and \eqref{eq:au5} were proved in
Shibata \cite{S1} (cf also \cite{SS4}) 
and G\"otz and Shibata \cite{DS}, respectively. 
In fact, we know the following two theorems.
\begin{thm}[\cite{DS}]\label{thm:comp} 
Let $1 < q < \infty$, $0 < \epsilon
< \pi/2$ and $\lambda_0 > 0$.  Set 
\begin{align*}
&Y_{q+} = \{(\bff_+, g) \mid
\bff_+ \in L_q(\BR^N_+)^N, g_{N+1} \in W^1_q(\BR^N_+)\} \\
&\CY_{q+} = \{F=(F_{+1}, \tilde F_{+1}, \tilde F_{+2}) \mid 
F_{+1} \in L_q(\BR^N)^N, \enskip \tilde F_{+1} \in L_q(\BR^N_+),
\enskip \tilde F_{+2} \in L_q(\BR^N_+)^N\}.
\end{align*}
Then, there exist an operator family 
$\CA_{+1}(\lambda) \in \Hol(\Sigma_{\epsilon, \lambda_0}, 
\CL(\CY_{q+}, W^2_q(\BR^N_+)^N))$
such that for any $\lambda \in \Sigma_{\epsilon, \lambda_0}$ and 
$(\bff_+, g_{N+1}) \in Y_{q+}$,
$\bu_+ = \CA_{+1}(\lambda)(\bff_+, \lambda^{1/2}g_{N+1}, \nabla g_{N+1})$ 
is a unique solution of problem \eqref{eq:au4} and we have 
\begin{align*}
\CR_{\CL(\CY_{q+}, L_q(\BR^N_\pm)^{N+N^2+N^3})}(\{(\tau\pd_\tau)^\ell
G^1_\lambda\CA_{+1}(\lambda) \mid \lambda \in \Sigma_{\epsilon, \lambda_0}\})
\leq c \quad(\ell = 0, 1)
\end{align*}
with some constant $c$.  
\end{thm}
\begin{remark} $\tilde F_{+1}$ and $\tilde F_{+2}$ are corresponding
variables to $\lambda^{1/2}g_{N+1}$ and $\nabla g_{N+1}$, respectively.
\end{remark}
\begin{thm}[\cite{S1}]\label{incomp} Let $1 < q < \infty$ and $0 < \epsilon
< \pi/2$.  Set 
\begin{align*}
&Y_{q-} = \{(\bff_-, f_\dv, \bff_\dv, \tilde\bg)\mid
\bff_-, \bff_\dv \in L_q(\BR^N_-), \enskip 
f_\dv \in W^1_q(\BR^N_-), \enskip
 \tilde\bg=(g_1, \ldots, g_N) \in W^1_q(\BR^N)^N\}, \\
&\CY_{q-} = \{F=(F_{-1}, F_{-2}, F_{-3}., F_{-4}, 
\tilde F_{-1}, \tilde F_{-2}) 
\mid F_{-1} \in L_q(\BR^N_-)^N, \enskip \\
&\quad F_{-2}\in L_q(\BR^N_-), F_{-3}, F_{-4} \in L_q(\BR^N_-)^N, \enskip 
\tilde F_{-1} \in L_q(\BR^N)^N, \tilde F_{-2} \in L_q(\BR^N)^{N^2}\}.
\end{align*}
Then, there exist operator families 
$\CA_{-1}(\lambda)$ and $\CP_{-1}(\lambda)$ with
$$\CA_{-1}(\lambda)\in \Hol(\Sigma_{\epsilon},  
\CL(\CY_{q-}, W^2_q(\BR^N_-)^N)), \quad 
\CP_{-1}(\lambda)
\in \Hol(\Sigma_{\epsilon}, \CL(\CY_{q-}, \hat W^1_q(\BR^N_-)))$$ 
such that for any $\lambda \in \Sigma_{\epsilon}$ and 
$\bF = (\bff_-, f_\dv, \bff_\dv, \tilde\bg) \in Y_{q-}$,
$\bu_- = \CA_{-1}(\lambda)\tilde\bF_\lambda$ and  
$\pi_- = \CP_{-1}(\lambda)\tilde\bF_\lambda$
are unique solutions of problem \eqref{eq:au5} and we have 
\begin{align*}
\CR_{\CL(\CY_{q-}, L_q(\BR^N_-)^{N+N^2+N^3})}(\{(\tau\pd_\tau)^\ell
G^1_\lambda\CA_{-1}(\lambda) \mid \lambda \in \Sigma_{\epsilon}\})
\leq c \quad(\ell = 0, 1),\\
\CR_{\CL(\CY_{q-}, L_q(\BR^N_-)^N)}(\{(\tau\pd_\tau)^\ell
\nabla\CP_{-1}(\lambda) \mid \lambda \in \Sigma_{\epsilon}\})
\leq c \quad(\ell = 0, 1)
\end{align*}
with some constant $c$.  Here, 
$\tilde\bF_\lambda = (\bff_-, \lambda^{1/2}f_\dv, \nabla f_\dv, 
\lambda\bff_\dv, \lambda^{1/2}\tilde\bg, \nabla\tilde\bg)$.
\end{thm}
\begin{remark} $\tilde F_{-1}$ and $\tilde F_{-2}$ 
are corresponding
variables to $\lambda^{1/2}\tilde\bg$ and $\nabla \tilde\bg$,
 respectively.
\end{remark}
Thus, it is sufficient to consider problem \eqref{eq:2} with
$\bff_\pm=0$, $f_\dv = 0$, $\bff_\dv=0$ and 
$g_j=0$ for $j = 1, \ldots, N+1$. 
Finally, we consider one more auxiliary problem:
\begin{alignat}2
&\rho_+\lambda\bu_+ - \DV\bS_+(\bu_+) = 0 
&\quad&\text{in $\BR^N_+$},
\nonumber\\
&\rho_-\lambda\bu_- - \DV\bS_-(\bu_-) + \nabla\pi_-
= 0, \quad \dv \bu_- = 0 &\quad&\text{in $\BR^N_-$},
\nonumber\\
&\mu_-\bD_{mN}(\bu_-)|_- - \mu_+\bD_{mN}(\bu_+)|_+ = 0, 
\nonumber \\
&(\mu_-D_{NN}(\bu_-)-\pi_-)|_- = 0, \nonumber \\
&(\mu_+\bD_{NN}(\bu_+) + (\nu_+-\mu_+)\dv\bu_+)|_+
=  0, \nonumber\\
&u_{-m}|_- - u_{+m}|_+ = h_m, \label{eq:prob6}
\end{alignat}
From Sect.3 through Sect.5, we prove the following theorem. 
\begin{thm}\label{thm:au} Let $1 < q < \infty$ and $0 < \epsilon
< \pi/2$.  Set 
\begin{align*}
&Z_{q} = \{\bh = (h_1, \ldots, h_{N-1}) \in W^1_q(\BR^N)^{N-1}\}, \\
&\CZ_{q} = \{F=(F_3, F_4, F_5) \mid 
 F_3 \in L_q(\BR^N)^{N-1}, \quad  F_4\in L_q(\BR^N)^{(N-1)N}, 
F_5 \in L_q(\BR^N)^{(N-1)N^2}\}.
\end{align*}
Then, there exist operator families 
$\CA_{\pm2}(\lambda)$ and $\CP_{-2}(\lambda)$ with
$$\CA_{\pm2}(\lambda)\in \Hol(\Sigma_{\epsilon} 
\CL(\CZ_{q}, W^2_q(\BR^N_-)^N)), \quad 
\CP_{-2}(\lambda)
\in \Hol(\Sigma_{\epsilon}, \CL(\CZ_{q}, \hat W^1_q(\BR^N_-)))$$ 
such that for any $\lambda \in \Sigma_{\epsilon}$ and 
$(\bh, H) \in Z_{q}$,
$\bu_\pm = \CA_{\pm2}(\lambda)\bH_\lambda$ and  
$\pi_- = \CP_{-2}(\lambda)\bH_\lambda$
are unique solutions of problem \eqref{eq:prob6} and we have 
\begin{align*}
\CR_{\CL(\CZ_{q}, L_q(\BR^N_-)^{N+N^2+N^3})}(\{(\tau\pd_\tau)^\ell
G^1_\lambda\CA_{\pm2}(\lambda) \mid \lambda \in \Sigma_{\epsilon}\})
&\leq c \quad(\ell = 0, 1),\\
\CR_{\CL(\CZ_{q}, L_q(\BR^N_-)^N)}(\{(\tau\pd_\tau)^\ell
\nabla\CP_{-2}(\lambda) \mid \lambda \in \Sigma_{\epsilon, \lambda_0}\})
&\leq c \quad(\ell = 0, 1)
\end{align*}
with some constant $c$.  Here, 
$\bH_\lambda = (\lambda \bh, \lambda^{1/2}\nabla\bh, 
\nabla^2\bh)$.
\end{thm}

%%%%%%%%%%%%%%%%%%wit componentwise %%%%%%
\section{Solution formulas for problem \eqref{eq:prob6}} \label{sec:2}

%Setting $A = (\mu_-D_{NN}(\bu_-)-\pi_-)|_-$ and $B = 
%(\mu_+\bD_{NN}(\bu_+) + (\nu_+-\mu_+)\dv\bu_+)|_+$ for short, we have
%$A-B = g_N + \sigma\Delta'H$ and $\rho_-^{-1}A - \rho_+^{-1}B = g_{N+1}$.
%Solving this linear equations, we have
%$$A = \frac{\rho_-\sigma}{\rho_--\rho_+}\Delta'H 
%+ \frac{\rho_-}{\rho_--\rho_+}(g_N - \rho_+g_{N+1}), \quad
%B = \frac{\rho_+\sigma}{\rho_--\rho_+}\Delta'H 
%+ \frac{\rho_+}{\rho_--\rho_+}(g_N - \rho_-g_{N+1}).
%$$

In this section, we consider the following equations:
\begin{alignat}2
&\rho_+\lambda\bu_+ - \DV\bS_+(\bu_+) = 0 &\quad&\text{in $\BR^N_+$},
\nonumber\\
&\rho_-\lambda\bu_- - \DV\bS_-(\bu_-) + \nabla\pi_-
= 0, \quad \dv \bu_- = 0 &\quad&\text{in $\BR^N_-$},
\nonumber\\
&\mu_-\bD_{mN}(\bu_-)|_- - \mu_+\bD_{mN}(\bu_+)|_+ = 0, \nonumber \\
&(\mu_-D_{NN}(\bu_-)-\pi_-)|_- = \sigma_-\Delta'H, \nonumber \\
&(\mu_+\bD_{NN}(\bu_+) + (\nu_+-\mu_+)\dv\bu_+)|_+
=  \sigma_+\Delta'H, \nonumber\\
&u_{-m}|_- - u_{+m}|_+ = h_m. \label{eq:6}
\end{alignat}
where, we have added $\sigma_\pm\Delta'H$ 
with  $\sigma_\pm = \frac{\rho_\pm\sigma}{\rho_--\rho_+}$ to
\eqref{eq:prob6} for the latter use.
Let $\hat v = \CF_{x'}[v](\xi', x_N)$  denote
the partial Fourier transform with respect to the tangential 
variable $x' = (x_1, \ldots, x_{N-1})$
with $\xi' = (\xi_1, \ldots, \xi_{N-1})$ defined by 
$\CF_{x'}[v](\xi', x_N)
= \int_{\BR^{N-1}}e^{-ix'\cdot\xi'}v(x', x_N)\,dx'$. 
Using the formulas:
$$\DV\bS_{+}(\bu_+) = \mu_{+}\Delta\bu_+ + \nu_{+}\nabla\dv\bu_+,
\quad \DV\bS_{-}(\bu_-) = \mu_{-}\Delta\bu_- $$
and applying 
 the partial Fourier transform to \eqref{eq:6}, we transfer problem  
\eqref{eq:6} to  the ordinary differential equations: 
\begin{equation}\label{eq:7}
\begin{cases}
\rho_{+}\lambda \hat u_{+j} + \mu_{+}|\xi'|^2\hat u_{+}
-\mu_{+}D_N^2\hat u_{+j} 
-\nu_{+}i\xi_j
(i\xi'\cdot \hat u'_+ + D_N\hat u_{+N}) = 0
&\quad\text{\,for $x_N > 0$}, \\
\rho_{+}\lambda \hat u_{+N} + \mu_{+}|\xi'|^2\hat u_{+N}
-\mu_{+}D_N^2\hat u_{+N} 
-\nu_{+}D_N
(i\xi'\cdot \hat u'_+ + D_N\hat u_{+N}) = 0
&\quad\text{\,for $x_N > 0$}, \\
\rho_{-}\lambda \hat u_{-j} + \mu_{-}|\xi'|^2\hat u_{-j}
-\mu_{-}D_N^2\hat u_{-j} + i\xi_j\hat \pi_- = 0 
&\quad\text{\,for $x_N < 0$},\\
\rho_{-}\lambda \hat u_{-N} + \mu_{-}|\xi'|^2\hat u_{-N}
-\mu_{-}D_N^2\hat u_{-N} + D_N\hat \pi_- = 0 
&\quad\text{\,for $x_N < 0$}, \\
i\xi'\cdot\hat u'_- +D_N\hat u_{-N} = 0 &\quad\text{\,for $x_N < 0$}, 
\end{cases}
\end{equation}
subject to the interface condition:
\begin{equation}\label{eq:bdy1}
\begin{cases}
\mu_{-}(D_N\hat u_{-m} + i\xi_m\hat u_{-N})|_-
- \mu_{+}(D_N\hat u_{-m} + i\xi_m\hat u_{-N})|_+
= 0,\\
(2\mu_{-}D_N\hat u_{-N} - \hat \pi_-)|_- = -\sigma_-A^2\hat H(0), \\
(2\mu_{+}D_N\hat u_{-N} + (\nu_{+}-\mu_{+})
(i\xi'\cdot\hat u'_+ + D_N\hat u_{+N})|_+ 
= -\sigma_+A^2\hat H(0), \\
\hat u_{-m}|_- - \hat u_{+m}|_+ = \hat h_m(0) 
\end{cases}
\end{equation}
where $D_N = d/dx_N$ and $i\xi'\cdot \hat v' 
= \sum_{\ell=1}^{N-1} i\xi_\ell \hat v_\ell$ for 
$\bv = (v_1, \ldots, v_{N-1}, v_N)$. 
Here and in the sequel, $j$ also runs from $1$ through $N-1$.
Applying the divergence to the first and second equations 
in \eqref{eq:6}, 
we have $\rho_{+}\lambda\dv \bu_+ - (\mu_{+} + \nu_{+})
\Delta \dv \bu_+ = 0$ in $\BR^N_+$ and $\Delta p_- = 0$ in $\BR^N_-$, 
so that 
$$(\rho_{+}\lambda - (\mu_{+} + \nu_{+})\Delta)
(\rho_{+}\lambda - \mu_{+}\Delta)\bu_+ = 0\enskip\text{ in $\BR^N_+$},
\quad (\rho_{-}\lambda -\Delta)\Delta \bu_-=0\quad\text{ in $\BR^N_-$}.
$$
Thus, the characteristic roots of 
\eqref{eq:7} are
\begin{equation}\label{eq:ch1}
A_+ = \sqrt{\rho_{+}(\mu_{+} + \nu_{+})^{-1}\lambda + A^2}, \quad 
B_\pm = \sqrt{\rho_{\pm}(\mu_{\pm})^{-1}\lambda + A^2}, \quad A = |\xi'|.
\end{equation}

To state our solution formulas of problem: \eqref{eq:7}-
\eqref{eq:bdy1}, we introduce some classes of multipliers. 
\begin{dfn}\label{def:symbols} Let $0 < \epsilon < \pi/2$, $\lambda_0 \geq 0$, 
 and let $s$ be a real number.  Set
\begin{gather*}
%\Sigma_\epsilon = \{\lambda \in \BC\setminus\{0\} \mid 
%|\lambda| \leq \pi-\epsilon\}, \quad
%\Sigma_{\epsilon, \lambda_0} = \{\lambda \in \Sigma_\epsilon
%\mid |\lambda| > \lambda_0\}, \\ 
\tilde\Sigma_{\epsilon, \lambda_0}
= \{(\lambda, \xi') \mid 
\lambda = \gamma + i\tau \in \Sigma_{\epsilon, \lambda_0}, \enskip
\xi' = (\xi_1, \ldots, \xi_{N-1}) \in \BR^{N-1}\setminus\{0\}\}.
\end{gather*}
%Note that $\Sigma_\epsilon = \Sigma_{\epsilon, 0}$ in this notation. 
Let $m(\lambda, \xi')$ be a function defined on 
$\tilde\Sigma_{\epsilon, \lambda_0}$.
\begin{itemize}
\item[\thetag1]~$m(\lambda, \xi')$ 
is called a multiplier of order $s$ 
with type $1$ if for any multi-index 
$\kappa' = (\kappa_1, \ldots, \kappa_{N-1})
\in \BN_0^{N-1}$ and $(\lambda, \xi') 
\in \tilde\Sigma_{\epsilon, \lambda_0}$
there exists a constant $C_{\kappa'}$ depending on $\kappa'$, 
$\epsilon$, $\mu_\pm$, $\nu_+$ and $\rho_{\pm}$
 such that there holds the estimates: 
\begin{equation}\label{eq:3-0*}
|\pd_{\xi'}^{\kappa'}m(\lambda, \xi')| \leq C_{\alpha'}
(|\lambda|^{1/2} + A)^{s-|\kappa'|}, \quad 
|\pd_{\xi'}^{\kappa'}(\tau\frac{\pd m}{\pd\tau}(\lambda, \xi'))| 
\leq C_{\kappa'}(|\lambda|^{1/2} + A)^{s-|\kappa'|}.
\end{equation}
\item[\thetag2]~$m(\lambda, \xi')$ 
is called a multiplier of order $s$ 
with type $2$ if for any multi-index 
$\kappa' = (\kappa_1, \ldots, \kappa_{N-1})
\in \BN_0^{N-1}$ and $(\lambda, \xi') 
\in \tilde\Sigma_{\epsilon, \lambda_0}$
there exists a constant $C_{\kappa'}$ depending on $\kappa'$, 
$\epsilon$, $\mu_\pm$, $\nu_+$, and $\rho_{\pm}$ 
such that there holds the estimates: 
\begin{equation}\label{eq:3-0}
|\pd_{\xi'}^{\kappa'}m(\lambda, \xi')| \leq C_{\kappa'}
(|\lambda|^{1/2} + A)^sA^{-|\kappa'|}, \quad 
|\pd_{\xi'}^{\kappa'}(\tau\frac{\pd m}{\pd\tau}(\lambda, \xi'))| 
\leq C_{\kappa'}(|\lambda|^{1/2} + A)^sA^{-|\kappa'|}.
\end{equation}
\end{itemize}
Let $\bM_{s,i}(\lambda_0)$ 
be the set of all multipliers of order $s$ with type $i$
($i = 1, 2$). 
\end{dfn}
Obviously, $\bM_{s,i}(\lambda_0)$ are vector spaces on $\BC$
and  
for $0 \leq \lambda_0 < \lambda_1$, $\bM_{s,i}(\lambda_0)
\supset \bM_{s,i}(\lambda_1)$.  Moreover, 
by the fact: $||\lambda|^{1/2} + A|^{-|\alpha'|}
\leq A^{-|\alpha'|}$ and the Leibniz rule, 
we have the following lemma immediately. 
\begin{lem}\label{lem:sym} Let $s_1$, $s_2$ be two real numbers. 
Then, the following three assertions hold.
\begin{itemize}
\item[\thetag1]~ Given $m_i \in \bM_{s_i, 1}(\lambda_0)$ $(i = 1 ,2)$, we have
$m_1m_2 \in \bM_{s_1+s_2, 1}(\lambda_0)$. 
\item[\thetag2]~ Given $\ell_i \in \bM_{s_i, i}(\lambda_0)$ $(i = 1 ,2)$, 
we have
$\ell_1\ell_2 \in \bM_{s_1+s_2, 2}(\lambda_0)$. 
\item[\thetag3]~ Given $n_i \in \bM_{s_i, 2}(\lambda_0)$ 
$(i = 1 ,2)$, we have
$n_1n_2 \in \bM_{s_1+s_2, 2}(\lambda_0)$. 
\end{itemize}
\end{lem}
\begin{remark} \thetag1~We see easily that $i\xi_j \in \bM_{1, 1}(0)$ 
($j=1, \ldots, N-1$), so that $A^2 \in \bM_{2,1}(0)$.  On the other hand, 
$A \in \bM_{1,2}(0)$ and $A^{-1} \in \bM_{-1,2}(0)$.  Especially, 
$i\xi_j/A \in \bM_{0,2}(0)$. \\
\thetag2~$\bM_{s, 1}(\lambda_0)
\subset \bM_{s,2}(\lambda_0)$ for any $s \in \BR$
and $\lambda_0 \geq0$.
\end{remark}
In this section, first of all we show the following solution formulas for 
problem \eqref{eq:7}-\eqref{eq:bdy1}:
\begin{equation}\label{eq:sol-form}\begin{split}
\hat u_{+J} &= \sum_{k=1}^3 \hat u^+_{Jk}, \quad 
\hat u_{-J} = \sum_{k=1}^3 \hat u^-_{Jk}, \quad 
\hat p_-  = e^{Ax_N}\Bigl\{
\sum_{m=1}^{N-1}p^-_{m, 1}\hat h_m(0) + Ap^-_{N,1}\hat H(0)\Bigr\},
\\
\hat u^\pm_{J1} & = AM_\pm(x_N)\Bigl\{
\sum_{m=1}^{N-1}R^\pm_{Jm, 0}\hat h_m(0)
+ AR^\pm_{JN,0}\hat H(0)\Bigr\} \\
\hat u^\pm_{J2} & = Ae^{\mp B_\pm x_N}\Bigl\{
\sum_{m=1}^{N-1}S^\pm_{Jm, -1}\hat h_m(0)
+ AS^\pm_{JN,-1}\hat H(0)\Bigr\} \\
\hat u^\pm_{j3} & =e^{\mp B_\pm x_N}T^\pm_{j,0}\hat h_j(0),
\quad
\hat u^\pm_{N,3} = 0
\end{split}\end{equation}
with 
\begin{gather}
R^\pm_{Jm, 0} \in \bM_{0, 2}(0), 
\enskip R^\pm_{JN,0} \in \bM_{0,2}(0),
\enskip 
S^\pm_{Jm, -1} \in \bM_{-1,2}(0), \enskip S^\pm_{JN,-1}
\in \bM_{-1,2}(0) \nonumber \\
T^\pm_{j, -1} \in \bM_{-1, 1}(0), \enskip 
T^\pm_{j, 0}
\in \bM_{0,1}(0), \enskip 
 p^-_{m, 1} \in
\bM_{1, 2}(0), \enskip p^-_{N, 1} \in 
\bM_{1,2}(0). \label{sym-est:0}
\end{gather}
Here and in the sequel, $J$  runs from $1$
through $N$. Recall that $j$ and $m$ run from $1$ through $N-1$, 
respectively. Moreover,   
 $M_\pm(x_N)$
denote the Stokes kernels defined by 
\begin{equation}\label{symb:1}
M_+(x_N) = \frac{e^{-B_+x_N} - e^{-A_+x_N}}{B_+-A_+}, \quad
M_-(x_N) = \frac{e^{B_-x_N} - e^{Ax_N}}{B_--A}.
\end{equation}

In the sequel, we prove \eqref{eq:sol-form}. 
We look for solutions $\hat u_{\pm J}$
and $\hat p_-$ to problem: \eqref{eq:7}-
\eqref{eq:bdy1} of the forms:
\begin{equation}\label{eq:sol1}\begin{split}
\hat u_{+J} & 
= \alpha_{+J}(e^{-B_+x_N} - e^{-A_+x_N}) + \beta_{+J}e^{-B_+x_N}, \\
\hat u_{-J} &  
= \alpha_{-J}(e^{B_-x_N} - e^{Ax_N}) + \beta_{-J}e^{B_-x_N},
\quad \hat p_- = \gamma_-e^{Ax_N}.
\end{split}\end{equation}
Using the symbols $B_\pm$, we write \eqref{eq:7} as follows:
\begin{equation}\label{eq:eq2}
\begin{cases}
\mu_{+}B_+^2\hat u_{+j} - \mu_+D_N^2\hat u_{+j}
-\nu_{+}i\xi_j(i\xi'\cdot\hat u_+' + 
D_N\hat u_{+N}) = 0&\quad (x_N > 0), \\
\mu_{+}B_+^2\hat u_{+N} - \mu_{+}D_N^2\hat u_{+N}
-\nu_{+}D_N(i\xi'\cdot\hat u_+' + 
D_N\hat u_{+N})= 0 &\quad (x_N > 0),\\
\mu_{-}B_-^2\hat u_{-j} - \mu_{-}D_N^2\hat u_{-j}
+i\xi_j\hat \pi_- = 0 &\quad (x_N < 0), \\
\mu_{-}B_-^2\hat u_{-N} - \mu_{-}D_N^2\hat u_{-N}
+D_N\hat \pi_- = 0 &\quad (x_N < 0), \\
i\xi'\cdot \hat u'_- + D_N \hat u_{-N}  = 0
&\quad (x_N < 0).
\end{cases}\end{equation}
Substituting the  formulas of $\hat u_{\pm J}$ 
 into \eqref{eq:eq2} and equating the coefficients of 
$e^{\mp B_\pm x_N}$,
$e^{-A_+x_N}$ and $e^{Ax_N}$, we have 
\allowdisplaybreaks{
\begin{align}
&\mu_{+}(A_+^2-B_+^2)\alpha_{+j} + \nu_{+}
i\xi_j(i\xi'\cdot\alpha_+' - A_+\alpha_{+N}) = 0, \nonumber\\
&\mu_{+}(A_+^2-B_+^2)\alpha_{+N} - \nu_{+}
A_+(i\xi'\cdot\alpha_+' - A_+\alpha_{+N}) = 0, \nonumber \\
&i\xi'\cdot\alpha'_+ - \alpha_{+N}B_+ + i\xi'\cdot\beta'_+
- \beta_{+N}B_+ = 0, \nonumber \\
&\mu_{-}(A^2 - B_-^2)\alpha_{-j} + i\xi_j\gamma_- = 0, 
\quad \mu_{-}(A^2 - B_-^2)\alpha_{-N} + A\gamma_- = 0, \nonumber\\
&i\xi'\cdot\alpha'_- + \alpha_{-N}B_- + i\xi'\cdot\beta'_-
+ \beta_{-N}B_- = 0,
\quad i\xi'\cdot\alpha'_- + A\alpha_{-N} = 0.
\label{eq:eq3}
\end{align}
First, we represent $i\xi'\cdot\alpha'_\pm$, $\alpha_{\pm N}$
and $\gamma_-$ by $i\xi'\cdot\beta'_\pm$ and $\beta_{\pm N}$. 
It follows from \eqref{eq:eq3} that 
\begin{alignat}2
i\xi'\cdot\alpha'_+ & = \frac{A^2}{A_+B_+ - A^2}(i\xi'\cdot\beta'_+
-B_+\beta_{+N}), &\quad
\alpha_{+N} &= \frac{A_+}{A_+B_+ - A^2}(i\xi'\cdot\beta'_+ - B_+\beta_{+N}), 
\nonumber \\
i\xi'\cdot\alpha'_- & = \frac{A}{B_--A}(i\xi'\cdot\beta'_- + B_-\beta_{-N}), 
&\quad \alpha_{-N} &= \frac{-1}{B_--A}(i\xi'\cdot\beta'_- + B_-\beta_{-N}), 
\label{eq:eq4}\\
\gamma_- & = -\frac{\mu_{-}(A+B_-)}{A}(i\xi'\cdot\beta'_- + B_-\beta_{-N}).
\nonumber
\end{alignat}
%%%%%%%%%%%%%%%%%%%%%%%%%%%
Substituting the relations: 
$$\hat u_{\pm J}(0) = \beta_{\pm J}, \enskip  
D_N\hat u_{+J}(0) = (A_+ - B_+)\alpha_{+J} - B_+\beta_{+J},
\enskip 
D_N\hat u_{-J}(0) = (B_--A)\alpha_{-J} + B_-\beta_{-J}
$$
into \eqref{eq:bdy1}, we have 
\begin{equation}\label{eq:eq4*}\begin{split}
&\beta_{+m} = \beta_{-m} - \hat h_m(0), \\
&\mu_{+}((B_+-A_+)\alpha_{+m} + B_+\beta_{+m} - i\xi_m\beta_{+N})
+\mu_{-}((B_--A)\alpha_{-m} + B_-\beta_{-m} 
+i\xi_m\beta_{-N}) = 0, \\
&2\mu_{-}((B_--A)\alpha_{-N}
+ B_-\beta_{-N}) - \gamma_- = -\sigma_-A^2\hat H(0), \\
&2\mu_{+}((B_+- A_+)\alpha_{+N} + B_+\beta_{+N})
+(\nu_{+} - \mu_{+})(-i\xi'\cdot\beta'_+ 
+ (B_+-A_+)\alpha_{+N}
+ B_+\beta_{+N})=\sigma_+A^2\hat H(0).
\end{split}\end{equation}
Using \eqref{eq:eq4} and \eqref{eq:eq4*},  we have
\begin{equation}\label{eq:eq5**}\begin{split}
0 & = L_{11}^+(i\xi'\cdot\beta'_+) + 
L_{11}^-(i\xi'\cdot\beta'_-) + 
L_{12}^+\beta_{+N} + L_{12}^-\beta_{-N}, \\
-\sigma_-A^3\hat H(0) & = L_{21}^-i\xi'\cdot\beta'_- + L_{22}^-\beta_{-N}, \\
-\sigma_+A^2\hat H(0) & = -L_{21}^+i\xi'\cdot\beta_+' - L_{22}^+\beta_{+N}
\end{split}\end{equation}
%%%%%%%%%%%%%%%%%%
with 
\begin{alignat}2
L_{11}^+ & = \mu_{+}\frac{A_+(B_+^2 - A^2)}{A_+B_+ - A^2}, 
&\quad L^-_{11} &= \mu_{-}(A+B_-), \nonumber \\
L_{12}^+ & = \mu_{+}\frac{A^2(2A_+B_+ - A^2-B_+^2)}{A_+B_+ - A^2},  
&\quad L^-_{12} &=  \mu_{-}A(B_--A), \nonumber \\ 
L_{21}^+ & = 2\mu_{+}\frac{A_+(B_+-A_+)}{A_+B_+ - A^2} 
- (\nu_{+}-\mu_{+})\frac{A_+^2-A^2}{A_+B_+ - A^2},
&\quad L^-_{21} &= \mu_-(B_--A),  \nonumber\\
L^+_{22} & = (\mu_{+}+\nu_{+})
\frac{B_+(A_+^2-A^2)}{A_+B_+ - A^2},
&\quad L^-_{22} &= \mu_-(A+B_-)B_-.
\label{eq:eq5}
\end{alignat}
As is seen in Sect. 4, we have 
\begin{equation}\label{sym-est:1}
L_{11}^+, L_{22}^+ \in \bM_{1,1}(0), \enskip
L_{12}^+ \in \bM_{2,1}(0), \enskip L_{21}^+ \in \bM_{0,1}(0), \enskip
L_{11}^-, L^-_{21} \in \bM_{1,2}(0), \enskip
 \enskip L_{12}^-, L_{22}^- \in \bM_{2,2}(0).
\end{equation}
Using $i\xi'\cdot\beta_+' 
= i\xi'\cdot\beta_-' - i\xi'\cdot \hat h'(0)$,
we write the linear equation \eqref{eq:eq5**} 
 in the following form:
\begin{equation}\label{eq:eq5*}
L\left( \begin{matrix}i\xi'\cdot\beta'_-\\
\beta_{+} \\ \beta_-\end{matrix}\right)
= \left(\begin{matrix}  
L_{11}^+i\xi'\cdot \hat h'(0)  \\
-\sigma_-A^3\hat H(0) \\ -\sigma_+A^2\hat H(0) - L_{21}^+i\xi'\cdot\hat h'(0)  
\end{matrix}\right)
\quad \text{with}\quad
L = \left(\begin{matrix}
L_{11}^++L_{11}^- & L_{12}^+ & L_{12}^- \\
L_{21}^- & 0 & L_{22}^- \\
-L_{21}^+ & -L_{22}^+ & 0
\end{matrix}\right).
\end{equation}
Moreover, we have
\begin{equation}\label{eq:eq6}
L^{-1} = \frac{1}{\det L}\left(\begin{matrix}
\CL_{11} & \CL_{12} & \CL_{13} \\
\CL_{21} & \CL_{22} & \CL_{23} \\
\CL_{31} & \CL_{32} & \CL_{33}
\end{matrix}\right)
\end{equation}
with
\begin{alignat*}3
\CL_{11} &= L_{22}^+L_{22}^-, 
&\quad 
\CL_{12} &= -L_{22}^+L_{12}^-, 
&\quad
\CL_{13} &= L_{12}^+L_{22}^- \\
\CL_{21} &= -L_{21}^+L_{22}^-,
&\quad
\CL_{22} &= L_{21}^+L_{12}^-, &\quad
\CL_{23} &= L_{12}^-L_{21}^--(L_{11}^++L_{11}^-)L_{22}^-, \\
\CL_{31} &= -L_{22}^+L_{21}^-, &\quad 
\CL_{32} &= (L_{11}^++L_{11}^-)L_{22}^+ - L_{12}^+L_{21}^+, &\quad
\CL_{33} &= -L_{12}^+L_{21}^-.
\end{alignat*}
By \eqref{sym-est:1} and Lemma \ref{lem:sym}, we have
\begin{alignat}3
\CL_{11} &\in \bM_{3,2}(0), &\quad \CL_{12} &\in \bM_{3,2}(0), &\quad
\CL_{13} &\in \bM_{4,2}(0) \nonumber \\
\CL_{i1} &\in \bM_{2,2}(0), &\quad \CL_{i2} &\in \bM_{2,2}(0), &\quad
\CL_{i3} &\in \bM_{3,2}(0)\quad(i = 2, 3). 
\label{sym-est:4}
\end{alignat}
The most important fact of this paper is that $\det L \not=0$
for any $(\lambda, \xi') \in \tilde\Sigma_{\epsilon, 0}$ and 
\begin{equation}\label{sym-est:3}
(\det L)^{-1} \in \bM_{-4,2}(0).
\end{equation}
%%%%%%%%%%%%%%%%%%%%%%%%
%%%%%%%%%%%%%%%%%%%%%%%%
From \eqref{eq:eq5*} and \eqref{eq:eq6} it follows that 
\begin{equation}\label{eq:eq7}\begin{split}
i\xi'\cdot\beta'_- &=
\frac{A}{\det L}\{(\CL_{11}L_{11}^+-\CL_{13}L_{21}^+)i\tilde\xi'\cdot\hat h'(0)
-(\CL_{12}\sigma_-A^2 + \CL_{13}\sigma_+A)\hat H(0)\},
\\
\beta_{+N} &=
\frac{A}{\det L}\{(\CL_{21}L_{11}^+ - \CL_{23}L_{21}^+)
i\tilde\xi'\cdot\hat h'(0)-
(\CL_{22}\sigma_-A^2 + \CL_{23}\sigma_+A)
\hat H(0)\}, \\
\beta_{-N} &=
\frac{A}{\det L}\{(\CL_{31}L_{11}^+ - \CL_{33}L_{21}^+)
i\tilde\xi'\cdot\hat h'(0)
-(\CL_{32}\sigma_-A^2 + \CL_{33}\sigma_+A)\hat H(0)\}
\end{split}\end{equation}
where $i\tilde\xi'\cdot \hat k'(0) = i\sum_{m=1}^{N-1}
\xi_mA^{-1}\hat k_m(0)$ with $k=g$ and $k=h$. 
%%%%%%%%%%%%%%%%
Using the relations: $\beta_{+m} = \beta_{-m} - \hat h_m(0)$,
 by \eqref{eq:eq7},  we have
\begin{equation}\label{rep:1-1}
i\xi'\cdot\beta'_\pm \mp B_\pm\beta_{\pm N} = 
 A\Bigl\{\sum_{m=1}^{N-1}P^\pm_{m, 0}\hat h_m(0)
 + AP^\pm_{N,0}\hat H(0)\Bigr\}
\end{equation}
with 
\begin{align*}
P^+_{m,0} & = \frac{1}{\det L}\{(\CL_{11}-B_+\CL_{21})L^+_{11}
- (\CL_{13}-B_+\CL_{23})L^+_{21})\}\frac{i\xi_m}{A}-\frac{i\xi_m}{A} \\
P^+_{N,0} & = \frac{-1}{\det L}\{(\CL_{12}-B_+\CL_{22})\sigma_- A
+ (\CL_{13}-B_+\CL_{23})\sigma_+\}, \\
P^-_{m,0} & = \frac{1}{\det L}\{(\CL_{11}+B_-\CL_{31})L^+_{11}
- (\CL_{13}+B_-\CL_{33})L^+_{21})\}\frac{i\xi_m}{A}, \\
P^-_{N,0} & = \frac{-1}{\det L}\{(\CL_{12}+B_-\CL_{32})\sigma_- A
+ (\CL_{13}+B_-\CL_{33})\sigma_+\}.
\end{align*}
By Lemma \ref{lem:sym}, \eqref{sym-est:1},
 and \eqref{sym-est:3}, we have 
\begin{equation}\label{sym-est:5}
P^\pm_{m, 0} \in \bM_{0, 2}(0), \quad P^\pm_{N,0} \in \bM_{0,2}(0).
\end{equation}
By \eqref{eq:eq4} we have
\begin{align*}
\hat p_-(x_N) &= -\mu_{-}\frac{(A+B_-)}{A}(i\xi'\cdot\beta'_- + 
B_-\beta_{-N})e^{Ax_N} \\
&= -\mu_-(A+B_-)\Bigl\{
\sum_{m=1}^{N-1}P^-_{m, 0}\hat h_\ell(0)
+ AP^-_{N,0}\hat H(0)\Bigr\}e^{Ax_N},
\end{align*}
so that setting 
$p^-_{m, 1} = -\mu_-(A+B_-)P^-_{m, 0}$ and $p^-_{N,1}
=-\mu_-(A+B_-)P^-_{N,0}$, we have the formula of
$\hat p_-(x_N)$ in \eqref{eq:sol-form}.
%%%%%%%%%%%%%%%%%%%%%%%%%%%%%%%
%%%%%%%%%%%%%%%%%%%%%%%%%%% 

By \eqref{eq:eq3}, we have
\begin{align*}
(B_+-A_+)\alpha_{+j} & = 
\frac{\nu_+ i\xi_j}{\mu_{+}(A_++B_+)}
(i\xi'\cdot \alpha'_+ - A_+\alpha_{+N}), \\
(B_+-A_+)\alpha_{+N} & = -
\frac{\nu_+A_+}{\mu_{+}(A_++B_+)}
(i\xi'\cdot \alpha'_+ - A_+\alpha_{+N}), \\
(B_--A)\alpha_{-j} & = -\frac{i\xi_j}{A}(i\xi'\cdot\beta'_- + B_-\beta_{-N}),
\quad
(B_- - A)\alpha_{-N}  = -(i\xi'\cdot\beta'_- + B_-\beta_{-N}).
\end{align*}
Since 
$
i\xi'\cdot\alpha'_+ - A_+\alpha_{+N} = \frac{A^2-A_+^2}{A_+B_+ - A^2}
(i\xi'\cdot\beta_+' - B_+\beta_{-N})$ 
as follows from \eqref{eq:eq4}, 
setting $A_-=A$, by \eqref{rep:1-1} we have 
\begin{equation}\label{rep:1-2}
(B_\pm-A_\pm)\alpha_{\pm J}
 = A\{\sum_{m=1}^{N-1}R^\pm_{Jm,0}\hat h_m(0) + AR^\pm_{JN,0}\hat H(0)\} 
\end{equation}
with
\begin{align*}
R^+_{jm, 0} 
&=\frac{\nu_+i\xi_jP^+_{m, 0}}{\mu_{+}(A_++B_+)}
\frac{A^2-A_+^2}{A_+B_+-A^2}, \quad
R^+_{jN,0} =\frac{\nu_+i\xi_jP^+_{N,0}}{\mu_{+}(A_++B_+)}
\frac{A^2-A_+^2}{A_+B_+-A^2}, \\
R^+_{Nm,0} &= -\frac{\nu_+A_+P^+_{m,0}}{\mu_{+}(A_++B_+)}
\frac{A^2-A_+^2}{A_+B_+-A^2},
\quad
R^+_{NN,0} = -\frac{\nu_+A_+P^+_{N,0}}{\mu_{+}(A_++B_+)}
\frac{A^2-A_+^2}{A_+B_+-A^2}, \\
R^-_{jm, 0} 
&=-\frac{i\xi_j}{A}P^+_{m, 0}, \quad
R^-_{jN,1} =-\frac{i\xi_j}{A}P^-_{N,1}, \quad
R^+_{Nm,0}  = -P^-_{m,0}, \quad 
R^-_{NN,0} = -P^-_{N,0}. 
\end{align*}
Recalling $A_- = A$, we have 
 $(e^{\mp B_\pm x_N} - e^{\mp A_\pm x_N})\alpha_{\pm J}
= M_\pm(x_N)(B_\pm-A_\pm)\alpha_{\pm J}$. 
Thus, if we set 
  $\hat u^\pm_{J1}= AM_\pm(x_N)(\sum_{m=1}^{N-1}R^\pm_{Jm,0}\hat h_m(0)
+ R^\pm_{JN,1}\hat H(0))$, 
then  $\hat u^\pm_{J1} = \alpha_{\pm J}(e^{\mp B_\pm x_N}- 
e^{\mp A_\pm x_N})$. 
As is seen in Sect. 4 below, we have
\begin{equation}\label{sym-est:6}
A_+ \in \bM_{1,1}(0), \enskip B_+ \in \bM_{1,1}(0), 
\enskip (A_+ + B_+)^{-1} \in 
\bM_{-1,1}(0), \enskip 
\frac{A^2-A_+^2}{A_+B_+ - A^2} \in \bM_{0, 1}(0),
\end{equation}
which, combined with \eqref{sym-est:5}, furnishes that 
$R^+_{Jm, 0} \in \bM_{0,2}(0)$ and $R^+_{JN,1}
\in \bM_{1,2}(0)$. And also, 
by \eqref{sym-est:5} and \eqref{sym-est:6}, we have 
 $R^-_{Jm, 0} \in \bM_{0, 2}(0)$ and $R^-_{JN,1} \in \bM_{1,2}(0)$.

%%%%%%%%%%%%%%%%

>From \eqref{eq:eq4*} it follows that 
\begin{align*}
\beta_{\pm j} &= \frac{\mp\mu_\mp B_\mp}{\mu_{+}B_++\mu_{-}B_-}\hat h_j(0)\\
&\quad
+\frac{1}{\mu_{+}B_++\mu_{-}B_-}\{- \mu_{+}(B_+-A)\alpha_{+j}
-\mu_{-}(B_--A)\alpha_{-j}+i\xi_j\mu_+\beta_{+N}
-i\xi_j\mu_-\beta_{-N}\}.
\end{align*}
We set 
$$T^\pm_{j,0} = 
\frac{\mp\mu_{\mp}B_{\mp}}{\mu_{+}B_++\mu_{-}B_-}$$
and in  view of \eqref{rep:1-2} and \eqref{eq:eq7} 
we set 
\begin{align*}
S^{\pm}_{jm,-1} &= \frac{-1}{\mu_{+}B_++\mu_{-}B_-}\Bigl(
\mu_{+}R^+_{jm,0} +\mu_{-}R^-_{jm,0} \\
&\qquad-\frac{\mu_{+}i\xi_j
(\CL_{21}L^+_{11} - \CL_{23}L_{21}^+)}{\det L}\frac{i\xi_m}{A}  
+\frac{\mu_{-}i\xi_j
(\CL_{31}L^+_{11} - \CL_{33}L_{21}^+)}{\det L}\frac{i\xi_m}{A}
\Bigr), \\
S^\pm_{jN,-1} & = \frac{-1}{\mu_+B_++\mu_-B_-}
\Bigl(\mu_+R^+_{jN,0} + \mu_-R^-_{jN,0} \\
&\qquad
+ \frac{\mu_+i\xi_j(\CL_{22}\sigma_-A + \CL_{23}\sigma_+)}{\det L}
- \frac{\mu_-i\xi_j(\CL_{32}\sigma_-A + \CL_{33}\sigma_+)}{\det L}
\Bigr).
\end{align*}
Thus, if we set  
$\hat u^\pm_{j2} = Ae^{\mp B_\pm x_N}
\{\sum_{m=1}^{N-1}S^\pm_{jm,-1}\hat h_m(0) + AS^\pm_{jN,-1}
\hat H(0)\}$
and $\hat u^\pm_{j3} = e^{\mp B_\pm x_N}T^\pm_{j,0}\hat h_j(0))$,    
then we have 
$\beta_{\pm j}e^{\mp B_\pm x_N} = \hat u^\pm_{j2} + \hat u^\pm_{j3}$. 
Moreover, by \eqref{sym-est:1}, \eqref{sym-est:4}, \eqref{sym-est:3}, 
 \eqref{sym-est:5}, \eqref{sym-est:6}, 
we have $S^\pm_{jm,-1} \in 
\bM_{-1,2}(0)$, $S^\pm_{jN,-1} \in \bM_{-1,2}(0)$, and $T^\pm_{j,0} \in 
\bM_{0,1}(0)$.  

Finally, in view of \eqref{eq:eq7}, 
we define $\hat u^\pm_{N2} = Ae^{\mp B_\pm x_N}\{
\sum_{m=1}^{N-1} S^\pm_{Nm, -1}\hat h_m(0) + AS^\pm_{NN,-1}
\hat H(0)\}$ with  
\begin{alignat*}2
S^+_{Nm,-1} & = \frac{\CL_{21}L_{11}^+ - \CL_{23}L^+_{21}}{\det L}
\frac{i\xi_m}{A}, &\quad 
S^+_{NN, -1} &= \frac{-(\CL_{22}\sigma_-A + \CL_{23}\sigma_+)}
{\det L}, \\
S^-_{Nm,-1} & = \frac{\CL_{31}L_{11}^+ - \CL_{33}L^+_{21}}{\det L}
\frac{i\xi_m}{A}, &\quad 
S^-_{NN, -1} &= \frac{-(\CL_{32}\sigma_-A + \CL_{33}\sigma_)}
{\det L}.
\end{alignat*}
Thus, we have  $\hat u^\pm_{N2}
= \beta_{\pm N}e^{\mp B_\pm x_N}$. Morevoer, by \eqref{sym-est:1},
\eqref{sym-est:3} and \eqref{sym-est:4}, we have $S^\pm_{Nm,-1} \in 
\bM_{-1,2}(0)$ and $S^\pm_{NN,-1} \in \bM_{-1,1}(0)$.  
 Summing up, we have obtained \eqref{eq:sol-form} and 
\eqref{sym-est:0}.

%%%%%%%%%%%%%

%%%%%%%%%%%%%%%%%%%%%%%
To prove Theorem \ref{thm:au}, we consider problem \eqref{eq:prob6},
namely problem \eqref{eq:6} with $H=0$.
To construct our solution operator from the solution formulas in
\eqref{eq:sol-form} with $\hat H=0$, 
first of all we observe the following formulas due to
Volevich:
$$
a(\xi', x_N)\hat h(0) = -\int^{\pm\infty}_0
\{(\pd_Na)(\xi', x_N+y_N)\hat h(y_N) 
+ a(\xi', x_N+y_N)\widehat{\pd_Nh}(\xi', y_N)\}\,dy_N,
$$ 
where $\pd_j  = \pd /\pd x_j$. 
%Using the identity:  $1 = \frac{\rho_{*\pm}\lambda}{\mu_{*\pm} B^2_\pm}
%- \sum_{k=1}^{N-1}\frac{(i\xi_k)(i\xi_k)}{B_\pm^2}$, we write
%\begin{align*}
%&a(\xi', x_N)\hat g(\xi', 0) 
% = -\int^{\pm\infty}_0a(\xi', x_N+y_N)
%\widehat{\pd_Ng}(\xi', y_N)\,dy_N \\
%&- \int^{\pm\infty}_0 \frac{(\pd_Ng)(\xi', x_N+y_N)
%\rho_{*\pm}\lambda^{1/2}}{\mu_{*\pm} B_\pm^2}
%\lambda^{1/2}\hat g(\xi', y_N)\,dy_N
%+ \sum_{k=1}^{N-1}\int^{\pm\infty}_0
%\frac{(\pd_Na)(\xi', x_N+y_N)i\xi_k}{B^2_\pm}
%\widehat{\pd_k g}(\xi', y_N)\,dy_N.
%\end{align*}
%Let $\CF^{-1}_{\xi'}$ denote the partial Fourier inverse 
%transform with respect to $\xi'$ variable and let $f_1$ and $f_2 = (f_{21},
% \ldots, f_{2N})$ be corresponding variables to $\lambda^{1/2}g$ and 
%$\nabla g = (\pd_1g, \ldots, \pd_Ng)$. If we define $A^\pm_1(f_1, f_2)$ by
%\begin{equation}\label{cont:1}\begin{split}
%A^\pm_1[a](f_1, f_2) &= -\int^{\pm\infty}_0
%\CF^{-1}_{\xi'}[a(\xi', x_N+y_N)
%\hat f_{2N}(\xi', y_N)]\,dy_N \\
%&- \int^{\pm\infty}_0\CF^{-1}\Bigl[\frac{(\pd_Na)(\xi', x_N+y_N)
%\rho_{*\pm}\lambda^{1/2}}
%{\mu_{*\pm} B_\pm^2}\hat f_1(\xi', y_N)\Bigr]\,dy_N \\
%&+ \sum_{k=1}^{N-1}\int^{\pm\infty}_0\CF^{-1}_{\xi'}
%\Bigl[\frac{(\pd_Na)(\xi', x_N+y_N)i\xi_k}{B_\pm^2}
%\hat f_{2k}(\cdot, y_N)\Bigr]\,dy_N, 
%\end{split}\end{equation}
%then we have
%\begin{equation}\label{eq:inv1}
%\CF^{-1}_{\xi'}[a(\xi', x_N)\hat g(\xi', 0)] = A^\pm_1[a](\lambda^{1/2}g,
%\nabla g).
%\end{equation}
%Analogously, 
Using the identity: $1 
= \frac{\rho_{\pm}\lambda}{\mu_{\pm} B^2_\pm}
- \sum_{k=1}^{N-1}\frac{(i\xi_k)(i\xi_k)}{B_\pm^2}$, we write
\begin{align*}
a(\xi', x_N)\hat h(\xi', 0) & = -\int^{\pm\infty}_0a(\xi',x_N+y_N)
\Bigl[\frac{\rho_{\pm}\lambda^{1/2}\widehat{\pd_Nh}(\xi', y_N)}
{\mu_{\pm} B^2_\pm}
- \sum_{k=1}^{N-1}\frac{i\xi_k \widehat{\pd_k\pd_Nh}(\xi', y_N)}
{B_\pm}\Bigr]\,dy_N\\
&- \int^{\pm\infty}_0 (\pd_Na)(\xi',x_N+y_N)
\Bigl[\frac{\rho_{\pm}\lambda \hat h(\xi', y_N)}{\mu_{\pm} B_\pm^2}
- \sum_{k=1}^{N-1}
\frac{\widehat{\pd_k\pd_k h}(\xi', y_N)}{B_\pm^2}\Bigr]\,dy_N.
\end{align*} 
Let $f_3$, $f_4 = (f_{41}, \ldots, f_{4N})$ and 
$f_5 = (f_{5JK} \mid J,K=1, \ldots, N)$ be the correspoding variables
to $\lambda h$, $\lambda^{1/2}\nabla h$ and $\nabla^2 h = (\pd_J\pd_Kh
\mid J, K=1, \ldots, N)$.  
If we define $A^\pm_1[a](f_3, f_4, f_5)$ by 
\begin{equation}\label{cont:1}\begin{split}
A^\pm_1[a](f_3, f_4, f_5) &  
= -\int^{\pm\infty}_0\CF^{-1}_{\xi'}\Bigl[a(\xi', x_N+y_N)
\Bigl\{\frac{\rho_{\pm}\hat f_{4N}(\xi', y_N)}{\mu_{\pm} B^2_\pm}
- \sum_{k=1}^{N-1}\frac{i\xi_k \hat f_{3kN}(\xi', y_N)}
{B_\pm}\Bigr\}\Bigr]\,dy_N\\
&- \int^{\pm\infty}_0 \CF^{-1}_{\xi'}\Bigl[(\pd_Na)(\xi', x_N+y_N)
\Bigl\{\frac{\rho_{\pm}\hat f_3(\xi', y_N)}{\mu_{\pm} B_\pm^2}
- \sum_{k=1}^{N-1}
\frac{\hat f_{3kk}(\xi', y_N)}{B_\pm^2}\Bigr\}\Bigr]\,dy_N,
\end{split}\end{equation}
then we have
\begin{equation}\label{eq:inv1}
\CF^{-1}_{\xi'}[a(\xi', x_N)\hat h(\xi', 0)] = A^\pm_1[a](\lambda h,
\lambda^{1/2}\nabla h, \nabla^2h).
\end{equation}

Let us define $u^\pm_{Ji}$ $(i=1,2,3)$ and $p_-$ by 
$u^\pm_{Ji}  =\CF^{-1}_{\xi'}[\hat u^\pm_{Ji}]$ ($i=1,2,3$)
and $p_- = \CF^{-1}_{\xi'}[\hat p_-]$
with $\hat H=0$.     
Setting $u_{J\pm} = \sum_{i=1}^3 u^\pm_{Ji}$,  
by \eqref{eq:sol-form} we see that $\bu_\pm=(u_{1\pm}, \ldots, 
u_{N\pm})$ and $p_-$ satisfy the equations \eqref{eq:prob6}.  
According to the formulas \eqref{cont:1} and \eqref{eq:inv1}, 
we define our solution operators
$\CS^\pm_{Ji}(\lambda)$ ($i = 1,2,3$)
and $\CP_{-2}(\lambda)$ of problem \eqref{eq:prob6} such that 
\begin{alignat}2
u^\pm_{Ji} & 
= \CS^\pm_{Ji}(\lambda)(\lambda \bh, \lambda^{1/2}\nabla \bh, \nabla^2 \bh) 
&\quad &\text{on $\BR^N_{\pm}$}\quad (i = 1,2,3), \nonumber \\
p_- & = \CP_{-2}(\lambda)(\lambda \bh,
 \lambda^{1/2}\nabla \bh, \nabla^2 \bh)
&\quad &\text{on $\BR^N_-$} \label{eq:r-1}
\end{alignat}
with 
$\bh = (h_1, \ldots, h_{N-1})$ as follows:  
Note that 
\begin{equation}\label{cont:3}\begin{split}
\pd_NM_\pm(x_N+y_N) &= \mp(e^{\pm B_\pm(x_N+y_N)} + A_\pm M_\pm(x_N + y_N)),
\quad \pd_Ne^{A(x_N+y_N)} = Ae^{A(x_N+y_N)}, \\
\pd_Ne^{\mp B_\pm(x_N+y_N)} &= \mp B_\pm e^{\mp B_\pm(x_N+y_N)}, 
\end{split}\end{equation}
where we have set $A_- = A$. 
Let 
$F_3= (F_{3m} \mid m =1, \ldots, N-1)$, $F_4 = (F_{4J m}
\mid J=1, \ldots, N,  m=1, \ldots, N-1)$ and $F_5 = (F_{5JKm}\mid
J, K=1, \ldots, N, m=1, \ldots, N-1)$ be the corresponding variables 
to $\lambda\bh = (\lambda h_1, \ldots, \lambda h_{N-1})$, 
$\lambda^{1/2}\nabla\bh 
= (\lambda^{1/2}\pd_J h_m \mid J=1, \ldots, N, m=1, \ldots, N-1)$
and $\nabla^2\bh = (\pd_J\pd_K h_m\mid J, K =1,\ldots, N,
 m=1, \ldots, N-1)$, 
respectively. 
Then, we define 
the operators $\CS^\pm_{J1}(\lambda)$, $\CS^\pm_{J2}(\lambda)$, 
$\CS^\pm_{J3}(\lambda)$, and 
$\CP_{-2}(\lambda)$ by 
\begin{align}
&\CS^\pm_{J1}(\lambda)(F_3, F_4, F_5) = \nonumber\\
&-\int^{\pm\infty}_0\!\!\!\!\! \CF^{-1}_{\xi'}
\Bigl[
AM_\pm(x_N+y_N)
\sum_{m=1}^{N-1}\Bigl(\frac{R^\pm_{Jm,0}\rho_{\pm}\lambda^{1/2}}
{\mu_{\pm}B^2_\pm}\hat F_{4Nm}(\xi', y_N)
-\sum_{k=1}^{N-1}\frac{R^\pm_{Jm,0}(i\xi_k)}{B^2_\pm}
\hat F_{5kNm}(\xi', y_N)\Bigr)\nonumber \\
&\qquad  \mp AM_\pm(x_N+y_N)
\sum_{m=1}^{N-1}\Bigl(
\frac{A_\pm R^\pm_{Jm,0}\rho_{\pm}}{\mu_{\pm}B^2_\pm}\hat F_{3m}(\xi', y_N)
-\sum_{k=1}^{N-1}\frac{A_\pm R^\pm_{Jm,0}}{B^2_\pm}
\hat F_{5kkm}(\xi', y_N)\Bigr)
\nonumber\\
&\qquad \mp Ae^{\mp B_\pm x_N}
\sum_{m=1}^{N-1}\Bigl(
\frac{R^\pm_{Jm,0}\rho_{\pm}}{\mu_{\pm}B^2_\pm}\hat F_{3m}(\xi', y_N)
-\sum_{k=1}^{N-1}\frac{R^\pm_{Jm,0}}{B^2_\pm}
\hat F_{5kkm}(\xi', y_N)\Bigr) 
\Bigr](x')\,dy_N;
\nonumber\\
%%%%%%%%%%%%%%%%%%%%%%%%%%%%%%%
%%%%%%%%%%%%%%%%%%%%%%%%%%%%%%%
&\CS^\pm_{J2}(\lambda)(F_3, F_4, F_5) =\nonumber\\
& -\int^{\pm\infty}_0\!\!\!\!\! 
\CF^{-1}_{\xi'}\Bigl[Ae^{\mp B_\pm(x_N+y_N)}
\sum_{m=1}^{N-1}
\Bigl(\frac{S^\pm_{Jm,-1}\rho_{\pm}\lambda^{1/2}}
{\mu_{\pm}B_\pm^2}\hat F_{4Nm}(\xi', y_N)
-\sum_{k=1}^{N-1}
\frac{S^\pm_{Jm,-1}i\xi_k}{B_\pm^2}\hat F_{5kNm}(\xi', y_N)
\Bigr)\nonumber\\
%&\qquad 
%+\frac{S^\pm_{JN,0}\rho_\pm^2\lambda}{\mu_\pm^2B_\pm^4}
%\hat F_{5N}(\xi', y_N) - \frac{2S^\pm_{JN,0}\rho_\pm A^2}{\mu_\pm B_\pm^4}
%\hat F_{5N}(\xi', y_N)
%+ \frac{S^\pm_{JN,0}A^2}{B_\pm^4}\hat F_{6N}(\xi', y_N)
%\Bigr\}\Bigr](x')\,dy_N
%\nonumber\\
%&\quad \pm \int^{\pm\infty}_0 \CF^{-1}_{\xi'}\Bigl[
\nonumber\\
&\qquad \mp 
Ae^{\mp B_\pm(x_N+y_N)}
\sum_{m=1}^{N-1}\Bigl(\frac{S^\pm_{Jm,-1}\rho_{\pm}}
{\mu_{\pm}B_\pm}\hat F_{3m}(\xi', y_N)
 - \sum_{k=1}^{N-1}\frac{S^\pm_{Jm,-1}}{B_\pm}\hat F_{5kkm}(\xi', y_N)
%\nonumber\\
%&\qquad
%+\frac{S^\pm_{JN,0}\rho_\pm^2\lambda^{1/2}}{\mu_\pm^2B_\pm^4}
%\hat F_4(\xi', y_N)
%-\sum_{k=1}^{N-1}\frac{2S^\pm_{JN,0}\rho_\pm(i\xi_k)}{\mu_\pm
%B_\pm^4}\hat F_{5k}(\xi', y_N)\nonumber\\
%&\phantom{\qquad
%+\frac{S^\pm_{JN,0}\rho_\pm^2\lambda^{1/2}}{\mu_\pm^2B_\pm^4}
%\hat F_4(\xi', y_N)}\,\,
%-\sum_{k=1}^{N-1}\frac{S^\pm_{JN,0}(i\xi_k)}{B^4_\pm}
%\hat F_{6k}(\xi', y_N)
\Bigr)\Bigr](x')\,dy_N;
\nonumber\\
%%%%%%%%%%%%%%%%%
%%%%%%%%%%%%%%%%%
&\CS^\pm_{j5}(\lambda)(F_3, F_4, F_5) = \nonumber\\
&\quad -\int^{\pm\infty}_0\CF^{-1}_{\xi'}\Bigl[
e^{\mp B_\pm(x_N+y_N)}\Bigl\{
\frac{T^\pm_{j,0}\rho_{\pm}}{\mu_{\pm}B_\pm^2}
\hat F_{4Nj}(\xi', y_N) - \sum_{k=1}^{N-1}
\frac{T_{j,0}^\pm i\xi_k}{B_\pm^2}\hat F_{5kNj}(\xi', y_N)
\nonumber\\
&\quad \phantom{\int^{\pm\infty}_0\CF^{-1}_{\xi'}\Bigl[e^{\mp B_\pm(x_N+y_N)}
\Bigl\{}
\pm\Bigl(\frac{T^\pm_{j,0}\rho_{\pm}}{\mu_{\pm}B_\pm}
\hat F_{3j}(\xi', y_N)
-\sum_{k=1}^{N-1}\frac{T^\pm_{j,0}}{B_\pm}
\hat F_{5kkj}(\xi', y_N)\Bigr)\Bigr\}\Bigr](x')\,dy_N; \nonumber \\
&\CS^\pm_{N3}(\lambda)(F_3, F_4, F_5) = 0 \nonumber\\
%%%%%%%%%%%%%%%
%%%%%%%%%%%%%%%%%%%%
&\CP_{-2}(\lambda)(F_3, F_4, F_5) =\nonumber\\
&\quad\int^0_{-\infty}\Bigl[e^{A(x_N+y_N)}\Bigl\{
\sum_{m=1}^{N-1}\Bigl(
\frac{p^-_{m,1}\rho_{-}\lambda^{1/2}}{\mu_{-}B_-^2}
\hat F_{4Nm}
-\sum_{k=1}^{N-1}\frac{p_{m,1}i\xi_k}{B_-^2}\hat F_{5kNm}(\xi', y_N)
\Bigr) \nonumber\\
%&\qquad +\frac{p^-_{N,2}\rho_-^2\lambda}{\mu_-B_-^4}
%\hat F_{5N}(\xi', y_N)
%-\frac{2 p^-_{N,2}\rho_-A^2}{\mu_-B_-^4}\hat F_{5N}(\xi', y_N)
%+\frac{p^-_{N,2}A^2}{B_-^4}\hat F_{6N}(\xi', y_N)
%\Bigr\}\Bigr](x')\,dy_N \nonumber\\
&\quad \phantom{+ \int^0_{-\infty}\CF^{-1}_{\xi'}\Bigl[aaa}
+\sum_{m=1}^{N-1}\Bigl(\frac{Ap^-_{m,1}\rho_{-}}{\mu_{-}B_-^2}
\hat F_{3m}(\xi', y_N)
-\sum_{k=1}^{N-1}\frac{Ap^-_{m,1}}{B_-^2}\hat F_{5kkm}(\xi', y_N)
\Bigr)
%\nonumber\\
%&\qquad
%+\frac{p^-_{N,2}\rho_-\lambda^{1/2}}{\mu_-B_-^4}\hat F_4(\xi', y_N)
%-\sum_{k=1}^{N-1}\frac{2p^-_{N,2}\rho_-(i\xi_k)}
%{\mu_-B_-^4}\hat F_{5k}(\xi', y_N) \nonumber\\
%&\phantom{\qquad
%+\frac{p^-_{N,2}\rho_-\lambda^{1/2}}{\mu_-B_-^4}\hat F_4(\xi', y_N)}\,\,
%-\sum_{k=1}^{N-1}\frac{p^-_{N,2}(i\xi_k)}{B_-^4}
%\hat F_{6k}(\xi', y_N)
\Bigr\}\Bigr](x')\,dy_N.
 \label{sol-form:1}
\end{align}
Obviously, by \eqref{eq:inv1}, we have  
\eqref{eq:r-1}.

%%%%%%%%%%%%%%%%%%%%
%%%%%%%%%%%%%%%%
%%%%%%%%%
%%%%%%%%%%%%%%%%%%%%%
%%%%%%%%%%%%%%%%%%%%%%

Given that   operators $\CA_{\pm2}(\lambda)$
are defined by 
$\CA_{\pm2}(\lambda)\bF' 
= \sum_{i=1}^3(\CS^\pm_{1i}(\lambda)\bF',
\ldots, \CS^\pm_{Ni}(\lambda)\bF')
$
with $\bF' = (F_3, F_4, F_5)$,
by \eqref{eq:r-1}  we have 
\begin{equation}\label{eq:r-1*}
\bu_\pm = \CA_{\pm2}(\lambda)(\lambda\bh, 
\lambda^{1/2}\nabla\bh, \nabla^2\bh),
\quad p_- = \CP_{-2}(\lambda)(\lambda\bh, 
\lambda^{1/2}\nabla\bh, \nabla^2\bh).
\end{equation}
Moreover, we have Theorem \ref{thm:au} with
the help of the following two lemmas:
%Moreover, if we set 
%\begin{align*}
%\CZ_q(\BR^N) = \{(F_3, F_4, F_5) \mid 
%&F_3 \in L_q(\BR^N)^{N-1}, \enskip F_4 \in L_q(\BR^N)^{N(N-1)},
%\enskip F_5
%\in L_q(\BR^N)^{(N-1)N^2}\},
%\end{align*}
%then, using Lemma \ref{lem:tech:comp} 
%and Lemma \ref{lem:tech:incomp}
%below,  we have
%\begin{equation}\label{eq:r-2}\begin{split}
%&\CR_{\CL(\CZ(\BR^N), L_q(\BR^N_{\pm})^{2N+N^2+N^3})}
%(\{(\tau\pd_\tau)^\ell G_\lambda^1 \CS_{\pm J}(\lambda) \mid
%\lambda \in \Lambda_{\epsilon, \lambda_0}\}) \leq C
%\quad(\ell = 0, 1), \\
%&\CR_{\CL(\CZ(\BR^N), L_q(\BR^N_{\pm})^N)}
%(\{(\tau\pd_\tau)^\ell \nabla\CP_-(\lambda) \mid
%\lambda \in \Lambda_{\epsilon, \lambda_0}\}) \leq C
%\quad(\ell = 0, 1).
%\end{split}\end{equation}
\begin{lem}\label{lem:tech:comp}
Let $1 < q < \infty$ and 
let $n_1^+$, $n_2^+$  and $n_3^+$ be multipliers belonging to $\bM_{-1,2}(0)$, 
$\bM_{-2,2}(0)$ and $\bM_{-1,1}(0)$, respectively.   
Let $K_i^+$ $(i = 1, 2, 3)$ 
be operators defined by 
\begin{align*}
K^+_1(\lambda)g & = \int^\infty_0\CF^{-1}_{\xi'}
[n_1^+(\lambda, \xi')AM_+(x_N+y_N)\hat g(\xi', y_N)](x')\,dy_N, \\
K^+_2(\lambda)g & = \int^\infty_0\CF^{-1}_{\xi'}
[n_2^+(\lambda, \xi')Ae^{-B_+(x_N+y_N)}\hat g(\xi', y_N)](x')\,dy_N, \\
K^+_3(\lambda)g & = \int^\infty_0\CF^{-1}_{\xi'}
[n_3^+(\lambda, \xi')e^{-B_+(x_N+y_N)}\hat g(\xi', y_N)](x')\,dy_N.
\end{align*}
Then, there exists a constant $C$ such that 
$$\CR_{\CL(L_q(\HSp), L_q(\HSp)^{1 + N+N^2})}
(\{(\tau\pd_\tau)^\ell G_\lambda^1 K_i^+(\lambda) \mid
\lambda \in \Lambda\}) \leq C
\quad(\ell = 0, 1, \enskip i = 1, 2, 3).
$$ 
\end{lem}
\begin{lem}\label{lem:tech:incomp} Let $1 < q < \infty$ and 
let $n_1^-$, $n_2^-$, $n_3^-$ and $n_4^-$ be multipliers belonging to
$\bM_{-1,2}(0)$, $\bM_{-2,2}(0)$, $\bM_{-1,1}(0)$ and $\bM_{0,2}(0)$, 
respectively. 
Let $K_i^+$ $(i = 1, 2, 3, 4)$ 
be operators defined by 
\begin{align*}
K^-_1(\lambda)g & = \int^0_{-\infty}\CF^{-1}_{\xi'}
[n_1^-(\lambda, \xi')AM_-(x_N+y_N)\hat g(\xi', y_N)](x')\,dy_N, \\
K^-_2(\lambda)g & = \int^0_{-\infty}\CF^{-1}_{\xi'}
[n_2^-(\lambda, \xi')Ae^{B_-(x_N+y_N)}\hat g(\xi', y_N)](x')\,dy_N, \\
K^-_3(\lambda)g & = \int^0_{-\infty}\CF^{-1}_{\xi'}
[n_3^-(\lambda, \xi')e^{B_-(x_N+y_N)}\hat g(\xi', y_N)](x')\,dy_N, \\
K^-_4(\lambda)g & = \int^0_{-\infty}\CF^{-1}_{\xi'}
[n_4^-(\lambda, \xi')e^{A(x_N+y_N)}\hat g(\xi', y_N)](x')\,dy_N.
\end{align*}
Then, there exists a constant $C$ such that 
\begin{equation}\label{eq:incomp-1}\begin{split}
&\CR_{\CL(L_q(\HSn), L_q(\HSn)^{1 + N + N^2})}
(\{(\tau\pd_\tau)^\ell G_\lambda K_i^-(\lambda) \mid
\lambda \in \Lambda\}) \leq C
\quad(\ell = 0, 1, \enskip i = 1, 2, 3), \\
&\CR_{\CL(L_q(\HSn), L_q(\HSn)^N)}
(\{(\tau\pd_\tau)^\ell \nabla K_4^-(\lambda) \mid
\lambda \in \Lambda\}) \leq C
\quad(\ell = 0, 1).
\end{split}\end{equation}
\end{lem}
\begin{remark}
Lemma \ref{lem:tech:comp} was proved in \cite[Sect.2]{KSS} and
Lemme \ref{lem:tech:incomp} was proved in \cite{SS4}.
\end{remark}

%%%%%%%%%%%%%%%%%%%%%
%%%%%%%%%%%%%%%%%%
\section{Some estimates of several multipliers}

In this section, we prove \eqref{sym-est:1} and 
\eqref{sym-est:6}.  
For this purpose, we use the following well-known estimate:
\begin{equation}\label{lem:3.1}
|\alpha\lambda + \beta| \geq (\sin\frac{\epsilon}{2})
(\alpha|\lambda| + \beta)
\end{equation}
for any $\lambda \in \Sigma_{\epsilon}$ and positive numbers
$\alpha$ and $\beta$. 
%\begin{proof}
% \thetag1 is well-known. In  case of \thetag{C1}, that is 
% $\delta = \gamma_{1+}\gamma_{2+}
% \lambda^{-1}$, 
%\thetag2 was proved in Shibata and 
%Tanaka \cite[Lemma 2.1]{ST}. In  cases of \thetag{C2} and \thetag{C3},
%we set $\Gamma^0_{\epsilon, \lambda_0} = \{ \lambda \in 
%\Gamma_{\epsilon, \lambda_0} \mid|\lambda|=\lambda_0\}$.  For
%$\lambda \in \Gamma^0_{\epsilon, \lambda_0}$, we see that 
%$(s\mu_+ + \nu_+ +\delta)^{-1}\lambda \not\in(-\infty, 0]$.
%In fact, assuming that $(s\mu_+ + \nu_+ +\delta)^{-1}\lambda
%=-\gamma^2$ with some $\gamma \in \BR\setminus\{0\}$, we have
%$\re\lambda=-\gamma^2(s\mu_++\nu_++\re\delta)$ and 
%$\im\lambda = -\gamma^2\im \delta$.  In case of \thetag{C2}, 
%it follows from $\re\lambda\geq|\re\delta/\im\delta||\im\lambda|$
%that $-\gamma^2(s\mu_++\nu_+) \geq 0$, which implies that 
%$s\mu_+ + \nu_+ \leq
% 0$.  This contradicts to 
%$s\mu_++\nu_+ > 0$.  In case of \thetag{C3}, it follows from 
%$\re\delta \geq 0$ that $\re\lambda < 0$, which 
%also contradicts to $\re\lambda \geq \lambda_0|\im\lambda| \geq 0$. 
%Since $\Gamma^0_{\epsilon, \lambda_0}$ is a compact set in
%$\BC$, there exists a $\sigma \in (0, \pi)$ such that 
%$|\arg \lambda| \leq \pi-\sigma$ for any 
%$\lambda \in \Gamma_{\epsilon, \lambda_0}$, from
%which we conclude that the assertion \thetag2 holds in cases of
%\thetag{C2} and \thetag{C3}.
%The assertion \thetag3 follows from \thetag1, \thetag2 and 
%assumptions: $|\lambda| \geq \lambda_0$ and $|\delta| 
%\leq \max(\delta_0, \gamma_{1+}\gamma_{2+}\lambda_0^{-1})$ (cf \eqref{est:0}). %\end{proof}

First we estimate $A_+^s$, $B_\pm^s$, $(A_+ + B_+)^s$
and $(\mu_{+}B_+ + \mu_{-}B_-)^s$.  
For this purpose,  we use the estimates:
\begin{equation}\label{5.1}
c(|\lambda|^{1/2} + A) \leq {\rm Re}\, M_1 \leq |M_1| \leq c'
(|\lambda|^{1/2} + A)\quad
(M_1 = A_+, \enskip B_\pm)
\end{equation} 
for any $(\lambda, \xi') \in \tilde \Sigma_{\epsilon}
=\Sigma_{\epsilon}\times(\BR^{N-1}\setminus\{0\})$
with some positive constants $c$ and $c'$, which immediately follows from
\eqref{lem:3.1}. Here and in the sequel, 
$c$ and $c'$ denote some positive 
constants essentially depending  on $\mu_{\pm}$, $\nu_{+}$, 
$\rho_{\pm}$ and $\epsilon$.  In particular, by \eqref{5.1}
we have
\begin{equation}\label{5.2}
c(|\lambda|^{1/2} + A) \leq {\rm Re}\, M_2 
\leq |M_2| \leq c'(|\lambda|^{1/2} + A)
\quad(M_2 = A_++B_+,\,\, \mu_{+}B_++\mu_-B_{-})
\end{equation}
for any $(\lambda, \xi') \in \tilde\Sigma_{\epsilon, 0}$. 
As was seen in Enomoto and Shibata 
\cite[Lemma 4.3]{ES}, using \eqref{5.1}, \eqref{5.2} and 
the Bell formula: 
\begin{equation}\label{5.3}
\pd_{\xi'}^{\kappa'}f(g(\xi'))
= \sum_{\ell=1}^{|\kappa'|}f^{(\ell)}(g(\xi'))
\sum_{\kappa'_1 + \cdots + \kappa_\ell' = \kappa' \atop |\kappa'_i| \geq 1}
\Gamma^{\kappa'}_{\kappa'_1, \ldots, \kappa'_\ell}
(\pd_{\xi'}^{\kappa_1'}g(\xi'))\cdots(\pd_{\xi'}^{\kappa'_\ell}g(\xi'))
\end{equation}
with  suitable coefficients 
$\Gamma^{\kappa'}_{\kappa'_1, \ldots, \kappa'_\ell}$,
where $f^{(\ell)}(t) = d^\ell f(t)/dt^\ell$, 
we see that 
\begin{equation}\label{5.4}
(M_3)^s \in \bM_{s, 1}(0) \quad(M_3 = A_+,\enskip B_+, \enskip
A_++B_+, \enskip \mu_{+}B_++\mu_{-}B_-). 
\end{equation}

Second, we estimate $(A_+B_+ - A^2)^{-1}$.  For this purpose, we write
\begin{equation}\label{5.5}
\frac{1}{A_+B_+ - A^2} = \frac{(\mu_{+}+\nu_{+})\mu_{+}}
{\rho_{+}(2\mu_{+} + \nu_{+})\lambda}P(\lambda, \xi')
\quad\text{with $P(\lambda, \xi') = \frac{A_+B_++A^2}
{\rho_{+}(2\mu_{+} + \nu_{+})^{-1}\lambda + A^2}$}.
\end{equation}
Noting that $A^2 \in \bM_{2,1}(0)$, 
by \eqref{lem:3.1}, \eqref{5.3} and \eqref{5.4} we have
\begin{equation}\label{5.6}\begin{split}
A_+B_+ + A^2 \in \bM_{2,1}(0), \enskip 
(\rho_{+}(2\mu_{+}+\nu_{+})^{-1}\lambda+ A^2)^s
\in \bM_{2s, 1}(0),
\end{split}\end{equation}
so that by Lemma \ref{lem:sym} we have 
\begin{equation}\label{5.7}
P \in \bM_{0,1}(0).
\end{equation}
Since $A^2 - A_+^2 = \rho_{+}(\mu_{+}+\nu_{+})^{-1}\lambda$,
by \eqref{5.5} and \eqref{5.7}, we have 
$\frac{A^2-A_+^2}{A_+B_+-A^2}\in \bM_{0,1}(0)$, 
which, combined with \eqref{5.4}, furnishes \eqref{sym-est:6}. 
 
Applying \eqref{5.5} to the formula in \eqref{eq:eq5}, we have
\begin{equation}\label{5.8}\begin{split}
L^+_{11} & = \frac{\mu_{+}(\mu_{+} + \nu_{+})}
{2\mu_{+}+\nu_{+}}A_+P, \quad 
L^+_{12}  =\mu_{+}A^2\Bigl(2 - \frac{\mu_{+}+\nu_{+}}
{2\mu_{+} + \nu_{+}}P\Bigr), \\
L^+_{21} & = \Bigl(\frac{2\mu_{+}\nu_{+}}{2\mu_{+}+\nu_{+}}
\frac{A_+}{B_++A_+} - \frac{\mu_{+}(\nu_{+}-\mu_{+})}
{2\mu_{+}+\nu_{+}}\Bigr)P, \quad
L_{22}^+  = \frac{\mu_{+}(\mu_{+}+\nu_{+})}
{2\mu_{+}+\nu_{+}}B_+P.
\end{split}\end{equation}
By Lemma \ref{lem:sym}, \eqref{eq:eq5},
\eqref{5.4}, \eqref{5.7} and \eqref{5.8}, we have 
$L_{11}^+ \in \bM_{1,1}(0)$, $L^+_{12} \in \bM_{2,1}(0)$, 
$L_{21}^+ \in \bM_{0,1}(0)$ and $L_{22}^+ \in 
\bM_{1,1}(0)$. In addition, since $A \in \bM_{1,2}(0)$
and $B_- \in \bM_{1,2}(0)$, by Lemma \ref{lem:sym} we have
$A\pm B_- \in \bM_{1,2}(0)$ and $(A+B_-)B_- \in \bM_{2,2}(0)$.
Summing up, we have proved \eqref{sym-est:1}. 
%%%%%%%%%%
%%%%%%%%%%%%%
\section{Analysis of Lopatinski determinant}
In this section, we show the following lemma which implies 
\eqref{sym-est:3}. 
\begin{lem}\label{lem:lop}
Let $L$ be the matrix defined in \eqref{eq:eq5*}.  Then,
there exists a positive constant $\omega$ depending solely on 
$\mu_\pm$, $\nu_+$, $\rho_{\pm}$, and $\epsilon$ such that 
\begin{equation}\label{lop:1} |\det L| \geq \omega(|\lambda|^{1/2} + A)^4
\end{equation}
for any $(\lambda, \xi') \in \tilde\Sigma_{\epsilon, 0}$.

Moreover, we have
\begin{equation}\label{lop:2}
|\pd_{\xi'}^{\kappa'}\{(\tau\pd_\tau)^\ell(\det L)^{-1}\}|
\leq C_{\kappa'}(|\lambda|^{1/2} + A)^{-4}A^{-|\kappa'|} \quad(\ell = 0, 1)
\end{equation}
for any multi-index $\kappa' \in \BN_0^{N-1}$ and 
$(\lambda, \xi') \in \tilde\Sigma_{\epsilon, 0}$.
Namely, $(\det L)^{-1} \in \bM_{-4,2}(0)$. 
\end{lem}
\begin{proof}
We see that 
\begin{equation}\label{det:1}
\det L = L_{22}^-\det L^+ +
L_{22}^+\det L^-
\end{equation}
with $L^\pm= \det \left(\begin{matrix}
L_{11}^\pm & L_{12}^\pm \\ L_{21}^\pm & L_{22}^\pm \end{matrix}\right)$. 
To prove \eqref{lop:1}, first we consider the case: $R_1|\lambda|^{1/2} \leq A$
with large $R_1 \geq 1$. Let $P$ be the function defined in \eqref{5.5}.
By \eqref{5.5} we see easily that $P = 2 + O(\delta_1)$, 
that $A_+ =A(1 + O(\delta_1))$ 
and  that $B_\pm = A(1 + O(\delta_1))$ when 
$|\rho_{+}(\mu_++\nu_+)^{-1}\lambda A^{-2}| 
\leq \rho_{+}(\mu_++\nu_+)^{-1}R_1^{-2}
\leq \delta_1^2$ and 
$|\rho_{\pm}\mu_\pm^{-1}\lambda A^{-2}| \leq 
\rho_{\pm}\mu_\pm^{-1}R_1^{-2} \leq \delta_1^2$ with 
very small positive number $\delta_1$. 
Thus, by \eqref{5.8} we have 
\begin{alignat}2
L^+_{11} & = \frac{2\mu_+(\mu_++\nu_+)}
{2\mu_++\nu_+}A(1 + O(\delta_1)), &\quad
L_{12}^+ & = \frac{2(\mu_+)^2}{2\mu_++\nu_+}
A^2(1 + O(\delta_1)), \nonumber \\
L^+_{21} & = \frac{2(\mu_+)^2}{2\mu_++\nu_+}(1 + O(\delta_1)),
&\quad 
L^+_{22} & = \frac{2\mu_+(\mu_++\nu_+)}
{2\mu_++\nu_+}A(1 + O(\delta_1)). \label{lop:0.1}
\end{alignat}
On the other hand, we have $B_--A = 
\frac{\gamma_{0-}\lambda}{\mu_{0-}(B_-+A)} = AO(\delta_1)$,
so that by \eqref{eq:eq5} we have
\begin{alignat}2
L^-_{11} & =2\mu_-A(1 + O(\delta_1)), &\quad 
L^-_{12} & = A^2O(\delta_1), \nonumber \\
L_{21}^- & = AO(\delta_1), &\quad L^-_{22} & = 2\mu_-A^2(1 + O(\delta_1)).
\label{lop:0.2}
\end{alignat}
Thus,  by \eqref{det:1} we have
\begin{equation}
\det L = \omega_1
A^4(1+O(\delta_1)) \quad\text{with 
$\omega_1 = \frac{8\mu_+\mu_-(\mu_+\nu_++\mu_-(\mu_++\nu_+))}{2\mu_++\nu_+}$}
\label{lop:1.1*}
\end{equation}
so that we can choose $R_1 \geq 1$ so large that 
\begin{equation}\label{lop:1.1}
|\det L| \geq \frac12\omega_1(|\lambda|^{1/2} + A)^4
\end{equation}
for any $(\lambda, \xi') \in \tilde\Sigma_{\epsilon,0}$ with
$R_1|\lambda|^{1/2}\leq A$.

Secondly, we consider the case: $R_2A \leq |\lambda|^{1/2}$ 
with large $R_2 \geq 1$.  In this case, we have 
\begin{equation}
A_+ = (\mu_++\nu_++\delta)^{-1/2}(\gamma_{0+}\lambda)^{1/2}(1 + 
O(\delta_2)), \quad
B_\pm = (\mu_\pm)^{-1/2}(\gamma_{0\pm}\lambda)^{1/2}(1 + O(\delta_2))
\label{lop:0.4}
\end{equation}
when $|(\mu_++\nu_+)(\rho_{+}\lambda)^{-1}A^2| \leq 
(\mu_++\nu_++\delta_1)\rho_+^{-1}R^{-2}_2\leq\delta_2^2$ 
and $|\mu_\pm(\rho_{\pm}\lambda)^{-1}A^2| 
\leq \mu\rho_{\pm}^{-1}R_2^{-2} \leq
\delta_2^2$ with some very small positive number $\delta_2$. 
By \eqref{eq:eq5} 
\begin{alignat}2
L_{11}^+ & = (\mu_+\rho_+\lambda)^{1/2}(1 + O(\delta_2))
&\quad
L_{12}^+ &= O(\delta_2)\lambda \nonumber \\
L_{21}^+ & = \frac{2\mu_+((\mu_++\nu_+)^{1/2}-\mu_+^{1/2})
-(\nu_+-\mu_+)\mu_+^{1/2}}{(\mu_++\nu_+)^{1/2}}
(1 + O(\delta_2)),  &\quad
L_{22}^+ & = ((\mu_++\nu_+)\rho_+\lambda)^{1/2}(1+O(\delta_2))
\nonumber \\
L_{11}^- & = (\mu_-\rho_-\lambda)^{1/2}(1+O(\delta_2))
&\quad
L_{12}^- &= O(\delta_2)\lambda \nonumber \\
L_{21}^- & = (\mu_-\rho_-\lambda)^{1/2}(1+O(\delta_2)), &\quad
L_{22}^- & = \rho_+\lambda(1+O(\delta_2))
\label{lop:0.3}
\end{alignat}
Thus, by \eqref{det:1} we have
\begin{equation}\label{lop:1.2*}
|\det L| = \omega_2
|\lambda|^2(1 + O(\delta_2))
\quad\text{with $\omega_2 = 
\mu_+^{1/2}(\mu_++\nu_+)^{1/2}
\rho_+\rho_- + 
\mu_-^{1/2}(\mu_++\nu_+)^{1/2}\rho_+^{1/2}\rho_-^{3/2}$}, 
\end{equation}
so that  we can choose $R_2 \geq 1$ so large that 
\begin{equation}\label{lop:1.2}
|\det L| \geq \frac12\omega_2(|\lambda|^{1/2} + A)^4
\end{equation}
for any $(\lambda, \xi') \in \tilde\Sigma_{\epsilon,0}$ with
$R_2A\leq |\lambda|^{1/2}$.

%%%%%%%%%%
Thirdly, we consider the case: $R_2^{-1}|\lambda|^{1/2} 
\leq A \leq R_1|\lambda|^{1/2}$.  Set 
\begin{align*}
\tilde\lambda &= \frac{\lambda}{(|\lambda|^{1/2} + A)^2}, \quad 
\tilde A = \frac{A}{|\lambda|^{1/2} + A}, \\
\tilde A_+ &= \sqrt{\rho_+(\mu_++\nu_+)^{-1}
\tilde\lambda + \tilde A^2}, \quad
\tilde B_\pm = \sqrt{\rho_\pm(\mu_\pm)^{-1}
\tilde\lambda + \tilde A^2}, \\
D_\epsilon(R_1, R_2) &= \{(\tilde\lambda, \tilde A) \mid 
(1 + R_1)^{-2}\leq|\tilde\lambda| \leq R_2^2
(1+R_2)^2,
 \enskip
(1 + R_2)^{-1} \leq \tilde A \leq R_1(1+R_1)^{-1},
\enskip \tilde\lambda \in \Sigma_\epsilon\}.
\end{align*}
If $(\lambda, \xi')$ satisfies the condition: 
$R_2^{-1}|\lambda|^{1/2} \leq A \leq R_1|\lambda|^{1/2}$ and
$\lambda \in \Sigma_\epsilon$, then 
$(\tilde\lambda, \tilde A) \in D_\epsilon(R_1, R_2)$. Note that 
$\lambda\not=0$ when 
$(\tilde\lambda, \tilde A) \in D_\epsilon(R_1, R_2)$.
We define $\tilde L_{ij}$ by replacing $A_+$, $A$ and $B_\pm$ by $\tilde A_+$,
$\tilde A$ and $\tilde B_\pm$ in \eqref{eq:eq5}, respectively. 
And the matrix $\tilde L$ is defined by replacing $L^\pm_{ij}$ 
by $\tilde L^\pm_{ij}$ in \eqref{eq:eq5*}. 
Setting $\det \tilde L = \tilde L_{11}\tilde L_{22} - 
\tilde A\tilde L_{12}\tilde L_{21}$, we have
\begin{equation}\label{hom:1}
\det L =(|\lambda|^{1/2} + A)^4\det \tilde L.
\end{equation}
We prove that $\det \tilde L\not= 0$ provided that 
$(\tilde\lambda, \tilde A) \in D_\epsilon(R_1, R_2)$ 
by contradiction. Suppose that 
$\det \tilde L = 0$.  By \eqref{hom:1} 
$\det L = 0$, so that in view of 
\eqref{eq:eq5*} there exist $w_{\pm J}$ and $p_-$ of 
the forms: $w_{\pm J}(x_N) 
= \alpha_{\pm J}(e^{\mp B_\pm x_N} - e^{\mp A_\pm x_N}) 
+ \beta_{\pm J}e^{\mp B_\pm x_N}$ and 
$p_-(x_N) = \gamma_-e^{Ax_N}$ with $A_- = A$ such that 
$\bw_\pm(x_N) 
= (w_{\pm1}(x_N), \ldots, w_{\pm N}(x_N)) \not=(0, \ldots, 0)$, and
$\bw_{\pm}(x_N)$ and $p_-(x_N)$ satisfy \eqref{eq:7}
and \eqref{eq:bdy1} with $\hat h_m(0) =0$ and $\hat H(0) = 0$,  
that is they satisfy the following homogeneous equations:
\allowdisplaybreaks{
\begin{alignat}2
&\rho_+\lambda w_{+j} - \sum_{m=1}^{N-1}\mu_+i\xi_m
(i\xi_j w_{+m} + i\xi_m w_{+j}) 
- \mu_+\pd_N(i\xi_j w_{+N} + \pd_N w_{+j}) \nonumber\\ 
&\quad - 
(\nu_+-\mu_+)i\xi_j(i\xi'\cdot w'_+ + \pd_Nw_{+N})=0
&\quad&\text{for $x_N > 0$}, \nonumber \\
&\rho_+\lambda w_{+N}
 - \sum_{m=1}^{N-1}\mu_+i\xi_m
(\pd_N w_{+m} + i\xi_m w_{+N}) 
- 2\mu_+\pd_N^2 w_{+N}\nonumber\\ 
&\quad - 
(\nu_+- \mu_{+})\pd_N(i\xi'\cdot w'_+ + \pd_Nw_{+N})=0
&\quad&\text{for $x_N > 0$},\nonumber\\
%%%%%%%%%%
&\rho_-\lambda w_{-j} - \sum_{m=1}^{N-1}\mu_-i\xi_m
(i\xi_j w_{-m} + i\xi_m w_{-j}) 
- \mu_-\pd_N(i\xi_j w_{-N} + \pd_N w_{j}) 
+ i\xi_jp_-=0
&\quad&\text{for $x_N < 0$}, \nonumber\\
&\rho_-\lambda w_{-N} - \sum_{m=1}^{N-1}\mu_-i\xi_m
(\pd_N w_{-m} + i\xi_m w_{-N}) 
- 2\mu_-\pd_N^2 w_{-N} + \pd_Np_- = 0 
&\quad&\text{for $x_N < 0$}, \nonumber\\
&i\xi'\cdot w'_- + \pd_Nw_{-N} = 0 
&\quad&\text{for $x_N < 0$}, \nonumber\\
&\mu_-(\pd_Nw_{-j} + i\xi_jw_{-})|_-
-\mu_+(\pd_Nw_{+j} + i\xi_jw_{+N})|_+ = 0,\nonumber\\
&(2\mu_-\pd_Nw_{-N} - p_-)|_- -
(2\mu_+\pd_Nw_{+N} + (\nu_+-\mu_+)(i\xi'\cdot w_+'+
\pd_Nw_{+N})|_+= 0. 
\label{homo:0}
\end{alignat}
}
Set $(a, b)_\pm =\pm\int^{\pm\infty}_0a(x_N)\overline{b(x_N)}\,dx_N$ and 
$\|a\|_\pm = (a, a)_\pm^{1/2}$. Multiplying the equations in \eqref{homo:0}
by $\overline{w_{\pm J}}$ and using the integration by parts
 and the jamp conditions
in \eqref{homo:0}, we have
\allowdisplaybreaks{
\begin{align}
0 & = \lambda(\rho_+\sum_{m=1}^{N}\|w_{+m}\|^2
+ \rho_-\sum_{m=1}^N\|w_{-m}\|^2
) 
+ \mu_+[\sum_{j,k=1}^{N-1} \|i\xi_kw_{+j}\|_+^2 
+ \|i\xi'\cdot w_+'\|_+^2 +
\sum_{j=1}^{N-1}\|\pd_Nw_{+j}\|_+^2\nonumber \\
&+ \sum_{j=1}^{N-1}(i\xi_jw_{+N}, \pd_Nw_{+N})_+
+ \sum_{j=1}^{N-1}\|i\xi_jw_{+N}\|^2_+
+ \sum_{j=1}^{N-1}(\pd_Nw_{+j}, i\xi_jw_{+N})_+ + 2\|\pd_Nw_{+N}\|_+^2]
\nonumber\\
& + (\nu_+-\mu_+)[\|i\xi'\cdot w'_+\|_+^2 
+ (\pd_Nw_{+N}, i\xi'\cdot w'_+)_+
+ (i\xi'\cdot w'_+, \pd_Nw_{+N})_+ + \|\pd_Nw_{+N}\|_+^2] \nonumber \\
%%%%%%%%%%%%%%%%%%%%%%%
&+ \mu_-[\sum_{j,k=1}^{N-1} \|i\xi_kw_{-j}\|_-^2 
+ \|i\xi'\cdot w'_-\|_-^2 +
\sum_{j=1}^{N-1}\|\pd_Nw_{-j}\|_-^2\nonumber \\
&+ \sum_{j=1}^{N-1}(i\xi_jw_{-N}, \pd_Nw_{-N})_-
+ \sum_{j=1}^{N-1}\|i\xi_jw_{-N}\|^2_-
+ \sum_{j=1}^{N-1}(\pd_Nw_{-j}, i\xi_jw_{-N})_- + 2\|\pd_Nw_{-N}\|_-^2]
\nonumber\\
& = \lambda(\gamma_{0+}\|w_+\|_+^2 + \gamma_{0-}\|w_-\|_-^2) 
+ \mu_+[\sum_{j,k=1}^{N-1} \|i\xi_kw_{+j}\|_+^2 
+ \|i\xi'\cdot w_+'\|_+^2 +
\sum_{j=1}^{N-1}\|\pd_Nw_{+j} + i\xi_jw_{+N}\|_+^2\nonumber \\
& + 2\|\pd_Nw_{+N}\|_+^2] + (\nu_+-\mu_+)
\|\pd_Nw_{+N} + i\xi'\cdot w'_+\|_+^2
\nonumber\\
&+ \mu_-[\sum_{j,k=1}^{N-1} \|i\xi_kw_{-j}\|_-^2 
+ \|i\xi'\cdot w_-'\|_-^2 +
\sum_{j=1}^{N-1}\|\pd_Nw_{-j} + i\xi_jw_{-N}\|_-^2
 + 2\|\pd_Nw_{-N}\|_-^2]. \label{homo:1}
\end{align}
}
Taking the real part and the imaginary part in 
\eqref{homo:1}, using the inequality:
$$
\sum_{j,k=1}^{N-1}\|i\xi_jw_{+k}\|_+^2 + \|i\xi'\cdot w'_+\|_+^2 
+ 2\|\pd_Nw_{+N}\|_+^2 
\geq 2(\|i\xi'\cdot w'_+\|_+^2 + \|\pd_Nw_{+N}\|_+^2)
\geq \|\pd_Nw_{+N} + i\xi'\cdot w'_+\|_+^2,
$$
and setting $K = \rho_+\|\bw_+\|_+^2 + \rho_-\|\bw_-\|_-^2$
and $L = \|\pd_Nw_{+ N} + i\xi'\cdot w'_+\|_+^2$ for short,
we have
\begin{equation}\label{new:1}
(\im \lambda)K=0, \quad
0 \geq (\re \lambda)K + \nu_+L.
\end{equation}
When ${\rm Im}\,\lambda\not=0$, obviously $\bw_\pm=0$ which 
leads to a contradiction.  When ${\rm Im}\,\lambda = 0$, 
$\lambda \geq 0$ and $\lambda\not=0$, because $\lambda \in \Sigma_\epsilon$.
By \eqref{new:1} we have $L=0$, because $\nu_+ > 0$. 
Thus, by \eqref{homo:1} we have
\begin{equation}\label{new:2*}
\pd_Nw_{\pm N} = 0,  \quad \pd_Nw_{\pm j} + i\xi_jw_{\pm N} = 0
\quad\text{on $\BR^\pm$}
\end{equation}
where $\BR^+ = (0, \infty)$ and $\BR^- = (-\infty, 0)$. 
By the first equation of
\eqref{new:2*}, $w_{\pm N}(x_N)$ are constants on $\BR^\pm$, 
but $w_{\pm N}(x_N) \to 0$ as $\pm x_N \to \infty$, so that
$w_{\pm N} = 0$.  Thus, by the second equation of \eqref{new:2*}
we have $\pd_Nw_{\pm j}(x_N) = 0$ on $\BR^\pm$.  But, 
$w_{\pm j}(x_N) \to 0$ as $\pm x_N\to \infty$, so that
$w_{\pm j}(x_N) =0$.  Thus, we have $\bw_\pm=0$ when
${\rm Im}\,\lambda=0$, which leads to a contradiciton. 
Therefore, we have $\det \tilde L\not=0$ for $(\tilde\lambda,
\tilde A) \in D_\epsilon(R_1, R_2)$.  Since 
$D_\epsilon(R_1, R_2)$ is compact, we have 
$$\inf_{(\tilde \lambda, \tilde A) \in D_{\epsilon}(R_1, R_2)}
|\det \tilde L| = c > 0,
$$
which, combined with \eqref{hom:1}, furnishes that
\begin{equation}\label{lop:1.3}
|\det L| \geq c(|\lambda|^{1/2} + A)^4
\end{equation}
provided that $R_2^{-1}|\lambda|^{1/2} \leq A \leq R_1|\lambda|^{1/2}$
and $\lambda \in \Sigma_\epsilon$. 
Setting $\omega = \min(c, \frac12\omega_1, \frac12\omega_2)$, by
\eqref{lop:1.1}, \eqref{lop:1.2} and \eqref{lop:1.3}, 
we have \eqref{lop:1}.

Next, we prove \eqref{lop:2}.  Recalling \eqref{sym-est:1} and 
the formula \eqref{det:1}, by Lemma \ref{lem:sym} we have
\begin{equation}\label{lop:1.4}
|\pd_{\xi'}^{\kappa'}\{(\tau\pd_\tau)^\ell \det L\}|
\leq C_{\kappa'}(|\lambda|^{1/2} + A)^4A^{-|\kappa'|}
\quad(\ell = 0, 1)
\end{equation}
for any multi-index $\kappa' \in \BN_0^{N-1}$ and 
$(\lambda, \xi') \in \tilde\Sigma_{\epsilon,0}$.
Thus, by the Bell formula \eqref{5.3}
with $f(t) =  1/t$ and $g(\xi') = \det L$, \eqref{lop:1.4} and 
\eqref{lop:1}, we have
$$|\pd_{\xi'}^{\kappa'}(\det L)^{-1}|
\leq C_{\kappa'}\sum_{\ell=1}^{|\kappa'|}
|\det L|^{-(\ell+1)}(|\lambda|^{1/2} + A)^{4\ell}A^{-|\kappa'|}
\leq C_{\kappa'}(|\lambda|^{1/2} + A)^{-4}A^{-|\kappa'|},
$$
which shows \eqref{lop:2} with $\ell=0$.  Analogously, we have
\eqref{lop:2} with $\ell=1$, which completes the proof of 
Lemma \ref{lem:lop}. 
\end{proof}
%%%%%%%%%%%%%%%%%%%%%%%%%%%%
\section{Problem with surface tension and height function}
In this section, we consider the problem:
\begin{alignat}2
&\rho_+\lambda\bu_+ - \DV\bS_+(\bu_+) = 0 &\quad&\text{in $\BR^N_+$},
\nonumber\\
&\rho_-\lambda\bu_- - \DV\bS_-(\bu_-) + \nabla\pi_-
= 0, \enskip \dv\bu_-=0  &\quad&\text{in $\BR^N_-$},
\nonumber\\
&\mu_-\bD_{mN}(\bu_-)|_- - \mu_+\bD_{mN}(\bu_+)|_+ = 0, \nonumber \\
&(\mu_-D_{NN}(\bu_-)-\pi_-)|_- = \sigma_-\Delta'H, \nonumber \\
&(\mu_+\bD_{NN}(\bu_+) + (\nu_+-\mu_+)\dv\bu_+)|_+
= \sigma_+\Delta'H, \nonumber\\
&u_{-m}|_- - u_{+m}|_+ =0, \nonumber\\
&\lambda H - \Bigl(\frac{\rho_-}{\rho_--\rho_+}u_{-N}|_-
- \frac{\rho_+}{\rho_--\rho_+}u_{+N}|_+\Bigr)
= d, \label{eq:1.2}
\end{alignat}
where, $\sigma_\pm = \frac{\rho_\pm\sigma}{\rho_--\rho_+}$,
and prove the following theorem.
\begin{thm}\label{thm:4.1} Let $1 < q < \infty$ and $0 < \epsilon < \pi/2$.
Then, there exist a $\lambda_0 > 0$ depending solely on 
$\mu_\pm$, $\nu_+$, $\rho_\pm$ and $\epsilon$ and operator families
$\CU_\pm(\lambda) \in \Hol(\Sigma_{\epsilon, \lambda_0}, 
\CL(W^2_q(\BR^N), W^2_q(\BR^N_\pm)^N))$,
$\CP_{-3}(\lambda) \in \Hol(\Sigma_{\epsilon, \lambda_0}, 
\CL(W^2_q(\BR^N), \hat W^1_q(\BR^N_-)))$ and $\CH(\lambda) 
\in \Hol(\Sigma_{\epsilon, \lambda_0}, \CL(W^2_q(\BR^N), W^3_q(\BR^N)))$
such that for any $\lambda \in \Sigma_{\epsilon,\lambda_0}$ and 
$d \in W^2_q(\BR^N)$, problem \eqref{eq:1.2} admits a unique solutions
$\bu_\pm = \CU_\pm(\lambda)d$, $\pi_-
= \CP_{-3}(\lambda)d$ and $H = \CH(\lambda)d$, and 
\begin{align*}
&\CR_{\CL(W^2_q(\BR^N), L_q(\BR^N_\pm)^{N+N^2+N^3})}
(\{(\tau\pd_\tau)^\ell (G_\lambda \CU_\pm(\lambda)) \mid
\lambda\in\Sigma_{\epsilon, \lambda_0}\}) \leq \gamma \quad(\ell=0, 1), \\
&\CR_{\CL(W^2_q(\BR^N), L_q(\BR^N_-)^N)}
(\{(\tau\pd_\tau)^\ell(\nabla\CP_{-3}(\lambda)) \mid 
\lambda\in\Sigma_{\epsilon, \lambda_0}\}) \leq \gamma \quad(\ell=0, 1), \\
&\CR_{\CL(W^2_q(\BR^N), L_q(\BR^N)^{N+1})}
(\{(\tau\pd_\tau)^\ell ((\lambda, \nabla) \CH(\lambda)) \mid
\lambda\in\Sigma_{\epsilon, \lambda_0}\}) \leq \gamma \quad(\ell=0, 1)
\end{align*}
with some constant $\gamma$ depending solely on 
$\mu_\pm$, $\nu_+$, $\rho_\pm$ and $\epsilon$.
\end{thm}
\begin{remark} Combining Theorem \ref{thm:4.1} with 
 Theorem \ref{thm:comp}, Theorem \ref{incomp} and 
Theorem \ref{thm:au}, we have Theorem \ref{main:r-bound} immediately.
\end{remark}
As was discussed in Sect.3, applying the partial Fourier transform to 
\eqref{eq:1.2}, we have the equations \eqref{eq:eq2} with interface condition:
\begin{equation}\label{eq:1.3}
\begin{cases}
\mu_{-}(D_N\hat u_{-m} + i\xi_m\hat u_{-N})|_-
- \mu_{+}(D_N\hat u_{-m} + i\xi_m\hat u_{-N})|_+
= 0,\\
(2\mu_{-}D_N\hat u_{-N} - \hat \pi_-)|_- = -\sigma_-A^2\hat H(0), \\
(2\mu_{+}D_N\hat u_{-N} + (\nu_{+}-\mu_{+})
(i\xi'\cdot\hat u'_+ + D_N\hat u_{+N})|_+ 
= -\sigma_+A^2\hat H(0), \\
\hat u_{-m}|_- - \hat u_{+m}|_+ = 0 
\end{cases}
\end{equation}
and the the resolvent equation for $H$:
\begin{equation}\label{eq:1.4}
\lambda \hat H(0) - \Bigl(\frac{\rho_-}{\rho_--\rho_+}\hat u_{-N}|_-
-\frac{\rho_+}{\rho_--\rho_+}\hat u_{+N}|_+\Bigr) = \hat d(0).
\end{equation}
We look for solutions $\hat u_{\pm J}$ and $\hat p_-$ of 
the form \eqref{eq:sol-form} with $h_j=0$, so that especially $\hat u^\pm_{N,3}
= 0$.  Our task is to represent $\hat H$ in terms of $\hat d(0)$. In view of
\eqref{eq:sol1}, we have $\hat u_{\pm J}(0) = \beta_\pm$, 
so that by \eqref{eq:eq7} with $\hat h_j(0) = 0$, we have 
\begin{equation}\label{eq:1.5} \begin{split}
\beta_{+N} &= -\frac{A}{\det L}
(\CL_{22}\sigma_-A^2 + \CL_{23}\sigma_+A)\hat H(0), \\
\beta_{-N} &= -\frac{A}{\det L}
(\CL_{32}\sigma_-A^2 + \CL_{33}\sigma_+A)\hat H(0).
\end{split}\end{equation}
Inserting these formulas into \eqref{eq:1.4}, we have
\begin{equation}\label{eq:1.6}
(\lambda + K)\hat H(0) = \hat d(0)
\end{equation}
with
\begin{equation}\label{eq:1.6*}
K = \frac{A^3}{\det L}\Bigl(\frac{\rho_-\sigma_-}{\rho_--\rho_+}\CL_{32}
-\frac{\rho_+\sigma_-}{\rho_--\rho_+}\CL_{22}\Bigr)
+ \frac{A^2}{\det L}\Bigl(\frac{\rho_-\sigma_+}{\rho_--\rho_+}\CL_{33}
-\frac{\rho_+\sigma_+}{\rho_--\rho_+}\CL_{23}\Bigr).
\end{equation}
We prove that 
\begin{lem}\label{lem:4.1} Let $0 < \epsilon < \pi/2$ and 
let $K$ be the function defined in \eqref{eq:1.6*}.  Then, there exists 
a positive constant $\lambda_0$ depending on $\epsilon$, $\mu_\pm$,
$\nu_+$ and $\rho_\pm$ such that 
\begin{equation}\label{eq:4.1}
|\pd_{\xi'}^{\kappa'}((\tau\pd_\tau)^\ell(\lambda+K)^{-1})|
\leq C_{\kappa'}(|\lambda|+A)^{-1}A^{-|\kappa'|}
\quad(\ell=0,1)
\end{equation}
for any multi-index $\kappa' \in \BN_0^{N-1}$ and $(\lambda, \xi')
\in \tilde\Sigma_{\epsilon, \lambda_0}$ with some constant $C_{\kappa'}$
depending on $\kappa'$, $\lambda_0$, $\epsilon$, $\mu_\pm$,
$\nu_+$ and $\rho_\pm$.
\end{lem}
\begin{proof}
To prove \eqref{eq:4.1} with $\kappa' = 0$ and $\ell=0$, 
first we consider the case 
where $R_1|\lambda|^{1/2} \leq A$ with large $R_1$. In the following,
$\delta_1$ is the same small number as in the proof of Lemma
\ref{lem:lop}. We know the 
asymptotic formula for $\det L$ given in \eqref{lop:1.1*}. On the other hand,
by \eqref{lop:0.1}, \eqref{lop:0.2} and 
\eqref{eq:eq6}, we have 
\begin{align*}
\CL_{32} &=(L_{11}^+ + L_{11}^-)L_{22}^+ - L_{12}^+L_{21}^+\\
&=\Bigl\{\Bigl(\frac{2\mu_+(\mu_++\nu_+)}{2\mu_++\nu_+} + 2\mu_-\bigr)
\frac{2\mu_+(\mu_++\nu_+)}{2\mu_++\nu_+}
-\Bigl(\frac{2(\mu_+)^2}{2\mu_++\nu_+}\Bigr)^2\Bigr\}A^2(1+O(\delta_1))
\\
& = \frac{4\mu_+((\mu_++\mu_-)\nu_++\mu_+\mu_-)}{2\mu_++\nu_+}
A^2(1 + O(\delta_1)), \\
\CL_{22} & = L_{21}^+L_{12}^- = A^2O(\delta_1), \quad
\CL_{33} = -L_{12}^+L_{21}^- = A^3O(\delta_1), \\
\CL_{23} & = L_{12}^-L_{21}^- - (L_{11}^++L_{11}^-)L_{22}^- \\
& = -4\Bigl(\frac{\mu_+(\mu_++\nu_+)}{2\mu_++\nu_+}
+\mu_-\Bigr)\mu_-A^3(1+O(\delta_1)).
\end{align*}
Since $\frac{\rho_\pm\sigma_\pm}{\rho_--\rho_+} 
= \Bigl(\frac{\rho_\pm\sigma}{\rho_--\rho_+}\Bigr)^2$, we have
$$
A^3\Bigl(\frac{\rho_-\sigma_-}{\rho_--\rho_+}\CL_{32} - 
\frac{\rho_+\sigma_-}{\rho_--\rho_+}\CL_{22}\Bigr)
+ A^2\Bigl(\frac{\rho_-\sigma_+}{\rho_--\rho_+}\CL_{33} - 
\frac{\rho_+\sigma_+}{\rho_--\rho_+}\CL_{23}\Bigr)\\
= \omega_3A^5(1 + O(\delta_1))
$$
with 
$$\omega_3
=4\frac{\mu_+((\mu_++\mu_-)\nu_++\mu_+\mu_-)}{2\mu_++\nu_+}
\Bigl(\frac{\rho_-}{\rho_--\rho_+}\Bigr)^2
+ 4\Bigl(\frac{\mu_+(\mu_++\nu_+)}{2\mu_++\nu_+}
+\mu_-\Bigr)\mu_-\Bigl(\frac{\rho_+}{\rho_--\rho_+}\Bigr)^2.
$$
Thus, by \eqref{lop:1.1*} we have
\begin{equation}\label{eq:1.8}
K = \frac{\omega_3}{\omega_1}A(1 + O(\delta_1)).
\end{equation}
Since $\lambda \in \Sigma_\epsilon$, by \eqref{eq:1.8} and \eqref{lem:3.1}
we have
$$|\lambda +K| \geq (\sin\frac{\epsilon}{2})(|\lambda| + 
\frac{\omega_3}{\omega_1}A) - \frac{\omega_3}{\omega_1}AO(\delta_1).$$
If we choose $\delta_1$ so small that 
$O(\delta_1) \leq \frac12\sin\frac{\epsilon}{2}$, we have 
\begin{equation}\label{eq:1.9}
|\lambda + K| \geq (\frac12\sin\frac{\epsilon}{2})(|\lambda| + 
\frac{\omega_3}{\omega_1}A)
\end{equation}
provided that $R_1|\lambda|^{1/2} \leq A$ with large $R_1 > 0$
and $\lambda \in \Sigma_\epsilon$.

Next, we consider the case where $A \leq R_1|\lambda|^{1/2}$. 
By \eqref{lop:1} we have  
\begin{equation}\label{eq:1.10}
|\det L|^{-1} \leq \omega^{-1}(1+R_1)^4|\lambda|^2
\end{equation}
for any $(\lambda, \xi') \in\tilde\Sigma_{\epsilon, 0}$
provided that $A \leq R_1|\lambda|^{1/2}$. On the other hand, 
by \eqref{lem:3.1} and \eqref{5.5}, we have
\begin{equation}\label{eq:1.11}
|P(\lambda, \xi')| \leq \frac{C(|\lambda| + A^2)}
{(\sin\frac{\epsilon}{2})(\rho_+(2\mu_++\nu_+)^{-1}|\lambda|
+A^2} \leq C.
\end{equation}
Here and in the sequel, $C$ denotes a generic constant depending on 
$R_1$, $\mu_\pm$, $\nu_+$, $\rho_\pm$ and $\epsilon$. 
By \eqref{eq:1.11} and \eqref{5.8}, we have
\begin{equation}\label{eq:1.12}
|L_{11}^+| \leq C|\lambda|^{1/2}, \enskip 
|L_{12}^+| \leq C|\lambda|, \enskip 
|L_{21}^+| \leq C, 
\enskip |L_{22}^+| \leq C|\lambda|^{1/2}.
\end{equation}
Moreover, by the definition of $L_{ij}^-$ given in \eqref{eq:eq5}
we have easily
$$ 
|L_{11}^-| \leq C|\lambda|^{1/2}, \enskip 
|L_{12}^-| \leq C|\lambda|, \enskip 
|L_{21}^-| \leq C|\lambda|^{1/2}, 
\enskip |L_{22}^-| \leq C|\lambda|,
$$
which, combined with \eqref{eq:1.10} and \eqref{eq:1.12}, 
furnishes that
\begin{equation}\label{eq:1.12*}
|K| \leq C|\lambda|^{1/2}
\end{equation}
for any $(\lambda, \xi') \in\tilde\Sigma_{\epsilon, 0}$
provided that $A \leq R_1|\lambda|^{1/2}$. 
Thus, we have $|\lambda + K| \geq |\lambda|^{1/2}(|\lambda|^{1/2}
 - C)$, so that choosing $\lambda_0 > 0$ so large that
$C\lambda_0^{1/2} \leq 1/2$, we have $|\lambda + K|
\geq \frac12|\lambda|$ for any 
$(\lambda, \xi') \in\tilde\Sigma_{\epsilon, \lambda_0}$
provided that $A \leq R_1|\lambda|^{1/2}$. Since $A \leq R_1|\lambda|^{1/2}$,
we observe that 
\begin{align*}
|\lambda + K| \geq \frac14|\lambda| + \frac14|\lambda|
\geq \frac14|\lambda| + \frac{\lambda^{1/2}_0}{4}|\lambda|^{1/2}
\geq \frac14(|\lambda| + \lambda^{1/2}_0R_1^{-1}A).
\end{align*}
Choosing $R_1 > 0$ so large that $\lambda^{1/2}_0R_1^{-1} \leq 
\frac{\omega_3}{\omega_1}$, we have 
\begin{equation}\label{eq:1.13}
|\lambda + K| \geq \frac14(|\lambda| + \frac{\omega_3}{\omega_1}A)
\end{equation}
for any $(\lambda, \xi') \in\tilde\Sigma_{\epsilon, 0}$
provided that $A \leq R_1|\lambda|^{1/2}$. 
Since $\Sigma_{\epsilon, \lambda_0} \subset \Sigma_{\epsilon}$, 
combining \eqref{eq:1.9} and \eqref{eq:1.13}, we have
\begin{equation}\label{eq:1.14}
|\lambda + K| \geq \omega_4(|\lambda| + A)
\end{equation}
for any $(\lambda, \xi') \in \tilde\Sigma_{\epsilon, \lambda_0}$
with $\omega_4 = \min(\frac14, \frac14\frac{\omega_3}{\omega_1}, 
\frac12\sin\frac{\epsilon}{2},  \frac12\sin\frac{\epsilon}{2}
\frac{\omega_3}{\omega_1})$. 

Next, we prove \eqref{eq:4.1} for any multi-index $\kappa' \in \BN_0^{N-1}$.
By \eqref{sym-est:1}, \eqref{sym-est:3} and Lemma \ref{lem:sym},
we have $K \in \bM_{1,2}(0)$, so that by Bell's formula
\eqref{5.3} with $f(t) = (\lambda + t)^{-1}$ and $g = K$, we have
$$|\pd_{\xi'}^{\kappa'}(\lambda + K)^{-1}| 
\leq C_{\kappa'}\sum_{\ell=1}^{|\kappa'|}
|\lambda+K|^{-(\ell+1)}(|\lambda|^{1/2} + A)^{\ell}A^{-|\kappa'|}
\leq C_{\kappa'}(|\lambda|+A)^{-1}A^{-|\kappa'|}.
$$
Analogously, we have 
$$|\pd_{\xi'}^{\kappa'}(\tau\pd_\tau(\lambda+K)^{-1})|
\leq C_{\kappa'}(|\lambda|+A)^{-1}A^{-|\kappa'|}.
$$
Summing up, we have obtained \eqref{eq:4.1}.
This completes the proof of Lemma \ref{lem:4.1}.
\end{proof}
By \eqref{eq:1.6} and Lemma \ref{lem:4.1}, we have
\begin{equation}\label{eq:h1}
\hat H(\xi', 0) = (\lambda + K)^{-1}\hat d(\xi', 0),
\end{equation}
so that we define $\hat H(\xi', x_N)$ by $\hat H(\xi', x_N) = 
e^{-(1+A^2)^{1/2}x_N}(\lambda + K)^{-1}\hat d(\xi', 0)$.
We have the following lemma.
\begin{lem}\label{lem:4.2}.  Let $1 < q < \infty$, $0 < \epsilon
< \pi/2$ and let $\lambda_0$ be the same constant as in 
Lemma \ref{lem:4.1}.
Given that the operator $\tilde\CH(\lambda)$ is defined by 
\begin{equation}\label{eq:h2}
[\tilde\CH(\lambda)d](x) = \CF^{-1}_{\xi'}
[e^{-(1+A^2)^{1/2}x_N}(\lambda + K)^{-1}\hat d(\xi', 0)](x')
\quad\text{for $x \in \BR^N_+$} 
\end{equation}
for any $d \in W^2_q(\BR^N)$,  $\tilde\CH(\lambda)
\in \Hol(\Sigma_{\epsilon, \lambda_0}, 
\CL(W^2_q(\BR^N), W^3_q(\BR^N_+)))$,
and 
\begin{equation}\label{eq:h3}
\CR_{\CL(W^2_q(\BR^N), W^2_q(\BR^N_+)^{1+N})}
(\{(\tau\pd_\tau)^\ell(\lambda, \nabla)\tilde\CH(\lambda) \mid 
\lambda \in \Sigma_{\epsilon, \lambda_0}\}) \leq \gamma
\quad(\ell=0,1)
\end{equation}
with some constant $\gamma$ depending on $\mu_\pm$, 
$\nu_+$, $\rho_\pm$,
$\epsilon$ and $\lambda_0$. 
\end{lem}
\begin{proof}
First, using the same idea as in \eqref{cont:1} and \eqref{eq:inv1}
and the identity: $1 = \frac{1+A^2}{1+ A^2} =
\frac{1}{1 + A^2} -\sum_{k=1}^{N-1} \frac{(i\xi_k)(i\xi_k)}{1 + A^2}$,
we rewrite $\tilde\CH(\lambda)d$ as follows:
\begin{align}
\tilde\CH(\lambda)d &= \CF^{-1}_{\xi'}\Bigl[\frac{e^{-(1+A^2)^{1/2}x_N}}
{\lambda+K}\hat d(\xi', 0)\Bigr] 
%\nonumber\\
%&
= -\int^\infty_0\CF^{-1}_{\xi'}\Bigl[
\frac{\pd}{\pd y_N}
\frac{e^{-(1+A^2)^{1/2}(x_N+y_N)}}
{\lambda+K}\hat d(\xi', y_N)\Bigr]\,dy_N \nonumber \\
%& = -\int^\infty_0\CF^{-1}_{\xi'}
%\Bigl[\frac{e^{-(1+A^2)^{1/2}(x_N+y_N)}}{\lambda+K}
%\widehat{\pd_Nd}(\xi', y_N)\Bigr](x')\,dy_N \nonumber \\
%&\quad
%+ \int^\infty_0\CF^{-1}_{\xi'}
%\Bigl[\frac{e^{-(1+A^2)^{1/2}(x_N+y_N)}(1 + A^2)^{1/2}}{\lambda+K}
%\hat d(\xi', y_N)\Bigr](x')\,dy_N \nonumber\\
& = -\int^\infty_0\CF^{-1}_{\xi'}
\Bigl[\frac{e^{-(1+A^2)^{1/2}(x_N+y_N)}}{(\lambda+K)(1+A^2)}
\widehat{\pd_Nd}(\xi', y_N)\Bigr](x')\,dy_N \nonumber\\
&\quad + \sum_{k=1}^{N-1}\int^\infty_0\CF^{-1}_{\xi'}
\Bigl[\frac{e^{-(1+A^2)^{1/2}(x_N+y_N)}(i\xi_k)}{(\lambda+K)(1+A^2)}
\widehat{\pd_k\pd_Nd}(\xi', y_N)\Bigr](x')\,dy_N \nonumber\\
&\quad + \int^\infty_0\CF^{-1}_{\xi'}
\Bigl[\frac{e^{-(1+A^2)^{1/2}(x_N+y_N)}}{(\lambda+K)(1 + A^2)^{1/2}}
\widehat{(1-\Delta') d}(\xi', y_N))\Bigr](x')\,dy_N \label{eq:h+}
\end{align}
for $x_N > 0$. 

To prove the $\CR$ boundedness of the operator family
$\{\tilde\CH(\lambda) \mid \lambda \in \Sigma_{\epsilon, \lambda_0}\}$, 
we prepare the following three lemmas.
\begin{lem}\label{lem:4.2*} Let $1 < q < \infty$ and
let $\ell(\xi')$ be a $C^\infty$
function defined on $\BR^{N-1}$ satisfying the estimates:
$$|\pd_{\xi'}^{\kappa'}\ell(\xi')| \leq C_{\kappa'}
(1 + A)^{1-|\kappa'|}$$
for any multi-index $\kappa' \in \BN_0^{N-1}$ and 
$\xi' \in \BR^{N-1}\setminus\{0\}$. Let $T$ be the operator defined
by 
$$Tf = \int^\infty_0\CF_{\xi'}^{-1}[e^{-(1+A^2)^{1/2}
(x_N+y_N)}\ell(\xi')\hat f(\xi', y_N)\Bigr](x')\,dy_N$$
for $f \in L_q(\BR^N_+)$. Then, 
 $T \in \CL(L_q(\BR^N_+))$.
\end{lem}
\begin{proof}
If we define a function $k(x)$ by
$k(x) = \CF^{-1}_{\xi'}[e^{-(1+A^2)^{1/2}x_N}\ell(\xi')](x')$, 
we have
$$[Tf](x) = \int^\infty_0\int_{\BR^{N-1}}k(x'-y', x_N+y_N)f(y', y_N)\,dy'dy_N.
$$
Our task is to prove that
\begin{equation}\label{eq:1.16}
|k(x)|\leq C|x|^{-N}.
\end{equation}
In fact, if we have \eqref{eq:1.16}, then by Young's inequality we have
$$\|[Tf](\cdot, x_N)\|_{L_q(\BR^{N-1})}
\leq \int^\infty_0\|k(\cdot, x_N+y_N)\|_{L_1(\BR^{N-1})}
\|f(\cdot, y_N)\|_{L_q(\BR^{N-1})}\,dy_N.
$$
By \eqref{eq:1.16}, we have $\|k(\cdot, x_N+y_N)\|_{L_1(\BR^{N-1})}
\leq C(x_N+y_N)^{-1}$, so that 
$$\|[Tf](\cdot, x_N)\|_{L_q(\BR^{N-1})}
\leq C\int^\infty_0\frac{\|f(\cdot, y_N)\|_{L_q(\BR^{N-1})}}{x_N+y_N}\,
dy_N.$$
By the Hardy inequality (cf. Stein \cite[p.271 A.3]{Stein}), we have
$\|Tf\|_{L_q(\BR^N_+)} \leq C\|f\|_{L_q(\BR^N_+)}$. Thus, we have proved 
$T \in \CL(L_q(\BR^N_+))$.

To prove \eqref{eq:1.16}, first we observe that 
\begin{equation}\label{eq:1.17}
k(x) = \sum_{|\alpha'| = N-1}\Bigl(\frac{ix'}{|x'|^2}\Bigr)^{\alpha'}
\int_{\BR^{N-1}}e^{ix'\cdot\xi'}\pd_{\xi'}^{\alpha'}(
e^{-(1+A^2)^{1/2}x_N}\ell(\xi'))\,d\xi'.
\end{equation}
For any multi-index $\beta' \in \BN_0^{N-1}$, we have
\begin{align*}
|\pd_{\xi'}^{\beta'}\{(\pd_{\xi'}^{\alpha'}
e^{-(1 + A^2)^{1/2}}\ell(\xi'))\}|
&\leq C_{\alpha', \beta'}e^{-c_{\alpha', \beta'}(1 + A^2)^{1/2}x_N}
(1 + A)^{1-|\alpha'| - |\beta'|} \\
&= C_{\alpha', \beta'}e^{-c_{\alpha', \beta'}(1 + A^2)^{1/2}x_N}
(1 + A)^{2-N-|\beta'|}
\end{align*}
with some positive constants $C_{\alpha', \beta'}$ and
$c_{\alpha', \beta'}$.  Since $2-N \leq 0$, we have 
$$|\pd_{\xi'}^{\beta'}\{(\pd_{\xi'}^{\alpha'}
e^{-(1 + A^2)^{1/2}}\ell(\xi'))\}|
\leq 
C_{\alpha', \beta'}e^{-c_{\alpha', \beta'}(1 + A^2)^{1/2}x_N}
A^{2-N-|\beta'|}
$$
for any multi-index $\beta' \in \BN_0^{N-1}$.  Thus, by
the result due to Shibata and Shimizu \cite{SS1}, we have 
$$
\Bigl|\int_{\BR^{N-1}}e^{ix'\cdot\xi'}\pd_{\xi'}^{\alpha'}(
e^{-(1+A^2)^{1/2}x_N}\ell(\xi'))\,d\xi'
\Bigr| \leq C|x'|^{-(N-1+(2-N)} = C|x'|^{-1},$$
which, combined with \eqref{eq:1.17},  furnishes that
\begin{equation}\label{eq:1.18}
|k(x)| \leq C|x'|^{-N}.
\end{equation}
On the other hand, 
\begin{align*}
|k(x)| &\leq \int_{\BR^{N-1}}e^{-(1+A^2)^{1/2}x_N}|\ell(\xi')|\,d\xi'
\leq C\int_{\BR^{N-1}}e^{-(1+A^2)^{1/2}x_N}(1+A)\,d\xi' \\
& \leq C\int_{\BR^{N-1}}\frac{d\xi'}{((1+A^2)^{1/2}x_N)^N}
+ C\int_{\BR^{N-1}}e^{-Ax_N}A\,d\xi' \\
& \leq \frac{C}{(x_N)^N}\int^\infty_0\frac{r^{N-2}}{(1+r^2)^{\frac{N}{2}}}
\,dr + \frac{C}{(x_N)^N}\int_{\BR^{N-1}}e^{-|\eta'|}|\eta'|\,d\eta'
\leq \frac{C}{(x_N)^N},
\end{align*}
which, combined with \eqref{eq:1.18}, furnishes \eqref{eq:1.16}.
This completes the proof of Lemma \ref{lem:4.2*}.
\end{proof}
The following lemma was proved in Enomoto-Shibata \cite{ES}
\begin{lem}\label{lem:4.3} Let $1 < q < \infty$ and let 
$\Lambda$ be a subset of $\BC$.  Let $m(\lambda, \xi')$ be a 
function defined on $\Lambda \times (\BR^{N-1}\setminus\{0\})$ such 
that for any multi-index $\alpha \in \BN^N_0$ 
 there exists a constant $C_{\alpha'}$ depending on
$\alpha'$ and $\Lambda$ such that 
\begin{equation}\label{3.6} |\pd_\xi^{\alpha'} m(\lambda, \xi')|
\leq C_\alpha A^{-|\alpha'|}
\end{equation}
for any $(\lambda, \xi') \in \Lambda\times(\BR^{N-1}\setminus\{0\})$.
Let $K_\lambda$ be an operator defined by $K_\lambda f
= \CF^{-1}_{\xi'}[m(\lambda, \xi')\hat f(\xi')]$.  Then, the set 
$\{K_\lambda \mid \lambda \in \Lambda\}$ is $\CR$-bounded on 
$\CL(L_q(\BR^{N-1}))$ and 
\begin{equation}\label{3.6*}
\CR_{\CL(L_q(\BR^{N-1}))}(\{K_\lambda \mid \lambda \in \Lambda\})
\leq C_{q, N}\max_{|\alpha'| \leq N} C_\alpha
\end{equation}
with some constant $C_{q, N}$ that depends solely on $q$ and $N$. 
\end{lem}
%%%%%%%%%%%%%%%%%%%
%%%%%%%%%%%%%%%%%%%
From the definition of $\CR$ boundedness 
we have the following lemma immediately.
\begin{lem}\label{lem:4.4} Let $1 < q < \infty$ and let 
$\Lambda$ be a subset of $\BC$.  Let $\{\CS_\lambda \mid \lambda \in \Lambda\}$
be an $\CR$ bounded operator family on $\CL(L_q(\BR^{N-1})$ and 
$T$ a bounded linear operator in $\CL(L_q(\BR^N_+))$.  Then, 
$\{\CS_\lambda T \mid \lambda \in \Lambda\}$
is an $\CR$-bounded operator family on $\CL(L_q(\BR^N_+))$ and 
\begin{equation}\label{eq:1.19}
\CR_{\CL(L_q(\BR^N_+))}(\{\CS_\lambda T 
\mid \lambda \in \Lambda\}) \leq \CR(\CS_\lambda)\|T\|_{\CL(L_q(\BR^N_+))}
\end{equation}
where we have set $\CR(\CS_\lambda) 
= \CR_{L_q(\BR^{N-1})}(\{\CS_\lambda \mid
\lambda \in \Lambda\}$ for short. 
\end{lem}
%\begin{proof}
%For any $n \in \BN$, let $\{f_j\}_{j=1}^n$ 
%be any $n$ elements of $L_q(\BR^N_+)$,
%$\{\lambda_j\}_{j=1}^n$ any $n$ elements of $\Lambda$ and 
%$\{\epsilon_j\}_{j=1}^n$ any $n$ independent ,symmetric, 
%$\{-1,1\}$ valued random variables on $[0, 1]$. We observe that 
%\begin{align*}
%&\int^1_0\|\sum_{j=1}^n\epsilon_j(\sigma)\CS_{\lambda_j}Tf_j\|_{L_q(\BR^N_+)}^q
%\,d\sigma
%= \int^\infty_0\int^1_0
%\|\sum_{j=1}^n\epsilon_j(\sigma)\CS_{\lambda_j}Tf_j\|_{L_q(\BR^N_+)}^q
%\,d\sigma\,dx_N\\
%&\quad \leq \CR(\CS_\lambda)^q
%\int^\infty_0\int^1_0
%\|\sum_{j=1}^n\epsilon_j(\sigma)Tf_j\|_{L_q(\BR^N_+)}^q
%\,d\sigma\,dx_N
%= \CR(\CS_\lambda)^q
%\int^1_0\|\sum_{j=1}^n\epsilon_j(\sigma)Tf_j\|_{L_q(\BR^N_+)}^q\,d\sigma \\
%&\quad
%= \CR(\CS_\lambda)^q
%\int^1_0\|T(\sum_{j=1}^n\epsilon_j(\sigma)f_j)\|_{L_q(\BR^N_+)}^q\,d\sigma 
%\leq \CR(\CS_\lambda)^q\|T\|_{\CL(\BR^N_+)}^q
%\int^1_0\|\sum_{j=1}^n\epsilon_j(\sigma)f_j\|_{L_q(\BR^N_+)}^q\,d\sigma.
%\end{align*}
%Thus, we have the lemma.
%\end{proof}
Under these preparations, we finish proving Lemma \ref{lem:4.2}. For any 
multi-index $\alpha \in \BN_0^N$ with $|\alpha| \leq 2$, 
using \eqref{eq:h+}, we write 
\begin{align*}
&\pd_x^\alpha(\lambda, \nabla)\tilde\CH(\lambda)d  
= -\int^\infty_0\CF^{-1}_{\xi'}
\Bigl[\frac{e^{-(1+A^2)^{1/2}(x_N+y_N)}\ell_\alpha(\xi')}{(1+A^2)}
\frac{(\lambda, i\xi', -(1+A^2)^{1/2})}{(\lambda+K)}
\widehat{\pd_Nd}(\xi', y_N)\Bigr](x')\,dy_N \nonumber\\
&\quad + \sum_{k=1}^{N-1}\int^\infty_0\CF^{-1}_{\xi'}
\Bigl[\frac{e^{-(1+A^2)^{1/2}(x_N+y_N)}(i\xi_k)\ell_\alpha(\xi')}
{(1+A^2)}\frac{(\lambda, i\xi', -(1+A^2)^{1/2})}{(\lambda+K)}
\widehat{\pd_k\pd_Nd}(\xi', y_N)\Bigr](x')\,dy_N \nonumber\\
&\quad + \int^\infty_0\CF^{-1}_{\xi'}
\Bigl[\frac{e^{-(1+A^2)^{1/2}(x_N+y_N)}\ell_\alpha(\xi')}
{(1 + A^2)^{1/2}}\frac{(\lambda, i\xi', -(1+A^2)^{1/2})}{(\lambda+K)}
\widehat{(1-\Delta') d}(\xi', y_N)\Bigr](x')\,dy_N
\end{align*}
where $\ell_\alpha(\xi')$ is an symbol satisfying the estimate:
$$|\pd_{\xi'}^{\kappa'}\ell_\alpha(\xi')| \leq C_{\kappa'}
(1 + A)^{|\alpha|-|\kappa'|}$$
for any multi-index $\kappa' \in \BN_0^{N-1}$. If we define the operators
$\CS_\lambda$, $T_1^\alpha$, $T_{2k}^\alpha$, and $T_3^\alpha$ by 
\begin{align*}
[\CS_\lambda g](x') &
= \CF^{-1}_{\xi'}\Bigl[\frac{(\lambda, i\xi', -(1+A^2)^{1/2})}{(\lambda+K)}\hat g(\xi')\Bigr](x'), \\
[T_1^\alpha f](x) & = \int^\infty_0\CF^{-1}_{\xi'}\Bigl[\frac{e^{-(1+A^2)^{1/2}(x_N+y_N)}\ell_\alpha(\xi')}{(1+A^2)}\hat f(\xi', y_N)\Bigr](x')\,dy_N,\\
[T_{2k}^\alpha f](x) & = \int^\infty_0\CF^{-1}_{\xi'}\Bigl[\frac{e^{-(1+A^2)^{1/2}(x_N+y_N)}(i\xi_k)\ell_\alpha(\xi')}
{(1+A^2)}\hat f(\xi', y_N)\Bigr](x')\,dy_N, \\
[T_3^\alpha f](x) & = \int^\infty_0\CF^{-1}_{\xi'}\Bigl[\frac{e^{-(1+A^2)^{1/2}(x_N+y_N)}\ell_\alpha(\xi')}
{(1 + A^2)^{1/2}}\hat f(\xi', y_N)\Bigr](x')\,dy_N
\end{align*}
for $g \in L_q(\BR^{N-1})$ and $f \in L_q(\BR^N_+)$, 
then, we have 
$$\pd_x^\alpha(\lambda, \nabla)\tilde\CH(\lambda)d 
= -\CS_\lambda T_1^\alpha(\pd_Nd) 
+ \sum_{k=1}^{N-1} \CS_\lambda T_{2k}^\alpha(\pd_k\pd_Nd)
+\CS_\lambda T_3^\alpha((1-\Delta')d).
$$
Since 
\begin{equation}\label{eq:1.20} \begin{split}
&|\pd_{\xi'}^{\kappa'}\Bigl(\frac{\ell_\alpha(\xi')}{1+A^2}\Bigr)| 
\leq C_{\kappa'}(1 +A)^{|\alpha|-2-|\kappa'|},  \quad 
|\pd_{\xi'}^{\kappa'}\Bigl(\frac{\ell_\alpha(\xi')(i\xi_k)}{1+A^2}\Bigr)| 
\leq C_{\kappa'}(1 +A)^{|\alpha|-1-|\kappa'|}, \\
&|\pd_{\xi'}^{\kappa'}\Bigl(\frac{\ell_\alpha(\xi')}{(1+A^2)^{1/2}}\Bigr)| 
\leq C_{\kappa'}(1 +A)^{|\alpha|-1-|\kappa'|},
\end{split}\end{equation}
for any multi-index $\kappa' \in \BN_0^{N-1}$ with some constant $C_{\kappa'}$,
by Lemma \ref{lem:4.2*} 
 $T_1^\alpha$, $T_{2k}^\alpha$ and $T_3^\alpha$ are bounded linear 
operators from $L_q(\BR^N_+)$ into itself.  By Lemma \ref{lem:4.1}, we have
$$|\pd_{\xi'}^{\kappa'}((\tau\pd_\tau)^\ell
\frac{(\lambda, i\xi', -(1+A^2)^{1/2})}{(\lambda+K)})|
\leq C_{\kappa'}A^{-|\kappa'|}\quad(\ell=0,1), $$
so that by Lemma \ref{lem:4.3} $(\tau\pd_\tau)^\ell\CS_\lambda$ ($\ell=0,1$)
are $\CR$ bounded operator families on $L_q(\BR^{N-1})$ and 
$$\CR_{\CL(L_q(\BR^{N-1}))}(\{(\tau\pd_\tau)^\ell\CS_\lambda \mid
\lambda \in \Sigma_{\epsilon, \lambda_0}\}) \leq \gamma \quad(\ell=0, 1)$$
with some constant $\gamma$ depending solely on $\mu_\pm$, $\nu_+$,
$\rho_\pm$, $\epsilon$ and $\lambda_0$.　Thus, Lemma \ref{lem:4.2} follows
from Lemma \ref{lem:4.4} immediately.  This completes the proof of
Lemma \ref{lem:4.2}. 
\end{proof}
By using the Lions method, we extend $\tilde\CH(\lambda)d$ to $x_N < 0$.
Namely, we define $\CH(\lambda)$ by 
$$[\CH(\lambda)d](x) = \begin{cases}
[\tilde\CH(\lambda)d](x) &\quad (x_N > 0) \\
\sum_{j=1}^4 a_j[\tilde\CH(\lambda)d](x', -jx_N) \quad(x_N < 0)
\end{cases}$$
where $a_j$ are constants satisfying the equations: 
$\sum_{j=1}^4 a_j(-j)^k = 1$ for $k=0,1,2,3$. By Lemma \ref{lem:4.2},
we have the following corollary of Lemma \ref{lem:4.2} immediately. 
\begin{cor}\label{cor:4.1} Let $1 < q < \infty$, $0 < \epsilon
< \pi/2$ and let $\lambda_0$ be the same constant as in Lemma \ref{lem:4.1}.
Then, there exists an operator family $\CH(\lambda) \in 
\Hol(\Sigma_{\epsilon, \lambda_0}, \CL(W^2_q(\BR^N), W^3_q(\BR^N)))$ such that
for any $\lambda \in \Sigma_{\epsilon, \lambda_0}$ and $d \in W^2_q(\BR^N)$, 
$\CF_{\xi'}^{-1}\Bigl[\frac{\hat d(\xi',0)}{\lambda+K}\Bigr](x')
= \CH(\lambda)d|_{x_N=0}$ and 
$$
\CR_{\CL(W^2_q(\BR^N), W^2_q(\BR^N)^{1+N})}
(\{(\tau\pd_\tau)^\ell(\lambda, \nabla)\CH(\lambda) \mid 
\lambda \in \Sigma_{\epsilon, \lambda_0}\}) \leq \gamma
\quad(\ell=0,1)
$$
with some constant $\gamma$ depending on $\mu_\pm$, $\nu_+$, $\rho_\pm$,
$\epsilon$ and $\lambda_0$. 
\end{cor}
Next, we give a solution operator of \eqref{eq:1.2} for velocity field. 
By \eqref{eq:sol-form} and \eqref{eq:1.6}, we have
$$\hat u_{\pm J} = AM_\pm(x_N)\frac{AR^\pm_{JN,0}}{\lambda+K}\hat d(0) 
+ Ae^{-\mp B_\pm x_N}\frac{AS^\pm_{JN,-1}}{\lambda +K}\hat d(0),
\quad p_- = e^{Ax_N}\frac{Ap^-_{N,1}}{\lambda+K}\hat d(0).
$$
Employing the same argument as in \eqref{cont:1} and \eqref{eq:inv1} 
and using \eqref{cont:3} and the identity: 
$1 = \frac{1 + A^2}{1 + A^2} = \frac{1}{1 + A^2} 
-\sum_{k=1}^{N-1}\frac{(i\xi_k)(i\xi_k)}{1 + A^2}$, we have 
\begin{align*}
& u_{\pm J}(x) = \CF^{-1}_{\xi'}[\hat u_{\pm J}(\xi', x_N)](x')  \\
&= -\int^{\pm\infty}_0\CF^{-1}_{\xi'}\Bigl[
\pd_N(AM_\pm(x_N+y_N)\frac{AR^\pm_{JN,0}}{\lambda + K}
\hat d(\xi', y_N))\Bigr](x')\,dy_N \\
&\quad- \int^{\pm\infty}_0\CF^{-1}_{\xi'}\Bigl[
\pd_N(Ae^{\mp B_\pm(x_N+y_N)}\frac{AS^\pm_{JN,-1}}{\lambda + K}
\hat d(\xi', y_N))\Bigr](x')\,dy_N\\
%& = -\int^{\pm\infty}_0\CF^{-1}_{\xi'}\Bigl[
%AM_\pm(x_N+y_N)\frac{AR^\pm_{JN,0}}{\lambda + K}
%\widehat{\pd_Nd}(\xi', y_N)\Bigr](x')\,dy_N \\
%&\quad\pm\int^{\pm\infty}_0\CF^{-1}_{\xi'}\Bigl[
%(Ae^{\mp B_\pm(x_N+y_N)} +AA_\pm M_\pm(x_N+y_N))
%\frac{AR^\pm_{JN,0}}{\lambda + K}\hat d(\xi', y_N)\Bigr](x')\,dy_N \\
%&\quad -\int^{\pm\infty}_0\CF^{-1}_{\xi'}\Bigl[
%Ae^{\mp B_\pm(x_N+y_N)}\frac{AS^\pm_{JN,-1}}
%{\lambda + K}\widehat{\pd_Nd}(\xi', y_N)\Bigr](x')\,dy_N\\
%&\quad \pm \int^{\pm\infty}_0\CF^{-1}_{\xi'}\Bigl[
%AB_\pm e^{\mp B_\pm(x_N+y_N)}\frac{AS^\pm_{JN,-1}}
%{\lambda + K}\hat d(\xi', y_N)\Bigr](x')\,dy_N\\
&=-\int^{\pm\infty}_0\CF^{-1}_{\xi'}\Bigl[
AM_\pm(x_N+y_N)\Bigl\{\frac{AR^\pm_{JN,0}}
{(\lambda+K)(1 + A^2)}\widehat{\pd_Nd}(\xi', y_N)\\
&\quad\quad -\sum_{k=1}^{N-1}\frac{AR^\pm_{JN,0}(i\xi_k)}{(\lambda+K)(1+A^2)}
\widehat{\pd_k\pd_Nd}(\xi', y_N) \mp
\frac{A_\pm AR^\pm_{JN,0}}{(\lambda + K)(1 + A^2)}
\widehat{(1-\Delta')d}(\xi', y_N)\Bigr\}\Bigr](x')\,dy_N\\
&\quad -\int^{\pm\infty}_0Ae^{\mp B_\pm(x_N+y_N)}
\Bigl\{\mp \frac{AR^\pm_{JN,0}}
{(\lambda + K)(1 + A^2)}\widehat{(1-\Delta')d}(\xi', y_N)
+ \frac{AS^\pm_{JN,-1}}{(\lambda + K)(1 + A^2)}\widehat{\pd_Nd}(\xi', y_N)\\
&\quad\quad -\sum_{k=1}^{N-1}
\frac{AS^\pm_{JN,-1}(i\xi_k)}{(\lambda + K)(1 + A^2)}\widehat{\pd_k\pd_Nd}
(\xi', y_N) \mp
\frac{B_\pm AS^\pm_{JN,-1}}{(\lambda + K)(1 + A^2)}\widehat{(1-\Delta')d}
(\xi', y_N)\Bigr\}\Bigr](x')\,dy_N
\end{align*}
Analogously, we have
\begin{align*}
p_-(x) & = \int^0_{-\infty}\CF^{-1}_{\xi'}\Bigl[e^{A(x_N+y_N)}
\frac{p^-_{N,1}}{\lambda + K}\Bigl\{\widehat{(1-\Delta')d}(\xi', y_N)
-\sum_{k=1}^{N-1}\widehat{\pd_k\pd_Nd}(\xi', y_N)\Bigr\}\Bigr](x')\,dy_N.
\end{align*}
By \eqref{sym-est:0}, \eqref{5.4}, Lemma \ref{lem:4.1} and 
\eqref{eq:1.20} with $\alpha=0$, we have 
\begin{equation}\label{eq:1.21}\begin{split}
&\frac{p^-_{N,1}}{\lambda + K} \in \bM_{0,2}(\lambda_0), \enskip 
\frac{AR^\pm_{JN,0}(i\xi_k)}{(\lambda + K)(1 + A^2)}, \enskip 
\frac{A_\pm AR^\pm_{JN,0}}{(\lambda + K)(1 + A^2)} 
\in \bM_{-1,2}(\lambda_0), \\
&\frac{AR^\pm_{JN,0}}{(\lambda + K)(1 + A^2)}, 
\enskip \frac{AS^\pm_{JN,-1}}{(\lambda + K)(1 + A^2)}, \enskip
\frac{AS^\pm_{JN,-1}(i\xi_k)}{(\lambda + K)(1 + A^2)}, \enskip 
\frac{B_\pm AS^\pm_{JN,-1}}{(\lambda + K)(1 + A^2)} \in \bM_{-2,2}(\lambda_0).
\end{split}\end{equation}
In fact,  by Leibniz's rule we have
$$|\pd_{\xi'}^{\alpha'}
\Bigl(\frac{A_+AR^+_{JN,0}}{(\lambda+K)(1 + A^2)}\Bigr)| 
\leq C_{\alpha'}\frac{|\lambda|^{1/2}+A}{(|\lambda| + A)(1 + A)}
A^{-|\alpha'|}.
$$
Since $(|\lambda| + A)(1 + A) = |\lambda| + |\lambda|A + A + A^2 \geq 
\frac12(|\lambda|^{1/2} + A)^2$, we have 
$\frac{A_\pm AR^\pm_{JN,0}}{(\lambda + K)(1 + A^2)} 
\in \bM_{-1,2}(\lambda_0)$. Analogously, we have other assertions 
in \eqref{eq:1.21}. Therefore, by Lemma \ref{lem:tech:comp}, 
Lemma \ref{lem:tech:incomp} and Corollary \ref{cor:4.1} we have 
Theorem \ref{thm:4.1}, which completes the proof of Theorem \ref{main:r-bound}.

%%%%%%%%%%%%%%%%
\section{A proof of theorem \ref{thm:stokes}}

First, we transfer problem \eqref{stokes:1}, \eqref{stokes:2*} and
\eqref{stokes:3} to the zero intial data case.  Let $e^{At}$ be the operator 
defined by
$e^{At}f = \CF^{-1}_\xi[e^{-(1+|\xi|^2)^{1/2}t}\,CF[f](\xi)]$. 
Here and in the sequel, $\CF[f](\xi)$ and $\CF^{-1}$
denote the Fourier transform $f$ on
$\BR^N$ and  the inverse  transform of $g(\xi)$ defind by 
$$\CF[f](\xi) = \int_{\BR^N} e^{-ix\cdot\xi}f(x)\,dx,\quad
\CF_\xi^{-1}[g(\xi)] = (2\pi)^{-N}\int_{\BR^N} e^{ix\cdot\xi}g(\xi)\,d\xi,$$
respectively. Since 
\begin{align*}
\|[\pd_t e^{At}f](\cdot, t)\|_{L_q(\BR^N)} & 
\leq Ct^{-1}e^{-t/2}\|f\|_{L_q(\BR^N)}, \\
\|[\pd_t e^{At}f](\cdot, t)\|_{L_q(\BR^N)} & 
\leq Ce^{-t/2}\|f\|_{W^1_q(\BR^N)}, \\
\|[e^{At}f](\cdot, t)\|_{L_q(\BR^N)} & \leq Ce^{-t/2}\|f\|_{L_q(\BR^N)}
\end{align*}
for $1 < q < \infty$ and $t>0$, employing the same real interpolation argument
as in the proof of Theorem 3.9 in \cite{SS2}, we have
$$\|\pd_te^{At}f\|_{L_p((0, \infty), L_q(\BR^N))}
+ \|e^{At}f\|_{L_p((0, \infty), W^1_q(\BR^N))} 
\leq C\|f\|_{B^{1-1/p}_{q,p}(\BR^N)}.$$
Thus, setting $I = e^{At}H_0$, we have $I|_{t=0} = H_0$ in $\BR^N$, and 
\begin{equation}\label{7.1}
\|\pd_tI\|_{L_p((0, \infty), W^2_q(\BR^N))}
+ \|I\|_{L_p((0, \infty), W^3_q(\BR^N))} 
\leq C\|H_0\|_{B^{3-1/p}_{q,p}(\BR^N)},
\end{equation}
where $1 < p, q < \infty$.  On the other hand, let $\tilde \bu_{0\pm}$ be the 
extension of  $\bu_{0\pm}$ to $\BR^N_{\mp}$ such that 
$\tilde\bu_{0\pm} = \bu_{0\pm}$ on $\BR^N_\pm$ and 
\begin{equation}\label{7.2}
\|\tilde\bu_{0\pm}\|_{B^{2(1-1/p)}_{q,p}(\BR^N)}
\leq C\|\bu_{0\pm}\|_{B^{2(1-1/p)}_{q,p}(\BR^N_\pm)}.
\end{equation}
Setting 
$$\bv_\pm = e^{-(1-\Delta)t}\tilde\bu_{0\pm} = 
\CF^{-1}_\xi[e^{-(1+|\xi|^2)t}\CF[\tilde\bu_{0\pm}](\xi)],
$$
obviously we see that $\bv_\pm|_{t=0} = \bu_{0\pm}$ in $\BR^N_\pm$,
and moreover 
by the same real interpolation argument in the proof of Theorem 3.9
in \cite{SS2}, we have 
\begin{equation}\label{7.3}
\|\pd_t\bv_\pm\|_{L_p((0, \infty) L_q(\BR^N))}
+ \|\bv_\pm\|_{L_p((0, \infty), W^2_q(\BR^N))} 
\leq C\|\bu_{0\pm}\|_{B^{2(1-1/p)}_{q,p}(\BR^N_\pm)}.
\end{equation} 
Setting $\bu_\pm= \bv_\pm +  \bw_\pm$ and $H = I + J$ in
\eqref{stokes:1} and \eqref{stokes:2*}, we have 
\begin{alignat}2
&\rho_{*+}\pd_t\bw_+ - \DV\bS_{*+}(\bw_+) = \tilde\bff_+
&\quad&\text{in $\BR^N_+\times(0, T)$}\nonumber \\
&\rho_{*-}\pd_t\bw_- - \DV\bS_{*-}(\bw_-)+\nabla p_- = \tilde\bff_-,
\quad\dv \bw_- = \tilde f_\dv = \dv\tilde\bff_\dv
&\quad&\text{in $\BR^N_-\times(0, T)$} \nonumber \\
&\mu_{*-}D_{iN}(\bw_-)|_- - \mu_{*+}D_{iN}(\bw_+)|_+ = \tilde g_i
\quad(i=1, \ldots, N-1)
&\quad&\text{in $\BR^N_0\times(0, T)$}, \nonumber \\
&(\mu_{*-}D_{NN}(\bw_-) - p_-)|_-=\frac{\rho_{*-}}{\rho_{*-}-\rho_{*+}}
\Delta'J&\quad&\text{in $\BR^N_0\times(0, T)$}, \nonumber \\
&(\mu_{*+}D_{NN}(\bw_+)+(\nu_{*+}-\mu_{*+})
\dv\bw_+\bI)|_+ = \frac{\rho_{*-}}{\rho_{*-}-\rho_{*+}}
\Delta'J + \tilde g_{N+1}&\quad&\text{in $\BR^N_0\times(0, T)$}, \nonumber \\
&w_{-i}|_- - w_{+i}|_+ = \tilde h_i\quad(i=1, \ldots, N-1)
&\quad&\text{in $\BR^N_0\times(0, T)$}, \nonumber\\
&\pd_tJ - \Bigl(\frac{\rho_{*-}}{\rho_{*-}-\rho_{*+}}w_{-N}|_-
-\frac{\rho_{*+}}{\rho_{*-} - \rho_{*+}}w_{+N}|_+
\Bigr) = \tilde d&\quad&\text{in $\BR^N_0\times(0, T)$}, \nonumber \\
&\bw_\pm|_{t=0} = 0 \quad\text{in $\BR^N_\pm$},\quad
J|_{t=0} = 0 \quad\text{in $\BR^N$} \label{stokes:4}
\end{alignat}
with
\begin{align*}
\tilde\bff_+ & = \bff_+ - (\rho_{*+}\pd_t\bv_+ - \DV\bS_{*+}(\bv_+)), \quad 
\tilde\bff_- = \bff_- - (\rho_{*-}\pd_t\bv_- - 
\DV\bS_{*-}(\bv_-)) + \nabla\tilde g_N, \\
p_- & = \pi_- + \tilde g_N, \quad
\tilde g_N = \frac{\rho_{*-}}{\rho_{*-}-\rho_{*+}}
(\sigma\Delta'H_1 + g_N- \rho_{*+}g_{N+1}) - 
\mu_{*-}D_{NN}(\bv_-), \\
\tilde f_\dv & = f_{\dv} - \dv \bv_-, \quad 
\tilde\bff_\dv = \bff_\dv - \bv_-, \\
\tilde g_i & = g_- -(\mu_{*-}D_{iN}(\bv_-)|_- -
\mu_{*+}D_{iN}(\bv_+)|_+) \quad(i = 1, \ldots, N-1), \\
\tilde g_{N+1} & = \frac{\rho_{*+}}{\rho_{*-}-\rho_{*+}}
(\sigma\Delta'H + g_N -\rho_{*-}g_{N+1})
-(\mu_{*+}D_{NN}(\bv_+) + (\nu_{*+}-\mu_{*+})\dv\bv_+)|_+),\\
\tilde H_- & = h_- - (v_{-i}|_- - v_{+i}|_+)
\quad(i=1, \ldots, N-1), \\
\tilde d & = d + \Bigl(\frac{\rho_{*-}}{\rho_{*-}-\rho_{*+}}v_{-N}|_-
- \frac{\rho_{*+}}{\rho_{*-}-\rho_{*+}}v_{+N}|_+\Bigr).
\end{align*}
Setting 
\begin{align*}
&\tilde\BF_{p,q}(t) 
= \|\tilde\bff_+\|_{L_p((0, t), L_q(\BR^N_+))} + 
\|\tilde\bff_-\|_{L_p((0, t), L_q(\BR^N_-))}
+ \|\tilde f_\dv\|_{L_p((0, t), W^1_q(\BR^N_-))}
+ \|\tilde f_\dv\|_{L_p((0, t), W^{-1}_q(\BR^N_-)}\\
&\quad +\|\pd_t\tilde\bff_\dv\|_{L_p((0, T), L_q(\BR^N_-))} 
+ \sum_{i=1}^{N+1}(\|\tilde g_i\|_{L_p((0, t), W^1_q(\BR^N))}
+ \|\pd_t\tilde g_i\|_{L_p((0, t), W^{-1}_q(\BR^N))})\\
&\quad
+ \sum_{j=1}^{N-1} (\|\tilde h_j\|_{L_p((0, t), W^2_q(\BR^N))}
+ \|\pd_t\tilde h_j\|_{L_p((0, t), L_q(\BR^N))})
+ \|\tilde d\|_{L_p((0, t), W^2_q(\BR^N))}\},
\end{align*}
by \eqref{7.1} and \eqref{7.2} we have
\begin{equation}\label{est:7.1}
\tilde\BF_{p,q}(t) \leq C\BF_{p,q}(t)
\quad\text{for any $t \in (0, T)$}.
\end{equation}
By \eqref{comp:1}
\begin{equation}\label{comp:2}
\tilde f_\dv|_{t=0}=0, \quad
\tilde\bff_\dv|_{t=0} = 0, \quad
\tilde g_i|_{t=0} = 0, \quad
\tilde h_j|_{t=0} = 0
\end{equation}
for $i = 1, \ldots, N-1$ and $N+1$, and $j = 1, \ldots, N-1$. 
Given function $f$ defined on $(0, T)$, the extension operator 
$E_tf$ is defined by
\begin{equation}\label{ext:7.1}
[E_tf](\cdot, s) = \begin{cases}
f_0(\cdot, s) &\quad -\infty < s < t, \\
f_0(\cdot, 2t-s) &\quad t < s < \infty,
\end{cases}
\end{equation}
where $f_0$ is the zero extension of $f$ to $(-\infty, 0)$, that is
$f_0(\cdot, s) = f(\cdot, s)$ for $0 < s < t$ and $f_0(\cdot, s) = 0$
for $-\infty < s < 0$.  Obviously, $[E_tf](\cdot, s) = 0$ for 
$s \not\in [0, 2t]$. Moreover, if $f|_{t=0} = 0$, then 
\begin{equation}\label{ext:7.2}
\pd_s[E_tf](\cdot, s) = 
\begin{cases} 0 &\quad -\infty < s < 0, \\
\pd_s f(\cdot, s) &\quad 0 < s < t, \\
-(\pd_sf)(\cdot, 2t-s)  &\quad t < s < 2t, \\
0 &\quad 2t < s < \infty.
\end{cases}
\end{equation}
Instead of \eqref{stokes:4}, we consider the equations:
\begin{alignat}2
&\rho_{*+}\pd_t\bw_+ - \DV\bS_{*+}(\bw_+) = \tilde\bff_+
&\quad&\text{in $\BR^N_+\times\BR$}\nonumber \\
&\rho_{*-}\pd_t\bw_- - \DV\bS_{*-}(\bw_-)+\nabla p_- = \tilde\bff_-,
\quad\dv \bw_- = \tilde f_\dv = \dv\tilde\bff_\dv
&\quad&\text{in $\BR^N_-\times\BR$} \nonumber \\
&\mu_{*-}D_{iN}(\bw_-)|_- - \mu_{*+}D_{iN}(\bw_+)|_+ = \tilde g_i
\quad(i=1, \ldots, N-1)
&\quad&\text{in $\BR^N_0\times\BR$}, \nonumber \\
&(\mu_{*-}D_{NN}(\bw_-) - p_-)|_-=\frac{\rho_{*-}}{\rho_{*-}-\rho_{*+}}
\Delta'J&\quad&\text{in $\BR^N_0\times\BR$}, \nonumber \\
&(\mu_{*+}D_{NN}(\bw_+)+(\nu_{*+}-\mu_{*+})
\dv\bw_+\bI)|_+ = \frac{\rho_{*-}}{\rho_{*-}-\rho_{*+}}
\Delta'J + \tilde g_{N+1}&\quad&\text{in $\BR^N_0\times\BR$}, \nonumber \\
&w_{-i}|_- - w_{+i}|_+ = \tilde h_i\quad(i=1, \ldots, N-1)
&\quad&\text{in $\BR^N_0\times\BR$}, \nonumber\\
&\pd_tJ - \Bigl(\frac{\rho_{*-}}{\rho_{*-}-\rho_{*+}}w_{-N}|_-
-\frac{\rho_{*+}}{\rho_{*-} - \rho_{*+}}w_{+N}|_+
\Bigr) = \tilde d&\quad&\text{in $\BR^N_0\times\BR$}.\label{stokes:5}
\end{alignat}
Let $\CL$ and $\CL^{-1}$ be the Laplace transform with respect to $t$
and its inverse transform defined by 
$$\CL[f](\lambda) = \int^\infty_{-\infty}e^{-\lambda t}
f(\cdot, t)\,dt,\quad
\CL^{-1}[g(\cdot,\lambda)](t) = \frac{1}{2\pi}
\int^\infty_{-\infty}e^{\lambda t}g(\gamma + i\tau)\,d\tau
$$
with $\lambda = \gamma + i\tau \in \BC$. If we denote the Fourier
transform with respect to $t$ and its inverse  transform by $\CF_t$ and 
$\CF^{-1}_\tau$, then we have
$\CL[f](\cdot,\lambda) = \CF_t[e^{-\gamma t}f](\tau)$ and 
$\CL^{-1}[g(\cdot, \lambda)](t) = e^{\gamma t}\CF^{-1}_\tau
[g(\cdot, \gamma + i\tau)](t)$. Let $\CA_\pm(\lambda)$, $\CP_-(\lambda)$
and $\CH(\lambda)$ be $\CR$-bounded solution operators of 
problem \eqref{eq:2} given in Theorem \ref{main:r-bound}.  If we 
apply the Laplace transform with respect to $t$
to \eqref{stokes:5}, then we have
the generalized resolvent equation \eqref{eq:2} with 
\begin{gather*}
\bff_\pm = \CL[E_t\tilde\bff_\pm], \enskip 
\tilde f_- = \CL[E_t\tilde f_\dv], \enskip \tilde\bff_- = 
\CL[E_t\tilde\bff_\dv], \\
g_m = \CL[E_t\tilde g_m], \enskip g_N=0, \enskip 
g_{N+1} = \CL[E_t\tilde g_{N+1}], \enskip 
h_m = \CL[E_t\tilde h_m].
\end{gather*}
Thus, setting
$$\bw_\pm = \CL^{-1}[\CA_\pm(\lambda)\bG_\lambda], \quad
p_- = \CL^{-1}[\CP_-(\lambda)\bG_\lambda], \quad
J= \CL^{-1}[\CH(\lambda)\bG_\lambda)]$$
with 
\begin{align*}
\bG_\lambda = &(\CL[E_t\tilde\bff_+], \CL[E_t\tilde\bff_-],
\lambda^{1/2}\CL[E_t\tilde f_\dv], \nabla \CL[E_t\tilde f_\dv], 
\lambda\CL[E_t\tilde\bff_\dv], \\
&\lambda^{1/2}\CL[E_t\tilde\bg], \nabla\CL[E_t\tilde\bg],
\lambda\CL[E_t\tilde\bh], \lambda^{1/2}\nabla\CL[E_t\tilde\bh],
\nabla^2\CL[E_t\tilde\bh]),
\end{align*}
where $\tilde\bg = (\tilde g_1, \ldots, \tilde g_{N-1}, 0, \tilde g_N)$
and $\tilde\bh = (\tilde h_1, \ldots, \tilde h_{N-1})$, we see that 
$\bw_\pm$, $p_-$ and $J$ satisfy the equations \eqref{stokes:5} and 
the estimate
\begin{align}
&\CE_\gamma \leq C\{
\|e^{-\gamma t}E_t\bff_+\|_{L_p(\BR, L_q(\BR^N_+))}
+ \|e^{-\gamma t}(E_t\bff_-,
\Lambda^{1/2}_\gamma E_t\tilde f_\dv,
\nabla  E_t\tilde f_\dv, \pd_tE_t\tilde\bff_\dv)\|_{L_p(\BR, L_q(\BR^N_-))}
\nonumber\\
&\quad + \|e^{-\gamma t}
(\Lambda^{1/2}_\gamma E_t\tilde \bg,
\nabla  E_t\tilde \bg, \pd_t E_t\tilde \bh,
\Lambda_\gamma^{1/2}\nabla  E_t\tilde \bh,
\nabla^2 E_t\tilde \bh)
\|_{L_p(\BR, L_q(\BR^N))}\}.
\label{est:7.2}
\end{align}
for any $\gamma \geq \gamma_0$ with some constants $C$ and $\gamma_0$ 
with the help of Weis's operator valued Fourier multiplier theorem
\cite{Weis}, where we have set 
$\Lambda_\gamma^{1/2}
= \CL^{-1}[\lambda^{1/2}\CL[f]]$ and 
\begin{align*}
\CE_\gamma &=
\sum_{\ell=\pm} \|e^{-\gamma t}(\pd_t\bw_\ell, 
\Lambda_\gamma^{1/2}\nabla\bw_\ell, \nabla^2\bw_\ell)
\|_{L_p(\BR, L_q(\BR^N_\ell)} + \|e^{-\gamma t}\nabla p_-
\|_{L_p(\BR, L_q(\BR^N_-))} \\
&\quad  + 
\|e^{-\gamma t}(\pd_t J, \nabla J)\|_{L_p(\BR, W^2_q(\BR^N))},
\end{align*}
and we have used the fact that
$$(\pd_t, \Lambda^{1/2}_\gamma\nabla, \nabla^2)\bw_\pm
= \CL^{-1}[G^1_\lambda\CA_\pm(\lambda)\bG_\lambda],
\quad \nabla p_-= \CL^{-1}[\nabla \CP_-(\lambda)\bG_\lambda],
\quad (\pd_t, \nabla)J = \CL^{-1}[G^2_\lambda \CH(\lambda)\bG_\lambda].$$ 
As was seen in Shibata and Shimizu \cite{SS2}, we know that 
\begin{align}
&\gamma\|e^{-\gamma t}f\|_{L_p(\BR, L_q(\Omega)}  \leq 
C\|e^{-\gamma t}\pd_tf\|_{L_p(\BR, L_q(\Omega))}, \nonumber\\
&\|\Lambda_\gamma^{1/2}f\|_{L_p(\BR, L_q(\Omega))}
 \leq C\{\|\pd_tf\|_{L_p(\BR, W^{-1}_q(\Omega))} 
+ \|f\|_{L_p(\BR, W^1_q(\Omega))}, \\
&\|e^{-\gamma t}\pd_t f\|_{L_p(\BR, L_q(\Omega))}
+ \|e^{-\gamma t}f\|_{L_p(\BR, W^2_q(\Omega))}
\leq C\|e^{-\gamma t}(\pd_tf, \Lambda^{1/2}_\gamma\nabla f, 
\nabla^2f)\|_{L_p(\BR, L_q(\Omega))} \label{est:7.3}
\end{align}
for any $\gamma \geq \gamma_0$ with some constant $C>0$, where
$\Omega = \BR^N_\pm$ or $=\BR^N$.  Thus, combining \eqref{est:7.2}
and \eqref{est:7.3}, we have
\begin{equation}
\CE_\gamma + 
\sum_{\ell=\pm}\gamma\|e^{-\gamma t}\bw_\pm\|_{L_p(\BR, L_q(\BR^N_\ell))}
+ \gamma\|e^{-\gamma t}J\|_{L_p(\BR, W^2_q(\BR^N))}
\leq C\CM_\gamma
\label{est:7.4}
\end{equation}
for any $\gamma \geq \gamma_0$ with some constant $C > 0$ 
with 
\begin{align*}
\CM_\gamma & = \sum_{\ell=\pm}\|e^{-\gamma t}E_t\tilde\bff_\ell
\|_{L_p(\BR, L_q(\BR^N_\ell))} 
+ \|e^{-\gamma t}\pd_t E_t\tilde f_\dv\|_{L_p(\BR, W^{-1}_q(\BR^N_-))}
+ \|e^{-\gamma t}E_t\tilde f_\dv\|_{L_p(\BR, W^1_q(\BR^N_-))} \\
&\quad + \|e^{-\gamma t}\pd_tE_t\tilde\bff_\dv\|_{L_p(\BR, L_q(\BR^N_-))}
+ \|e^{-\gamma t}\pd_tE_t\tilde\bg\|_{L_p(\BR, W^{-1}_q(\BR^N))}
+ \|e^{-\gamma t}E_t\tilde\bg\|_{L_p(\BR^N, W^{-1}_q(\BR^N))}\\
&\quad + \|e^{-\gamma t}\pd_tE_t\tilde\bh\|_{L_p(\BR, L_q(\BR^N))}
+ \|e^{-\gamma t}E_t\tilde\bh\|_{L_p(\BR, W^2_q(\BR^N))}
\end{align*}
for any $\gamma \geq \gamma_0$ with some constant $C > 0$.
Especially, it follows from \eqref{est:7.4} and 
the fact that $E_tf$ vanishes for $t < 0$ that 
$\bv\pm$ and $J$ also vanish for $t < 0$.  In fact, we observe that
\begin{align*}
&\gamma(\|\bw_+\|_{L_p((-\infty, 0), L_q(\BR^N_+))}
+ \|\bw_-\|_{L_p((-\infty, 0), L_q(\BR^N_-))}
+ \|J\|_{L_p((-\infty, 0), W^2_q(\BR^N))} \\
&\leq \gamma(
\|e^{-\gamma t}\bv_+\|_{L_p(\BR, L_q(\BR^N_+))}
+ \|e^{-\gamma t}\bv_-\|_{L_p(\BR, L_q(\BR^N_-))}
+ \|e^{-\gamma t}J\|_{L_p(\BR, W^2_q(\BR^N))}
\leq C\CM_\gamma \leq CM_{\gamma_0}
\end{align*}
for any $\gamma \geq \gamma_0$, so that letting $\gamma \to \infty$
we have $\|\bw_\pm\|_{L_p((-\infty, 0), L_q(\BR^N\pm))}=0$
and $\|J\|_{L_p((-\infty, 0), W^2_q(\BR^N))} = 0$, which implies that
$\bw_\pm$ and $J$ vanishe for $t < 0$.  By the equation, we also have
that $\nabla p_-$ vanishes for $t < 0$. Summing up, we have seen that
$\bw_\pm \in L_p(0, \infty), W^2_q(\BR^N_\pm)) \cap 
W^1_p((0, \infty), L_q(\BR^N_\pm))$, $p_- 
\in L_p((0, \infty), \hat W^1_q(\BR^N_-))$ and 
$J \in L_p((0, \infty), W^3_q(\BR^N)) \cap 
W^1_p((0, \infty), W^2_q(\BR^N))$, $\bw_\pm$, $p_-$ and 
$J$ satisfy \eqref{stokes:4} with $T =t$, because $E_t f= f$ for 
$t \in (0, t)$.  Moreover, by 
\eqref{ext:7.1}, \eqref{ext:7.2}, and \eqref{est:7.4}, we have
\begin{align}
%&\sum_{\ell=\pm}\{\|\pd_t\bw_\ell\|_{L_p((0, t), L_q(\BR^N_\ell))}
%+ \|\bw_\ell\|_{L_p((0, t), W^2_q(\BR^N_\ell))}\}
%+ \|\nabla p_-\|_{L_p((0, t), L_q(\BR^N_-))}\nonumber\\
%&\quad+ \|\pd_t J\|_{L_p((0, t), W^2_q(\BR^N))}
%+ \|J\|_{L_p((0, t), W^3_q(\BR^N))}
\CI_{p,q}(\bw_\pm, p_-,J)(t)
\leq Ce^{\gamma t}\tilde \BF_{p,q}(t) \leq Ce^{\gamma t}\BF_{p,q}(t).
\label{est:7.5}
\end{align}
Since the uniqueness follows from the existence of solutions to the dual 
problem, the equation \eqref{stokes:4} admits unique soltuions for
$(0, t)$.  Especially, we can construct unique solutions of problem 
\eqref{stokes:4} for $(0, T)$ 
and these solutions are also solutions of problem 
\eqref{stokes:4} with $T =t$  for any $0 < t < T$,
so that by \eqref{est:7.5} we have completed the proof of 
Theorem \ref{thm:stokes}. 

\end{document}